\documentclass[opre,sglanonrev]{informs3_noinforms2c}
\RequirePackage{tgtermes}
\RequirePackage{newtxtext}
\RequirePackage{newtxmath}
\RequirePackage{bm}
\RequirePackage{endnotes}
\usepackage{hyperref}
\usepackage[linesnumbered,ruled,vlined]{algorithm2e}
\SetKwInput{KwParam}{Hyperparameters}
\SetKwInput{KwSub}{Subroutines}
\SetKwInput{KwGlobal}{Global}
\SetKwInput{KwPredef}{Predefined functions}

\OneAndAHalfSpacedXII 

\usepackage{algorithm2e}
\usepackage{tikz}

\usepackage{natbib}
 \bibpunct[, ]{(}{)}{,}{a}{}{,}%
 \def\bibfont{\small}%
\newtheorem{observation}{Observation}

\EquationsNumberedThrough    

\TheoremsNumberedThrough     
\ECRepeatTheorems  %

\MANUSCRIPTNO{MOOR-0001-2024.00}

\begin{document}


\RUNAUTHOR{Cetin, Chen, and Goldberg}

\RUNTITLE{Online Stochastic Packing with General Correlations}

\TITLE{Online Stochastic Packing and Combinatorial Allocation with General Correlations}

\ARTICLEAUTHORS{%
\AUTHOR{Sabri Cetin}
\AFF{Cornell ORIE,\EMAIL{sc2955@cornell.edu}}

\AUTHOR{Yilun Chen}
\AFF{CUHK Shenzhen,\EMAIL{chenyilun@cuhk.edu.cn}}

\AUTHOR{David A. Goldberg}
\AFF{Cornell ORIE,\EMAIL{dag369@cornell.edu}}

} 

\ABSTRACT{
There has been a growing interest in studying online stochastic packing and combinatorial allocation under more general correlation structures, motivated by the complex data sets and models driving modern applications.  Several past works either assume correlations are weak or have a particular structure, have a complexity scaling with the number of Markovian “states of the world" (which may be exponentially large in the case of full history dependence), scale poorly with the horizon $T$, or make additional continuity assumptions.  We show that for all $\epsilon$, the online stochastic packing (combinatorial allocation) linear programming problem with general correlations (suitably normalized and with sparse columns) has an approximately optimal policy (with optimality gap $\epsilon T$) whose per-decision runtime scales as the time to simulate a single sample path of the underlying stochastic process (assuming access to a Monte Carlo simulator), multiplied by a constant independent of the horizon or number of Markovian states.  We derive analogous results for network revenue management (also with flexible products), and online bipartite matching and independent set in bounded-degree graphs, by rounding.  Our algorithms implement stochastic gradient methods in a novel on-the-fly/recursive manner for the associated massive deterministic-equivalent linear program.
}




 

\maketitle


\section{Introduction}\label{introsec1}
We consider the problem of online stochastic packing (o.s.p.) and combinatorial allocation (c.a.) under general correlations, in which at each time period the decision maker (DM) must select one of $q$ options, each with a different resource requirement and reward, while satisfying all resource constraints.  In the simplest case $q = 2$, the two options may correspond to either accepting the arriving item (earning the corresponding reward and using the required resources) or rejecting the item.  In more complex cases of c.a., one option may correspond to rejecting the arriving item, with each other option corresponding to a different way to satisfy that demand.  We assume that option 0 always results in zero reward and zero resources being utilized.  More precisely, we assume that there are $m$ resources $\lbrace \textrm{RES}_i, i = 1,\ldots,m \rbrace$, each with a known budget $\lbrace b_i, i = 1,\ldots,m \rbrace$, and a known time horizon of $T$ periods.  In each period $t$ the DM observes new information, including the rewards $\lbrace Z^r_t , r = 0,\ldots,q-1 \rbrace$ and resource consumption vectors (r.c.v.s) $\lbrace \overline{a}^r_t, r = 0,\ldots,q-1 \rbrace$ (with $a^r_{i,t}$ the non-negative amount of $\textrm{RES}_i$ utilized if the DM selects option $r$ at time $t$).  
\\\indent It may be that the DM also observes additional information (e.g. updated market indicators), or it may be that the only new information consists of the new rewards and r.c.v.s for the $q$ options at time $t$.  We abstract this concept to assuming the DM observes a stochastically evolving (in general non-Markovian and high-dimensional) \textbf{observable information process} (o.i.p.), and that the realized rewards and r.c.v.s of the $q$ options for all arrivals through time $t$ are part of this o.i.p. (through time $t$).  In many settings of interest this o.i.p. will simply be the set of realized rewards and r.c.v.s. (which may have a general joint distribution).  Thus in each period $t$ the o.i.p. is updated, and the DM must then decide which option they will take (i.e. set $X^r_t = 1$ for exactly one $r \in \lbrace 0,\ldots,q-1 \rbrace$, and set $X^{r'}_t$ equal 0 for all other $r'$).  This repeats for the $T$ periods of the horizon.  
\\\indent The DM's goal is to maximize $\mathbb{E}\big[ \sum_{t=1}^T \sum_{r=1}^{q-1} Z^r_t X^r_t \big]$ subject to $\sum_{t=1}^T \sum_{r=1}^{q-1} a^r_{i,t} X^r_t \leq b_i$ w.p.1 for $i = 1,\ldots,m.$  The maximization is over all policies, where a policy is equivalent to the selection of a stochastic process $\lbrace X^r_t, r = 0,\ldots,q-1, t = 1,\ldots,T \rbrace$ adapted to the filtration generated by the o.i.p.  When the o.i.p. is simply the set of realized rewards and r.c.v.s., this is equivalent to the decision at time $t$ being a function of the realized rewards and r.c.v.s through time $t$.  Such high-dimensional o.s.p./c.a. problems find an array of applications across Operations Research (OR) and Computer Science (CS), and may exhibit complex structure with long horizon and non-Markovian dynamics (which can lead to severe algorithmic challenges), as in the following examples.
\begin{itemize}
    \item \textbf{Network Revenue Management (NRM).\ }  In the $\textrm{NRM}$ problem, each item corresponds to a customer who desires some product, and the DM must decide whether to make the sale (using different amounts of resources $i=1,\ldots,m$) or not.  While the literature has typically made strong independence assumptions on the customer requests (and associated profits), real-world $\textrm{NRM}$ problems may exhibit a more complex dependency structure (\cite{gallego2019revenue}).  Several works have studied $\textrm{NRM}$ with more general dependencies, but either assume a particular model of dependency, or exhibit a computational complexity scaling with the number of Markovian states of the world (which may be exponential in $T$ in the history-dependent setting), see\ \cite{Demiguel2006,jiang2023,lietal2025,aouadetal2022,baietal2023}.
    \item \textbf{Online Matching.\ } In the online matching problem, each item corresponds to an edge which is revealed (sequentially) in an unknown random graph.  In each time period the DM must decide whether to accept the edge or not, subject to the constraint that at the end of the horizon the set of accepted edges constitutes a feasible matching in the graph (i.e. each node is incident to at most one edge).  Online matching has applications to a broad range of problems, including advertising, transportation, and healthcare (\cite{huangetal2024}).  Recently there has been a recognition that real-world applications may require stochastic models that go beyond the independence assumptions made in much of the literature (\cite{aouadetal2022,gaoetal2025,feldmanetal2025}).
\item \textbf{Combinatorial Allocation.\ }  In the c.a. problem, each item corresponds to a user request whose demand can be met by different subsets of resources (with each subset accruing a potentially different welfare for the user), subject to an overall budget for each resource.  Online combinatorial allocation has applications to problems including online auctions, adwords, and dynamic virtual machine allocation (\cite{tanetal2020}), as well as NRM with flexible and/or opaque products, in which a customer's demand can be met through different resource allocations (\cite{zhuetal2024,housnietal2025}).  Recently, there has been an interest in extending the results in this space to models with more complex dependency structures (\cite{feldmanetal2025}).  
\end{itemize}
\ \indent Although o.s.p./c.a. with general correlations is a fundamental problem studied across multiple communities (as we detail in our literature review Section\ \ref{litsec}), to the best of our knowledge all algorithms to date either : (1) incur an exponential dependence on the horizon $T$ (or dimension $D$ of the o.i.p. when that process is Markovian), (2) make additional assumptions regarding the correlations and/or o.i.p., or (3) do not come with theoretical guarantees.  The main question for this work is whether there exists an algorithm which simultaneously addresses these three points.  Our main contribution is a positive resolution of this question, as we now describe in greater detail. 
\\\textbf{Motivating examples for our approach and algorithm}.  Here we motivate our approach and algorithm by considering two illustrative examples of NRM with correlations.  For simplicity, we restrict our discussion to the case $m = 1$, $q = 2$, $T \geq 9$, and $a^1_t = 1$ w.p.1 for all $t$, i.e. multiple stopping, also assuming that the o.i.p. consists only of the realized rewards.  Since when $q = 2$ only option 1 is non-trivial, here we simply refer to option 1's reward as ``the reward".  It is well-known that when the rewards are independent or have weak correlations, resolve heuristics obtain good performance (\cite{balseiro2024survey}), but when those assumptions break different methods are needed (\cite{aouadetal2022,baietal2023,jiang2023,lietal2025}).  
\\\indent\textbf{Example 1.}  We consider an example where the reward in a given period may provide a strong signal about future rewards.  Similar models (e.g. in which current demand provides a signal about future demands) have been widely studied in the OR literature (\cite{xinetal2022,lietal2025b}).  Thus suppose the reward sequence is one of two sequences : $S^1$ (w.p. $.3$) or $S^0$ (w.p. $.7$).  Sequence $S^1$ is defined as : $S^1_t = .5$ for $t = 1,\ldots,\lfloor \frac{T}{3} \rfloor - 2$; $S^1_{\lfloor \frac{T}{3} \rfloor - 1} = .01$; $S^1_t = .45$ for  $t = \lfloor \frac{T}{3} \rfloor,\ldots,T - \lfloor \frac{T}{3} \rfloor + 2$; $S^1_t = 1$ for $t = T - \lfloor \frac{T}{3} \rfloor + 3,\ldots,T$.  Analogously, $S^0_t = .5$ for $t = 1,\ldots,\lfloor \frac{T}{3} \rfloor - 2$; $S^0_{\lfloor \frac{T}{3} \rfloor - 1} = .001$; $S^0_t = .45$ for  $t = \lfloor \frac{T}{3} \rfloor,\ldots,T - \lfloor \frac{T}{3} \rfloor + 2$; $S^0_t = 0$ for $t = T - \lfloor \frac{T}{3} \rfloor + 3,\ldots,T$.  $S^1$ and $S^0$ are indistinguishable until time $\lfloor \frac{T}{3} \rfloor - 1$, when it becomes known (based on whether that reward is $.01$ or $.001$) whether the sequence is $S^1$ (implying the final $\lfloor \frac{T}{3} \rfloor - 2$ rewards are 1) or $S^0$ (implying the final $\lfloor \frac{T}{3} \rfloor - 2$ rewards are 0).  We assume the budget $b = \lfloor \frac{T}{3} \rfloor - 2$, so in principle one could saturate the budget by either taking all items with reward $.5$, $b$ items with reward $.45$, or the final $b$ items (with rewards 0 or 1).  We denote the expected reward earned by an optimal policy by $\textrm{OPT}_T$.  
\\\indent First, let us see how three well-studied heuristics peform.  The Certainty-Equivalent (C-E) heuristic (denoted $CE$) solves (in each period) the associated C-E LP (which relaxes the problem to one involving the expected number of items of each reward value) for the to-go problem, and determines the current action by randomly rounding (\cite{jasinetal2012}).  The Fluid Bayes Selector heuristic (denoted $FBayes$) of \cite{veraetal2021} works similarly, but rounds using a threshold (\cite{veraetal2021}).  The (non-fluid) Bayes Selector heuristic (denoted $Bayes$) of \cite{veraetal2021} selects an action to minimize the probability of disagreement with the ``full information/hindsight" relaxation (h.r.).  For a policy $\pi$, let $\textrm{VAL}_{\pi,T}$ denote the expected reward earned by policy $\pi$ for Example 1.  In Appendix\ \ref{heuristicformalsec}, we prove that $\liminf_{T \rightarrow \infty} T^{-1} \big( \textrm{OPT}_T - \textrm{VAL}_{\pi,T} \big) \geq .007$ for $\pi \in \lbrace CE,FBayes,Bayes \rbrace$.  In contrast to the setting in which the reward sequence has weak correlations (where resolve heuristics match the performance even of the h.r., which is larger than the optimal value, up to error $o(T)$, see \cite{balseiro2024survey}), here these heuristics do not even match the optimal policy on the fluid scale (i.e. up to first-order scaling in $T$).  As we show in Appendix\ \ref{heuristicformalsec}, these heuristics underestimate the value of deferring decisions until more information is revealed, incorrectly accepting too many early rewards.
\\\indent Next, let us comment on the value of the h.r. (commonly used as a benchmark to compute regret), i.e. the best a DM can do in expectation if the (random) reward sequence is revealed at time 1, and denote this value by $\textrm{HR}_T$.  It is easily verified that $HR_T = .7 \times .5 \times (\lfloor \frac{T}{3} \rfloor - 2) + .3 \times 1 \times (\lfloor \frac{T}{3} \rfloor - 2)$.  It follows from our analysis in Appendix\ \ref{heuristicformalsec} that $\lim_{T \rightarrow \infty} T^{-1} \big( \textrm{HR}_T - \textrm{OPT}_T \big) = \frac{7}{600}$.  In contrast to the setting of weak correlations, here the optimal value and h.r. do not even match on the fluid scale, since the large fluctuations imply a greater benefit to the DM in learning the sequence upfront.
\\\indent Finally, let us compute the optimal policy.  As policies must be adapted, at time $t$ the DM must specify what action to take as a function of the partial sequence of rewards realized in periods $1,\ldots,t$.  Let ${\mathcal S}^t$ denote the set of all possible partial sequences of realized rewards in periods $1,\ldots,t$.  For $i = 0,1$, let $S^i_{[t]} = (S^i_1,\ldots,S^i_t)$ denote the partial sequence in the first $t$ periods.  Then ${\mathcal S}^t$ consists of the single partial sequence of exactly $t$ .5's for $t = 1,\ldots,\lfloor \frac{T}{3} \rfloor - 2$, and consists of the two partial sequences $S^1_{[t]}$ and $S^0_{[t]}$ otherwise.  Specifying a policy is equivalent to specifying a binary action $X(E)$ on every partial sequence $E \in {\mathcal E} \stackrel{\Delta}{=} \bigcup_{t=1}^T {\mathcal S}^t$.  There are $\lfloor \frac{T}{3} \rfloor - 2 + 2 \big( T - (\lfloor \frac{T}{3} \rfloor - 2) \big) = 2 T - \lfloor \frac{T}{3} \rfloor + 2$ such $X(E)$ variables in total, which scales as $\frac{5}{3} T$ for large $T$.  
\\\indent We now formalize the problem of computing an optimal policy (i.e. solving for $\lbrace X(E), E \in {\mathcal E} \rbrace$) as an explicit optimization, first reasoning about the objective function.  For $E \in {\mathcal E}$, let $\mu(E)$ denote the probability that partial reward sequence $E$ is realized.  Note that $\mu(E) = 1$ if $E$ is a sequence of $t$ $.5$'s for $t = 1,\ldots,\lfloor \frac{T}{3} \rfloor - 2$, while $\mu(S^1_{[t]}) = .3$ and $\mu(S^0_{[t]}) = .7$ for $t = \lfloor \frac{T}{3} \rfloor - 1,\ldots,T$.  For $t = 1,\ldots,T$ and $E \in {\mathcal S}^t$, let $Z(E)$ denote the reward realized at time $t$ on partial reward sequence $E$ (i.e. the final reward of $E$).  Thus $Z(E) = .5$ if $E$ is a sequence of $t$ $.5$'s for $t = 1,\ldots,\lfloor \frac{T}{3} \rfloor - 2$, and $Z(E)$ equals $S^1_t$ if $E = S_{[t]}$ and equals $S^0_t$ if $E = S^0_{[t]}$ for $t = \lfloor \frac{T}{3} \rfloor - 1,\ldots,T.$  If partial reward sequence $E$ actually appears as part of the realized sequence of rewards, the DM will earn reward $Z(E) X(E)$ (attributable to $E$).  Thus the expected total reward earned by the DM equals $\sum_{E \in {\mathcal E}} \mu(E) Z(E) X(E)$.  The feasibility constraints require $\sum_{t=1}^T X(S^1_{[t]}) \leq b$ and $\sum_{t=1}^T X(S^0_{[t]}) \leq b,$ where to be clear $S^i_{[t]}$ is in all cases equivalent to some $E \in {\mathcal E}.$  The above defines an integer program (IP) equivalent to the problem of computing an optimal policy.  In fact, it follows from the results of \cite{naoretal2025} (for uniform matroids) that the natural LP relaxation always has an integer optimal solution (for multiple stopping), and it thus suffices to solve the following LP on $2 T - \lfloor \frac{T}{3} \rfloor + 2$ variables :  
\begin{equation}\label{LP000}
\max_{X} \sum_{E \in {\mathcal E}} \mu(E) Z(E) X(E)\ s.t.\ \sum_{t=1}^T X(S^1_{[t]}) \leq b, \sum_{t=1}^T X(S^0_{[t]}) \leq b, X(E) \in [0,1]\ \forall\ E \in {\mathcal E}.
\end{equation}
\ \indent Although we defer the proof to Appendix\ \ref{heuristicformalsec}, it can be shown that an optimal solution to (\ref{LP000}) corresponds to the DM skipping all the initial rewards of value .5, and then : (1) if the reward at time $\lfloor \frac{T}{3} \rfloor - 1$ is .001 the DM uses their entire budget on rewards of value $.45$; and (2) if instead the reward at time $\lfloor \frac{T}{3} \rfloor - 1$ is .01 the DM uses their entire budget on rewards of value $1$.  Such a policy achieves an expected reward of $.3 \times (\lfloor \frac{T}{3} \rfloor - 2) \times 1 + .7 \times (\lfloor \frac{T}{3} \rfloor - 2) \times .45$, which scales as $\frac{41}{200} T$ for large $T$, which (on the fluid scale) is larger than that earned by the analyzed heuristics, and smaller than that earned by the (unachievable) h.r.  Intuitively, the optimal policy is able to use the fact that at time $\lfloor \frac{T}{3} \rfloor - 1$ the DM learns whether the remaining rewards will contain a large number of 1's or not, a fact the heuristics cannot exploit.  Although neither the resolve heuristics nor the h.r. can provide an approximation to the optimal policy accurate on even the fluid scale, an approach that does work is to simply solve LP (\ref{LP000}).  In the next example, we will see why this strategy is in general intractable, thus motivating our own approach.
\\\indent\textbf{Example 2.}  Suppose that at each time $t$, there are two values the reward can take, $l_t$ and $h_t$.  However, the probability that the reward at time $t$ takes value $l_t$ (as opposed to $h_t$) may be a complex function of the entire history of realized rewards up to that time.  Here there are $2^T$ possible reward sequences (one for each length-T sequence $S$ with $S_t$ either $l_t$ or $h_t$), and $|{\mathcal E}| = \sum_{t=1}^T 2^t = 2 \times (2^T - 1)$.  Thus specifying a policy requires assigning values to an \textbf{exponentially large} (in $T$) number of variables $\lbrace X(E), E \in {\mathcal E} \rbrace$.  The associated LP is now : 
\begin{equation}\label{pack00c}
\max_{X} \sum_{E \in {\mathcal E}} \mu(E) Z(E) X(E)\ \ s.t.\ \ \sum_{t=1}^T X(S_{[t]}) \leq b\ \forall\ S \in \mathcal{S}^T, X(E) \in [0,1]\ \forall\ E \in \mathcal{E}.
\end{equation}
It is not clear what we can do without incurring an exponential runtime, since just outputting the solution would take exponential time!  Although there is a vast literature on methods for solving large LPs (and associated multi-stage stochastic programming and control problems), to the best of our knowledge past approaches cannot overcome this issue (without additional assumptions).
\\\textbf{What problem are we solving?} It would seem impossible to efficiently solve (\ref{pack00c}).  However, the DM does not need to specify the policy for all $E \in {\mathcal E}$, and in reality need only compute a decision in real time for each partial history they actually encounter when solving the online problem, of which there will be exactly $T$ (on the realized reward sequence $S$).  Thus the DM need only compute the optimal values of very few variables in (\ref{pack00c}).  At least in principle such a task could be accomplished efficiently, and it would suffice to have an algorithm to efficiently compute an (approximately optimal) value $X(E)$ for any given $E \in {\mathcal E}$ (as the DM could then call it on each of $S_{[1]},\ldots,S_{[T]}$).  However, it is not clear whether it is possible to efficiently compute the optimal values of a few variables in\ (\ref{pack00c}) without computing many values.  We provide such an algorithm (from which our main results follow), and give its intuition below.
\\\textbf{Algorithmic Intuition.}  At a high level, our algorithm solves a penalty-based relaxation of (\ref{pack00c}) ``on-the-fly" using a recursive stochastic gradient ascent (SGA) method which only implements those calculations required, and rounds the solution (on-the-fly) to derive a policy.  In particular, for a convex penalty function $\phi$ (with derivative $\phi'$), we consider the following relaxation of (\ref{pack00c}) :
\begin{equation}\label{prob0}
\max_{X} \sum_{E \in \mathcal{E}} \mu(E) Z(E) X(E)  - \sum_{S \in \mathcal{S}^T} \mu(S) \phi\left( \sum_{t=1}^T X(S_{[t]}) - b \right)\ s.t.\ X(E) \in [0,1]\ \forall\ E\ \in \mathcal{E}.
\end{equation}
A key question is whether we can efficiently compute an (approximately optimal) value $X(E)$ for a given $E.$  In the language of SGA, can we efficiently compute ${\mathcal X}^K(E)$, with ${\mathcal X}^K(E)$ the value of $X(E)$ after $K$ iterations of SGA run on $(\ref{prob0})$ (for appropriate $K$).  For our notion of approximate optimality (i.e. $\epsilon T$ optimality gap), it will suffice for $K$ to scale as $\textrm{poly}(\frac{1}{\epsilon})$ (with a step-size $\alpha$ scaling as $\textrm{poly}(\epsilon)$), where throughout we use $\textrm{poly}(\cdot)$ to denote a generic polynomial function.  Standard definitions imply that ${\mathcal X}^K(E)$ can be expressed (recursively) as a simple function of ${\mathcal X}^k(E')$ for $k < K$ and $E' \in \mathcal{E}$.  We are thus led to the question of how many different ${\mathcal X}^k(E')$ values we actually need to calculate ${\mathcal X}^K(E)$, and (supposing the answer is much smaller than $|{\mathcal E}|$) what algorithm can efficiently implement such a calculation?  We now argue that the answer is indeed a very small number (constant depending only on $\epsilon$, independent of $|{\mathcal E}|$), and sketch our algorithm for efficiently implementing the calculation.  Let us first characterize the gradient of the objective of (\ref{prob0}) w.r.t. variables $\lbrace X(E), E \in {\mathcal E} \rbrace$.  Let $\zeta_E$ be a random reward sequence conditioned to start with partial sequence $E$.  Then the component of the gradient corresponding to $E \in \mathcal{E}$ equals $\mu(E) \left( Z(E) - \mathbb{E}_{S' \sim \zeta_E} \left[ \phi'\left( \sum_{t=1}^T X(S'_{[t]}) - b \right) \right] \right),$ where $\sim$ denotes equivalence in distribution.  We will account for the $\mu(E)$ implicitly using a weighted Euclidean norm, and thus let us ``pretend" that the gradient component corresponding to $E$ is $Z(E) - \mathbb{E}_{S' \sim \zeta_E} \left[ \phi'\left( \sum_{t=1}^T X(S'_{[t]}) - b_i \right) \right]$.  Using our simulator $\textrm{SIM}$, we can draw a single reward sequence $S'$ conditioned to start with partial sequence $E$, and form an unbiased estimate for gradient component $E$ equal to $Z(E) - \phi'\left( \sum_{t=1}^T X(S'_{[t]}) - b \right).$  Letting $\Pi(x)$ denote the projection of $x$ onto $[0,1]$ (for $x \in {\mathcal R}$), SGA sets 
$${\mathcal X}^K(E) = \Pi\Bigg( {\mathcal X}^{K-1}(E) + \alpha \times \bigg( 
Z(E) - \phi'\big( \sum_{t=1}^T {\mathcal X}^{K-1}(S'_{[t]}) - b \big) \bigg) \Bigg).$$
Thus to compute ${\mathcal X}^K(E)$ using SGA, one need only evaluate ${\mathcal X}^{K-1}$ at $\lbrace S'_{[t]}, t = 1,\ldots,T \rbrace$, and combine with a bit of computation.  \textbf{The key insight is that we may apply this logic recursively}.  Thus to compute ${\mathcal X}^{K-1}(S'_{[1]})$, we draw a sample path $S''^1 \sim \zeta_{S'_{[1]}}$, and by the same reasoning deduce that to compute ${\mathcal X}^{K-1}(S'_{[1]})$ we only need the values of ${\mathcal X}^{K-2}(S''^1_{[1]}), {\mathcal X}^{K-2}(S''^1_{[2]}),\ldots, {\mathcal X}^{K-2}(S''^1_{[T]})$ (and a bit of computation).  Applying the same logic to $S'_{[2]},\ldots,S'_{[T]}$ (and defining appropriately conditioned random sample paths $S''^2,\ldots,S''^T$), we deduce that to compute ${\mathcal X}^{K}(E)$ using SGA, we really only need the $T^2$ values ${\mathcal X}^{K-2}(S''^1_{[1]}),\ldots,{\mathcal X}^{K-2}(S''^1_{[T]}); {\mathcal X}^{K-2}(S''^2_{[1]}),\ldots,{\mathcal X}^{K-2}(S''^2_{[T]}); \ldots ; {\mathcal X}^{K-2}(S''^T_{[1]}),\ldots,{\mathcal X}^{K-2}(S''^T_{[T]})$ (and a bit of computation), which we use to compute ${\mathcal X}^{K-1}(S'_{[1]}),\ldots,{\mathcal X}^{K-1}(S'_{[T]})$, which we use to compute ${\mathcal X}^K(E)$.  By recursing this logic all the way down to ${\mathcal X}^1$ (which is easy to compute from initial conditions), we conclude that to compute ${\mathcal X}^{K}(E)$ using SGA, it suffices to compute ${\mathcal X}^k(E')$ for only $\sum_{j=1}^K T^j \leq K T^K$ different $(k,E')$ pairs (and apply a similarly bounded amount of computation).  Here $T$ appears because we exactly compute $\sum_{t= 1}^T X(S'_{[t]})$ in the gradients, and hence must recurse on $T$ terms.  We prove that SGA is sufficiently robust that we can instead sample $\textrm{poly}(\frac{1}{\epsilon})$ terms randomly from the sum $\sum_{t = 1}^T X(S'_{[t]})$ to compute a ``good enough" approximation.  Using this we show that for $K = \textrm{poly}(\frac{1}{\epsilon})$, to compute ${\mathcal X}^K(E)$, it suffices to compute ${\mathcal X}^k(E')$ for only $\big( \textrm{poly}(\frac{1}{\epsilon}) \big)^{\textrm{poly}(\frac{1}{\epsilon})}$ different $(k,S')$ pairs, and combine with a similarly bounded amount of computation.
\\\indent Inspired by the above observations, we now sketch our algorithm $\textrm{R}$ for computing ${\mathcal X}^K(E)$.  To ensure we never attempt to assign ${\mathcal X}^k(E')$ two different values, we utilize a memoization table (i.e. a matrix $\Upsilon$ with $|\mathcal{E}|$ rows and an infinite number of columns) to efficiently and consistently implement the required gradient calculations.  Here $\Upsilon(E,0)$ is initialized to 0, and all other entries are initialized to a dummy value $\emptyset.$  $\forall\ E' \in \mathcal{E}$ and $k \geq 1$, $\Upsilon(E', k)$ will store the value of ${\mathcal X}^k(E')$ (if computed).
\vspace{-.25in}
\begin{algorithm}[ht]
\caption{Sketch of Algorithm $\textrm{R}(E,k)$}
\label{alg:gd-implement00}
\DontPrintSemicolon
If $\Upsilon(E, k) = \emptyset$\\
   \ \ \ Call $\textrm{SIM}(E)$ to generate a random conditioned (on $E$) reward sequence $\hat{S}$\\ 
   \ \ \ Subsample $\textrm{poly}(\frac{1}{\epsilon})$ uniformly random time indices ${\mathcal T}$\\
   \ \ \ Recursively Call $R(\hat{S}_{[t]},k-1)$ for all $t \in {\mathcal T}$\\
   \ \ \ $\Upsilon(E,k) = \Pi\Bigg(\Upsilon(E,k-1) + \alpha \times \bigg( Z(E) - \phi'\big( T \times \frac{1}{|{\mathcal T}|} \times \sum_{t \in {\mathcal T}} \Upsilon(S'_{[t]},k-1) - b\big) \bigg) \Bigg)$
\end{algorithm} 
\vspace{-.25in}
$\textrm{R}$ implements the recursive logic outlined above.  If ${\mathcal X}^k(E)$ has not been previously computed, it simulates a conditional sample path $\hat{S}$, randomly samples $\textrm{poly}(\frac{1}{\epsilon})$ time indices (to approximate the linear sum appearing in the gradient), and implements a SGA update (using the appropriate ${\mathcal X}^{k-1}(\hat{S}_{[t]})$ values).  The most important step is the $\textrm{poly}(\frac{1}{\epsilon})$ recursive calls to $R(\hat{S}_{[t]},k-1)$ to ensure the necessary ${\mathcal X}^{k-1}(\hat{S}_{[t]})$ values are precomputed.  $\textrm{R}$'s runtime is driven by these recursive calls (as the actual calculations to implement a gradient step are otherwise trivial), leading to a runtime of $\big( \textrm{poly}(\frac{1}{\epsilon}) \big)^{\textrm{poly}(\frac{1}{\epsilon})}$ for $K = \textrm{poly}(\frac{1}{\epsilon})$.  $\textrm{R}$ can be used by the $\textrm{DM}$ to implement an $\epsilon T$-optimal multiple stopping policy as follows.  The DM selects the desired accuracy $\epsilon$ and sets $K = \textrm{poly}(\frac{1}{\epsilon})$.  At each time $t$, the DM calls $\textrm{R}(S_{[t]},K)$ (with $S_{[t]}$ the partial sequence of realized rewards up to time $t$), computes the value $K^{-1} \sum_{j=1}^K \Upsilon(S_{[t]},j)$ (in the unaccelerated setting) or $\Upsilon(S_{[t]},K)$ (in the accelerated setting), rounds this value on the fly to a binary decision, and makes sure not to stop if the budget has already been exhausted.  We prove that all the associated errors can be controlled, yielding a policy with expected suboptimality $\epsilon T$, and per-decision runtime $\big( \textrm{poly}(\frac{1}{\epsilon}) \big)^{\textrm{poly}(\frac{1}{\epsilon})}$ (multiplied by the cost ${\mathcal C}$ of using $\textrm{SIM}$).  Although our policy is not as interpretable as the heuristics analyzed in Example 1, its performance (up to error $\epsilon T$) matches that of the complex high-dimensional dynamic program which yields the optimal policy in the setting of general correlations.  Finally, we note that our algorithm is entirely data-driven/sample-based - everything is implemented using only calls to $\textrm{SIM}$.  Of course, our formal algorithm/policy and the above sketch have some key differences, including modifications to the penalty formulation, recursions, and gradient definitions to accommodate greater generality and more sophisticated gradient methods and rounding.

\subsection{Literature review}\label{litsec}
\ \indent Much of the o.s.p. literature studies either adversarial models, or stochastic models in which rewards and r.c.v.s are drawn independently or as a random permutation (\cite{balseiro2023}), although hybrid models have also been studied (including the recent work \cite{guptaetal2026b} which leverages the robustness of sampling-based methods as we do for SGA).  In many such models one can achieve a constant regret (\cite{arlottoetal2019, veraetal2021,bumpensantietal2020,jasinetal2012,chen2025beyond}).  Several of these models allow for non-stationarity, a very large/infinite number of item types, and/or sampling/data-driven algorithms, and may have good empirical performance beyond the independence assumption (\cite{balseiro2023b,jiangetal2023,banerjeeetal2020,bray2025, besbesetal2025, jiangetal2025, zhang2026, zhang2026b,chen2025beyond}).  As our Example 1 illustrates, the gap between the optimal value and h.r. may be $\Theta(T)$ in our setting, making our results incomparable to these works (which typically prove a regret scaling $o(T)$).
\\\indent The linear relaxations of the packing problems we consider are examples of multistage stochastic linear programs (MSLP), which can be viewed as \textit{deterministic equivalent} massive LPs with a tree-like structure (\cite{olsen1976multistage}), and a variable for each partial trajectory of the o.i.p.  In general the size of this problem will be exponential in $T$ and $D$ (\cite{carpentieretal2015}).  Much of this literature first performs conditional multistage sampling to construct a so-called scenario tree and reduce the size of the problem, although the size of these trees scales exponentially in $T$ (\cite{heitschetal2009}).  Most solution methodologies in this literature, including progressive hedging and stochastic dual dynamic programming, have a complexity that scales linearly in the size of the scenario tree, resulting in an exponential dependence on $T$ (\cite{rockafellar1991scenarios,fullner2023stochastic}).  Recently these methods have been combined with sampling, which still results in a complexity scaling exponentially in $T$ unless one assumes independence (\cite{zhang2024sampling,mu2020convergence, zhao2005lagrangian, aydin2012sampling,lan2022complexity}).  Approaches based on IP and robust optimization (\cite{bertsimasetal2023,bertsimasetal2023b}) have a similar exponential complexity in the worst case.  Gradient methods have also been applied (\cite{lanetal2024}), including older works on quasi-gradient methods (\cite{ermoliev1988}) and more recent works such as \cite{cheungetal2000} and \cite{lanetal2017}, \cite{ahmed2006} and \cite{bieletal2021} which apply Nesterov smoothing, and \cite{zhaoetal1999} and \cite{hubneretal2017} which apply interior point methods.  To the best of our knowledge these methods again have a complexity scaling exponentially in $T$, as does related work in stochastic composite optimization (\cite{yangetal2019,zhangetal2021c,chenetal2025b,ghadimietal2020,zhangetal2024b}), conditional stochastic optimization (\cite{huetal2020}), and conditional compositional optimization (\cite{senetal2026}).  Several works suggest that such a dependence is likely unavoidable (\cite{shapiro2006}), while the recent works \cite{parketal2024,yoonetal2025} question this premise, and derive algorithmic results incomparable to our own.  When $T$ is $O(1)$, polynomial-time algorithms have been developed in \cite{swamy2012sampling,nemirovskietal2006,bavejaetal2023}.  Several recent works use generative AI to build simulation models in stochastic programming (\cite{dengetal2022,wangetal2022}), and speak to the growing relevance of complex generative models and simulators.
\\\indent In independent and concurrent work, \cite{zhangjaillet2025} also formulates an on-the-fly methodology for multistage stochastic programming.  They study convex problems generally with mirror descent and take a saddle point approach.  We study online packing with a smoothed penalty approach and make different modeling assumptions.  These differences lead to the analysis of \cite{zhangjaillet2025} applying more broadly, but introducing certain norms which may scale unfavorably for packing (see Section 6.2 and Appendix C of \cite{zhangjaillet2025}).  Furthermore, our work studies rounding for various applications and accelerated methods for packing, while \cite{zhangjaillet2025} does not study rounding, implements acceleration in different settings, and makes other contributions.  The works are thus incomparable and complementary.
\\\indent Our o.s.p. model can also be viewed as a high-dimensional stochastic control problem.  Several works in this literature have proven a polynomial complexity by imposing additional continuity assumptions (\cite{rust1997,belomestnyetal2025}), or incurring an exponential dependence on other parameters (\cite{becketal2025}), while in general dynamic programming will scale exponentially in $D$ (\cite{carpentieretal2015}).  Recently, \cite{goldberg2018polynomial} derived a polynomial-time algorithm for optimal stopping under our computational model.  Although these results were extended to some special cases of the problems we consider (\cite{chen2021}), our results are stronger and use very different techniques.  For example, for multiple stopping, the results of \cite{chen2021} are only polynomial-time when the number of stops is $O(1)$, in contrast to our results which allow the number of stops to scale linearly with $T$.  Ideas from statistical learning (such as VC dimension) have also been applied.  However, related settings where this approach has led to efficient algorithms typically either assume independence and/or prove that (for the problems considered) it suffices to learn a much simpler family of policies (e.g. threshold or static policies, see \cite{duttingetal2024,jinetal2024,xieetal2024}), or heuristically restrict to a family of policies and/or basis functions with small VC-dimension (without provable performance guarantees, see \cite{egloff2005,goeletal2017}).  Since it is easily verified that the family of policies we optimize over to solve o.s.p./c.a. in the setting of general correlations has exponential (in $T$) VC-dimension (and will in fact have cardinality doubly exponential in $T$ even for the setting of Example 2), and as to the best of our knowledge there are no simple families of policies guaranteed to perform well, our approach and results are incomparable to these works.
\\\indent Another relevant set of results pertains to the complexity of tabular RL and MDP under a generative model.  State-of-the-art complexity results here typically scale with the size of the state-space (which can be exponential $T$ and $D$), see \cite{sidford2018near,zureketal2024}, and matching lower bounds are known (\cite{kakade2003sample,azar2012sample}).  However, these bounds assume one must output a (near)-optimal action for every state and allow for general state-action transitions, in contrast to our work which only requires the DM to output a (near)-optimal action ``on-the-fly" and applies to packing.  Other past works avoid an exponential dependence on the dimension using sparse sampling, but incur an exponential dependence on $T$ (\cite{kearns2002sparse}).  Gradient methods have also been applied, and relevant recent works include \cite{tiapkinetal2022,chenetal2024b,abbasietal2019}.  Let us also note the works \cite{archibald2020,duetal2013,geiersbachetal2023} which use gradient methods to compute the optimal solution of stochastic control problems.
\\\indent The variant of NRM we study is identical to o.s.p./c.a., and is well-understood when the underlying uncertainty is independent / has strong concentration properties (\cite{balseiro2024survey}), although such an assumption is understood to be restrictive and imposed for tractability (\cite{gallego2019revenue}).  To address this, several works have put NRM in the framework of multistage stochastic programming, although no polynomial-time algorithms are derived (\cite{Demiguel2006,molleretal2008}).  Other works have used approximate dynamic programming (\cite{fariasetal2007}) and martingale duality (\cite{akanetal2009}).  Several recent works have relaxed the independence assumption by considering restricted dependency structures and deriving constant-factor approximations (\cite{aouadetal2022}) and/or proving asymptotic optimality (\cite{baietal2023}), see also \cite{ahnetal2025}.  Other works have incoporated Markovian uncertainty (\cite{jiang2023,lietal2025}), either proving constant-factor approximations (\cite{jiang2023}) or asymptotic guarantees (\cite{lietal2025,lietal2025b,lanetal2024b}), where some of the LP formulations considered in \cite{lietal2025} are essentially the same as those we consider.  These algorithms generally have a runtime depending on the number of states, which can be exponential in $T$ and $D$.  Gradient methods have also been applied (\cite{bertsimas2005simulation, van2008simulation, topaloglu2008stochastic}), often to optimize suboptimal heuristics.  
\\\indent There is a vast literature on online combinatorial optimization, a special case of o.s.p (see \cite{huangetal2024} for a recent survey on online matching).  For any fixed maximum degree $\Delta$, algorithms with competitive ratio better than $1-\frac{1}{e}$ (the bound from the seminal paper \cite{karpetal1990}) are known for online matching, but as $\Delta \rightarrow \infty$ the ratio converges to $1 - \frac{1}{e}$ (\cite{buchbinderetal2007,albersetal2022}).  The works \cite{Srinivasan2007,byrkaetal2018,naoretal2025} study such problems with general correlations using the results of \cite{swamy2012sampling}, and hence have an exponential dependence on $T$.  \cite{chen2021} extends the approach of \cite{goldberg2018polynomial} to online maximum weight bipartite independent set, a problem we also study.  Although that work develops a polynomial-time algorithm for approximating the optimal value, it uses a very complicated flow-based extension of \cite{goldberg2018polynomial}, restricts to the setting of a known graph, and does not yield an efficient policy.  Other recent works going beyond the independent setting include
\cite{heuseretal2025,feldmanetal2025,gaoetal2025}, which prove results incomparable to our own.  A recent line of work on so-called philosopher inequalities aims to derive approximation algorithms directly for a given online stochastic problem (relative to the optimal value of the associated MDP), see \cite{papadimitriouetal2021}.  These results typically assume independence (\cite{papadimitriouetal2021}), or only yield polynomial-time algorithms for fixed horizon (\cite{naoretal2025}).  Our results can be viewed as providing such philosopher inequalities for models with general correlations and long horizon.  \cite{vuongetal2025} studies several online combinatorial problems under general correlations in a model where one is given the different possible scenarios explicitly and runtimes may depend polynomially on the number of scenarios.  They derive a $1-\frac{1}{e}$-approximation for \textsc{mmo} (and many other problems) in that setting, and pose as an open question whether efficient algorithms can be derived when one instead has simulator access and there may be exponentially many scenarios.  Our results for \textsc{mmo} make partial progress on this question, although are incomparable to those of \cite{vuongetal2025} due to our additional assumptions on the weights and maximum degree.
\subsection{Outline of paper}
We formalize o.s.p./c.a. in Section \ \ref{setupsec}. In Section \ref{mainresultssec} we state our main results, provide further intuition for our runtimes and guarantees (and their limitations), and outline our proofs.  In the next 3 sections we prove results analogous to our main results for 3 relaxations of o.s.p. (building on each other, moving towards our main results) : a smoothed-penalty formulation (Section\ \ref{pensec1}), a non-smooth penalty formulation (Section\ \ref{penpensec1}), and the linear relaxation of o.s.p. (Section\ \ref{lpsec}).  We prove our main results for $\textrm{NRM}$ in Section\ \ref{nrmain1sec1}, online bipartite independent set in Section\ \ref{issec1}, and online bipartite matching in Section\ \ref{matchsec1}.  We provide concluding remarks in Section\ \ref{concsec1}, and an Electronic Compendium in Appendix Sections\ \ref{techsec1} and\ \ref{suppsec1}.

\section{Problem setup}\label{setupsec}
\subsection{Problem formulation}\label{formulationsec}
We now formulate our o.s.p./c.a. model more precisely.  We suppose there is a general stochastic o.i.p. $\lbrace \textrm{M}_t, t = 1,\ldots,T \rbrace$ with potentially non-Markovian and high-dimensional dynamics.  We assume $\textrm{M}_t \in {\mathcal R}^D$ w.p.1 for\ $t \geq 1$, for some dimension $D$ (potentially very large, scaling with other problem parameters e.g. $T,m,q$).  To simplify exposition (and without loss of generality, w.l.o.g.) we assume in each period there are exactly $q$ options (as opposed to at most $q$), as if there are less than $q$ options one can simply create a suitable number of ``dummy" options with zero reward and r.c.v.  In the setting that the o.i.p. is simply the realized rewards and r.c.v.s. of the $q$ options, we would take $D = q \times (m + 1)$, although we note that if the o.i.p. is appropriately encoded (e.g. if at each time only the positive components of the r.c.v.s, along with which of the $m$ resources those correspond to, are represented) one could in fact take $D = q \times (2 L + 1)$ (with $L$ an upper bound on the number of positive components of the r.c.v. of any given option).  In all cases we assume the rewards and r.c.v.s at time $t$ can be directly inferred from $M_t$.  In general our decisions (and the relevant conditional distributions) will depend not just on $\textrm{M}_t$ (which is $D$-dimensional) but on $\textrm{M}_{[t]} \stackrel{\Delta}{=} (M_1,\ldots,M_t)$ (which is $t \times D$-dimensional), since $\textrm{M}$ may be non-Markovian.  Let ${\mathcal S}$ denote the support of $\textrm{M}_{[T]}$, i.e. the set of all potential trajectories of the process, ${\mathcal S}^t$ the support of $\textrm{M}_{[t]}$, and ${\mathcal E} \stackrel{\Delta}{=} \bigcup_{t=1}^T {\mathcal S}^t$.  In the case that the o.i.p. is simply the realized rewards and r.c.v.s, ${\mathcal E}$ represents the set of all possible partial trajectories (of all possible lengths) of the realized rewards and r.c.v.s.  We assume $|{\mathcal E}| < \infty$, not because our results depend on $|{\mathcal E}|$, but because some of our results use arguments from convex optimization which are simpler in the finite setting.  Let $\mu : {\mathcal E} \rightarrow [0,1]$ denote the distribution function associated with $\textrm{M}$, i.e. $\mu(E) = \mathbb{P}\big( \textrm{M}_{[t]} = E \big)$ for $E \in {\mathcal S}^t$, where we remind the reader that we do not assume any knowledge of $\mu$ (our algorithms are entirely data-driven/sample-based).  We assume there are $m$ resources $\lbrace \textrm{RES}_i, i = 1,\ldots,m \rbrace$, each with a known budget $\lbrace b_i, i = 1,\ldots,m \rbrace$, and a known time horizon of $T$ periods.  In every time period $t$, the o.i.p. updates (to $\textrm{M}_t$), the new item arrives, and the $\textrm{DM}$ must irrevocably decide which of the $q$ options to pursue (using the information observed in times $1,\ldots,t$).  This repeats in each of the $T$ time periods.  The \textrm{DM}'s objective is to maximize the expected reward subject to satisfying the resource constraints w.p.1.
\\\indent \textbf{Massive IP and LP formulation.} 
One conceptually important idea (adopted in the stochastic programming community) is that the \textrm{DM}'s online decision-making problem can be formulated as a massive IP, denoted by \textsc{pack} (with optimal value $\textrm{OPT}_{\textsc{pack}}$).  For $E \in {\mathcal E}$, let $Z^r(E), a^r_i(E), X^r(E)$ denote (respectively) the reward, component $i$ of the r.c.v., and policy decision (for option $r$ on $E$).
\begin{align}
\tag{\textsc{pack}} \label{pack}
\max_{X} \quad & \sum_{E \in \mathcal{E}}  \mu(E) \sum_{r=1}^{q-1} Z^r(E)\, X^r(E) \\
\text{s.t.} \quad 
& \sum_{t=1}^T \sum_{r=1}^{q-1} a^r_i(S_{[t]})\, X^r(S_{[t]}) \leq b_i 
&& \forall\, i = 1, \ldots, m;\  S \in \mathcal{S} \notag \\
& X^r(E) \in \{0,1\}\ \textrm{and}\ \sum_{j=0}^{q-1} X^j(E) = 1 
&& \forall\, r = 0,\ldots,q-1, E \in \mathcal{E} \notag
\end{align}
\ \indent Any feasible solution to \ref{pack} is equivalently a feasible o.s.p./c.a. policy, and an optimal solution to \ref{pack} corresponds to an optimal policy. We denote the LP relaxation of \ref{pack}, in which the constraints $X^r(E) \in \{0,1\}$ are replaced by the constraints $X^r(E) \in [0,1]$, by \textsc{lp} (and the corresponding optimal value $\textrm{OPT}_{\textsc{lp}}$).  Although as noted previously it follows from \cite{naoretal2025} that \ref{pack} and \textsc{lp} have the same value for multiple stopping, they will in general differ (\cite{huangetal2024}).
\\\indent \textbf{Policy formulation.}  Alternatively (and equivalently), we may formalize the $\textrm{DM}$'s problem as a policy optimization, as follows.  For \ref{pack} (\textsc{lp}), a policy is a (possibly randomized) mapping $A: \mathcal{E} \times [0, 1] \to \{0, 1\}^q\ ([0,1]^q)$, which outputs a length-q binary (fractional) vector, whose components sum to unity, specifying which option is to be taken (fractionally) given as input any partial history $E \in \mathcal{E}$ and an independent random seed $\xi$ distributed uniformly on $[0, 1].$  For \ref{pack}, a vector with component $r$ set equal to 1 indicates that option $r$ is taken for that item; for \textsc{lp} the fractional value of component $r$ indicates what fraction of the item is subjected to option $r$.  A policy is said to be admissible if w.p. $1$ all packing constraints are respected, namely w.p.1 $\sum_{t= 1}^T \sum_{r=1}^{q-1} a^r_i(S_{[t]}) A^r(S_{[t]}, \xi) \leq b_i$ for all $i \in \lbrace 1,\ldots,m \rbrace$ and $S \in {\mathcal S}$, where the ``w.p. $1$'' is with respect to $\xi.$  For simplicity, we will henceforth suppress the explicit dependence of A on $\xi$, with the understanding that the policy may be randomized.  The expected reward earned by a policy $A$ can be written as $\mathbb{E}\big[ \sum_{E \in {\mathcal E}} \mu(E) \sum_{r=1}^{q-1} Z^r(E) A^r(E) \big]$, where (since the sum over $E \in {\mathcal E}$ is written explicitly) the expectation is over any randomization utilized by $A$.  As a notational convenience we will at times switch a bit informally between referring to $\textrm{A}$ as either a randomized feasible policy/algorithm, a random mapping from ${\mathcal E}$ to the appropriate range, or as a random $|{\mathcal E}|$-dimensional vector with component $E$ denoted either $\textrm{A}(E)$ or $\textrm{A}_E.$
\\\textbf{Overview of Assumptions}.  In this paper we make the following assumptions, and provide some additional details and discussion in Section\ \ref{compmodelsec} and Appendix\ \ref{appassumptionsec}. 
\vspace{-.05in}
\begin{assumption}[Normalized rewards and r.c.v.]
$a^r_{i,t}, Z^r_t \in [0,1]$\ $\forall$\ $r,i,t$ w.p.1.  \label{assumption1}
\end{assumption}
\vspace{-.15in}
\begin{assumption}[Bounded column sparsity]
$|\lbrace i : a^r_{i,t} \neq 0 \rbrace| \leq L$ (for a fixed $L$)\ $\forall$\ $r,t$ w.p.1. In matching and independent set applications, we use $\Delta$ in place of $L$.
\label{assumption2}   
\end{assumption}
\vspace{-.15in}
\begin{assumption}[Non-infinitesimal consumption]
Either $a^r_{i,t} = 0$ or  $a^r_{i,t} \geq \iota$ for a fixed $\iota \in (0,1]$\ $\forall$\ $r,i,t$ w.p.1.   \label{assumption3}    
\end{assumption}
\vspace{-.15in}
\begin{assumption}[Simulator access]
 At a computational cost ${\mathcal C} \geq 1$, the $\textrm{DM}$ can input any partial trajectory (i.e. prefix) of the o.i.p., and get an independent draw of the remaining trajectory of the o.i.p., drawn from the appropriate conditional distribution.  In general ${\mathcal C}$ will depend on both $T$ and $D$, as dictated by the complexity of implementing the relevant conditional simulation.
\label{assumption4}   
\end{assumption}
\subsection{Applications to NRM and online combinatorial optimization}\label{applicationssec}
\ \indent \textbf{NRM} (with flexible products).  The problem of $\textsc{nrm}$ is identical to $\ref{pack}$, as discussed in Section\ \ref{introsec1}.  Here $q$ corresponds to the number of different ways to satisfy a given customer, and $L$ to the maximum number of distinct resources required under any given option to satisfy a customer.  Our modeling framework allows for no-shows (in which some periods have no arrival), as well as an effectively random time horizon, by having the o.i.p. dictate that in such time periods both the reward and r.c.v.s are all zeros.  We denote the optimal value of such an $\textsc{nrm}$ problem modeled as a packing problem by $\textrm{OPT}_{\textsc{nrm}}$, and note that at this level of generality \textsc{nrm} is equivalent to o.c.a. (with more classical non-flexible variants of \textsc{nrm} equivalent to o.c.a. with $q = 2$)
\\\indent \textbf{Online Bipartite Max Weight Independent Set (\textsc{is})}. \textsc{is} is a special case of \ref{pack}. Formally, suppose there is an unknown $n$-node bipartite graph $G$, with known partites (node sets) $\mathcal{L}$ and $\mathcal{R}$ satisfying $|\mathcal{L} \cup \mathcal{R}| = n.$ There is a known upper bound $\Delta$ on the degree of any node, but the edges and node weights are initially unknown and will be revealed sequentially (node-by-node, in a known fixed order).  In each period $t$, $\textrm{M}_{t}$ is realized, and the weight of node $t$, as well as the set of edges containing node $t$ are revealed.  When the set of edges containing a given node is revealed, the $\textrm{DM}$ will not necessarily learn the identities of the other nodes appearing in those edges (which may not be revealed until those nodes are themselves visited).  The \textrm{DM} must determine whether to include node $t$, subject to the constraint that no two included nodes belong to the same edge. 
\\\indent In the language of \ref{pack}, $q = 2$, and there are $T = n$ time periods. For each partial trajectory $S_{[t]}$, there is a binary variable determining whether node $t$ is included given $S_{[t]}$, and $Z^1(S_{[t]})$ denotes the weight of node $t$ given $S_{[t]}$.  There are $m = \lfloor \frac{1}{2}\Delta n\rfloor$ resource constraints, one for each potential edge (as a graph with maximum degree $\Delta$ can have at most $\lfloor \frac{1}{2} \Delta n \rfloor$ edges).  For each potential edge $i$ and trajectory $S$, $a^1_i(S_{[t]}) = 1(0)$ if node $t$ belongs (does not belong) to edge $i$ on trajectory $S$, and $\sum_{t = 1}^T a^1_i(S_{[t]})$ equals either 2 or 0  (as required by the graph structure), where $\sum_{t = 1}^T a^1_i(S_{[t]}) = 0$ indicates that edge $i$ is not realized on trajectory $S$.  Setting $b_i = 1$ for all $i$ enforces the independent set constraints.  We denote the optimal value of such an \textsc{is} problem by $\textrm{OPT}_{\textsc{is}}$.
\\\indent \textbf{Online Maximum Weight Bipartite Matching (\textsc{mwm}) and its fractional relaxation (\textsc{mwmlp})}. 
\textsc{mwm}(\textsc{mwmlp}) is a special case of \ref{pack} (\textsc{lp}), and the basic setup is similar to \textsc{is}, with $G = \mathcal{L} \cup \mathcal{R}$ satisfying $|G| = n$.  $T = \lfloor \frac{1}{2}\Delta n\rfloor$ (potential) edges arrive sequentially.  In each period $t$, $\textrm{M}_t$ is realized, either revealing the two nodes that constitute potential edge $t$, or revealing that potential edge $t$ is never realized (in which case potential edges $t+1, \ldots, T$ are also never realized).  The \textrm{DM} must then determine whether to include the edge in the matching or not (in the fractional problem one must assign a fractional value to that edge), subject to the constraint that the sum of the (fractional) values assigned to the edges incident to any given node is at most 1.
\\\indent In the language of \ref{pack}, $q = 2$, and there are $T = \lfloor\frac{1}{2} \Delta n\rfloor$ time periods. For each partial trajectory $S_{[t]}$, there is a variable determining the value assigned to edge $t$ given $S_{[t]}$, and $Z^1(S_{[t]})$ denotes the weight of potential edge $t$ given $S_{[t]}$ if realized (and 0 otherwise).  There are $n$ resource constraints, one for each node.  For each potential edge $t$ and trajectory $S$, $a^1_i(S_{[t]}) = 1$ if edge $t$ is realized and incident to node $i$, and equals 0 otherwise; and $\sum_{i = 1}^n a^1_i(S_{[t]})$ equals either 2 or 0  (as required by the graph structure), where $\sum_{i = 1}^n a^1_i(S_{[t]}) = 0$ indicates that edge $t$ is not realized on trajectory $S$.  Setting $b_i = 1$ for all $i$ enforces the matching constraints.  We denote the optimal value of such a \textsc{mwm} (\textsc{mwmlp}) problem by $\textrm{OPT}_{\textsc{mwm}}$ ($\textrm{OPT}_{\textsc{mwmlp}}$).
\\\indent\textbf{Online Maximum Weight Bipartite Matching with Online nodes (\textsc{mmo})}.  \textsc{mmo} is a special case of \textsc{mwm}, in which instead of edges being revealed one at a time, in each time period a new ``online" node in partite ${\mathcal R}$ arrives and all of its incident edges (and their weights) are revealed simultaneously, revealing which of the ``offline" nodes of partite ${\mathcal L}$ are incident to it.  At that time an irrevocable decision must be made about which one, if any, of those edges is selected into the matching.  The nodes in partite ${\mathcal R}$ are referred to as the ``online nodes", and such an information structure is common in the literature (\cite{gamlathetal2019}).  $\textsc{mmo}$ can be modeled as an instance of o.c.a.  In the language of \ref{pack}, there are $T = |{\mathcal R}|$ time periods, and in each period $q = \Delta + 1$ options (corresponding to which, if any, of the at most $\Delta$ revealed edges are selected).  For each partial trajectory $S_{[t]}$, there are $\Delta + 1$ variables determining the values assigned to the (at most) $\Delta$ edges incident to node $t$ given $S_{[t]}$ (the $+1$ is for the option of selecting no edges), with $Z^r(S_{[t]})$ the weight of potential edge $r$ (incident to node $t$) given $S_{[t]}$ if realized.  There are $|{\mathcal L}|$ resource constraints, one for each offline node.  For each node $t$, potential edge $r$, and trajectory $S$, $a^r_i(S_{[t]}) = 1$ if the $r$th edge incident to online node $t$ is incident to offline node $i$, and equals 0 otherwise.  Setting $b_i = 1$ for all $i$ enforces the matching constraints.  We denote the optimal value of \textsc{mmo} by $\textrm{OPT}_{\textsc{mmo}}.$  
\subsection{Computational model}\label{compmodelsec}
Our computational model is in line with several past works in stochastic optimization (\cite{goldberg2018polynomial, de2008complexity}), where standard arithmetic operations each require unit time, and memory access costs are ignored.  We assume that sampling $k$ indices without replacement from $\{1,\ldots,T\}$ takes $O(k)$ time (independent of $T$), as proven in \cite{ting2021}, and that sampling either a Uniform, Bernoulli, or Beta r.v. may be done in unit time.  Our model is defined by access to a simulator, operating at a cost of ${\mathcal C} \geq 1$ time units per call.  Simulator $\textrm{SIM}$ takes a partial trajectory $E \in \mathcal{S}^t$ (a $D \times t$ matrix, for any $t$) and returns a complete trajectory $\hat{S} \in \mathcal{S}$ drawn from the random distribution of trajectories conditional on the prefix $E$.
\section{Main results and proof outline}\label{mainresultssec}
\subsection{Main results}
Our results will be stated in terms of an absolute constant $c_0$, which is some number (independent of any problem specifics/parameters) that could in principle be made explicit.  Let $\nu \stackrel{\Delta}{=} \min_{i = 1,\ldots,m} \frac{b_i}{T}$ denote the ratio of the minimum budget to $T$, so $\nu$ is $\Theta(1)$ when all budgets scale linearly in $T$.  Let $\kappa$ denote an upper bound on the maximum number of options to utilize any one resource at any one time $t$ : $\kappa \stackrel{\Delta}{=} \max_{E \in {\mathcal E} , i \in \lbrace 1,\ldots,m \rbrace} |\lbrace r : a^r_i(E) > 0 \rbrace |$ (so $\kappa = 1$ if different options always utilize disjoint resources, and $\kappa = q$ if at some time all options may utilize the same resource).  When we refer to the per-decision time required by an algorithm, throughout this should be understood as describing the total time required for all computation and simulation.
\\\indent We first state our main result for $\textsc{nrm}$, also encompassing our main result for o.c.a. (and is stated at the level of generality required for multiple options).  Here we state a simplified result with worse dependence on certain parameters to maximize readability, and defer a statement of our stronger results to Section\ \ref{nrmain1sec1}.  Our results for $\textsc{NRM}$ utilize a policy $\textrm{A}_{\textsc{nrm}}$ (defined in Section\ \ref{nrmain1sec1}) which is the formalization and generalization (from multiple stopping to \textsc{nrm}) of the on-the-fly SGA policy sketched in Section\ \ref{introsec1} (combined with a randomized rounding scheme).  For clarity of exposition we state two runtime guarantees for $\textrm{A}_{\textsc{nrm}}$, instead of a single more complicated guarantee.

\begin{theorem}\label{NRMAIN1}
$\forall\ \epsilon \in (0, 1)$, $\textrm{A}_{\textsc{nrm}}$ is an admissible policy for $\textsc{nrm}$ with the following properties.
\begin{itemize}
\item On any trajectory $S \in \mathcal{S}$ one can compute decisions $\textrm{A}_{\textsc{nrm}}(S_{[t]})\ \forall\ t = 1, \dots, T$ on-the-fly, in per-decision time ${\mathcal C} \times \exp\bigg( c_0 (\frac{Lq}{\iota})^3 \lceil \log(\frac{m}{\epsilon}) \rceil \epsilon^{-1} \min(\sqrt{m},\epsilon^{-1}) \bigg),$
s.t. $\mathbb{E}\big[ \sum_{E \in {\mathcal E}} \mu(E) \sum_{r=1}^{q-1} Z^r(E) \textrm{A}_{\textsc{nrm}}^r(E) \big] \geq \textrm{OPT}_{\textsc{nrm}} - \epsilon T$ so long as $T \geq c_0 (\frac{Lq}{\iota})^3 m \epsilon^{-2}$.
\item On any trajectory $S \in \mathcal{S}$ one can compute decisions $\textrm{A}_{\textsc{nrm}}(S_{[t]})\ \forall\ t = 1, \dots, T$ on-the-fly, in per-decision time ${\mathcal C} \times 
\exp\bigg( c_0 (\frac{L \kappa}{\iota \nu})^2 \lceil \log(\frac{q}{\epsilon}) \rceil \epsilon^{-1}\bigg),$ s.t. $\mathbb{E}\big[ \sum_{E \in {\mathcal E}} \mu(E) \sum_{r=1}^{q-1} Z^r(E) \textrm{A}_{\textsc{nrm}}^r(E) \big] \geq \textrm{OPT}_{\textsc{nrm}} - \epsilon T$
so long as $T \geq c_0 (\frac{L \kappa}{\iota \nu})^2  \log(2m) \epsilon^{-2}$.
\end{itemize}
\end{theorem}
We next state our main results for $\textsc{is}, \textsc{mwmlp}$, and $\textsc{mmo}.$  These results utilize policies $\textrm{A}_{\textsc{is}}$ (defined in Section\ \ref{issec1}), $\textrm{A}_{\textsc{mwmlp}}$ (defined in Section\ \ref{matchsec1}), and $\textrm{A}_{\textsc{mmo}}$ (defined in Section\ \ref{matchsec1}), all of which are very similar to $\textrm{A}_{\textsc{nrm}}$ (sometimes combined with more sophisticated rounding schemes).
\begin{theorem}\label{NRMAINCOMBOPT}\ 
\begin{itemize}
\item $\forall\ \epsilon \in (0, 1)$, $\textrm{A}_{\textsc{is}}$ is an admissible policy for $\textsc{is}$ with the following properties.  On any trajectory $S \in \mathcal{S}$ one can compute decisions $\textrm{A}_{\textsc{is}}(S_{[t]})\ \forall\ t \in \lbrace 1, \dots,n \rbrace$ on-the-fly, in per-decision time ${\mathcal C} \times \big( c_0 \frac{\Delta}{\epsilon} \big)^{c_0 \Delta \epsilon^{-1}}$, and $\mathbb{E}\big[ \sum_{E \in {\mathcal E}} \mu(E) Z^1(E) \textrm{A}^1_{\textsc{is}}(E) \big] \geq \textrm{OPT}_{\textsc{is}} - \epsilon n$.
\item $\forall\ \epsilon \in (0, 1)$, $\textrm{A}_{\textsc{mwmlp}}$ is an admissible policy for $\textsc{mwmlp}$ with the following properties.  On any trajectory $S \in \mathcal{S}$ one can compute decisions $\textrm{A}_{\textsc{mwmlp}}(S_{[t]})\ \forall\ t \in \lbrace 1,\ldots,\lfloor \frac{1}{2} \Delta n \rfloor \rbrace$ on-the-fly, in per-decision time ${\mathcal C} \times \big( c_0 \frac{\Delta}{\epsilon} \big)^{c_0 \Delta \epsilon^{-1}}$, and $\mathbb{E}\big[ \sum_{E \in {\mathcal E}} \mu(E) Z^1(E) \textrm{A}^1_{\textsc{mwmlp}}(E) \big] \geq \textrm{OPT}_{\textsc{mwmlp}} - \epsilon n$.
\item $\forall\ \epsilon \in (0, 1)$, $\textrm{A}_{\textsc{mmo}}$ is an admissible policy for $\textsc{mmo}$ with the following properties.  On any trajectory $S \in \mathcal{S}$ one can compute decisions $\textrm{A}_{\textsc{mmo}}(S_{[t]})\ \forall\ t \in \lbrace 1,\ldots,|{\mathcal R}| \rbrace$ on-the-fly, in per-decision time ${\mathcal C} \times \big( c_0 \frac{\Delta}{\epsilon} \big)^{c_0 \sqrt{\Delta} \epsilon^{-1}}$, and $\mathbb{E}\big[ \sum_{E \in {\mathcal E}} \mu(E) Z^1(E) \textrm{A}^1_{\textsc{mmo}}(E) \big] \geq .68 \times \textrm{OPT}_{\textsc{mmo}} - \epsilon n$.
\end{itemize}
\end{theorem}
\subsection{Discussion of runtimes, $\epsilon T$ scaling of error, and limitations of our approach} We begin by discussing the runtime of $\textrm{A}_{\textsc{nrm}}$, and for this discussion ignore certain constants and logarithmic terms.  Supposing $L,\iota^{-1},q$ are small constants, our first result for $\textrm{A}_{\textsc{nrm}}$ yields a runtime roughly ${\mathcal C} \times \min\big(m^{\epsilon^{-2}}, m^{\sqrt{m} \epsilon^{-1}}\big)$, and $\epsilon T$ optimality gap so long as $m$ is sublinear in $T$ (i.e. $m \leq \epsilon^2 T$).  For multiple stopping, these complexity bounds reduce to a ${\mathcal C} \times \exp\big( \epsilon^{-1}\big)$ runtime so long as $T \geq \epsilon^{-2}$.  These results make no assumption about how the budgets scale with $T$ and have no dependence on $\nu$.  Supposing instead $L,\iota^{-1},\kappa,\nu$ are small constants (implying the resource budgets scale linearly in $T$), our second result for $\textrm{A}_{\textsc{nrm}}$ yields a runtime roughly ${\mathcal C} \times q^{\epsilon^{-1}}$, and $\epsilon T$ optimality gap so long as $m$ is subexponential in $T$ (i.e. $m \leq \exp(\epsilon^2 T)$).  Thus when the budgets are $\Theta(T)$, $\textrm{A}_{\textsc{nrm}}$ is efficient even when $m$ and $q$ are polynomially large in $T$ (as long as $\kappa$ is $O(1)$).  For $\textsc{is}$ and $\textsc{mwmlp}$, our algorithm's per-decision complexity scales as $(\frac{\Delta}{\epsilon})^{\Delta \epsilon^{-1}}$; while for $\textsc{mmo}$ it scales as $(\frac{\Delta}{\epsilon})^{\sqrt{\Delta} \epsilon^{-1}}$ (the $\sqrt{\Delta}$ speedup arising since all edges for an online node are revealed at once).
\\\indent Next, let us comment on the nonstandard $\epsilon T$ error.  Here $\epsilon$ should be thought of as a small fixed constant independent of $T$ (with $T$ large), although we note that $\epsilon$ can even go to zero as $T$ gets large at a slow rate (e.g. $O\big(\frac{\log(\log(T))}{\log(T)}\big)$) while ensuring our runtimes remain polynomial in $T$.  The $\Theta(T)$ gap between the optimal value and the h.r. (as shown in Example 1) explains why an $\epsilon T$ additive error (for arbitrarily small fixed $\epsilon$ and large $T$) relative to the h.r. is unachievable in our setting.  Let us also point out that since in Example 1 the budget scales linearly with $T$ and there is only one resource yet such a gap persists, a large budget assumption (as required in \cite{agrawaletal2014} to ensure closeness to the h.r.) will not remedy the situation.  There are results in the literature showing that under certain assumptions one can achieve a constant competitive ratio relative to the h.r. (e.g. in adversarial settings, see \cite{balseiro2023}, which would also apply in our setting), and such results are incomparable to our own.  For example, for multiple stopping this competitive ratio vanishes as $\frac{b_1}{T}$ goes to zero, and we note that for (the linear relaxation of) multiple stopping a competitive ratio of $\frac{b_1}{T}$ can anyways be achieved by setting $X^1_t = \frac{b_1}{T}\ \forall\ t$.
\\\indent Our results are most natural in settings where the optimal value itself is $\Theta(T)$ (with $T$ large and $\epsilon$ a fixed small constant), since then our algorithms also yield polynomial-time approximation schemes (PTAS, with additive error $\epsilon \times \textrm{OPT}_{\textsc{pack}}$) in essentially the same runtimes.  Much of the literature on \textsc{nrm} imposes assumptions which indeed imply this (e.g. budgets scaling linearly in $T$ with expected rewards bounded away from 0, see \cite{bumpensantietal2020}).  Even without large budgets, there are many other natural settings where the optimal value would also scale as $\Theta(T)$, such as maximum cardinality matching in bounded degree graphs with no isolated nodes, and more broadly settings in which for any given resource, the ratio of its budget to the maximum number of times it is requested is bounded away from zero.  Due to the data-driven nature of our model, such a PTAS-type result would be impossible in general, at least not without making some additional statistical assumptions, since otherwise such an algorithm would be able to distinguish between two exponentially unlikely events with only polynomially many samples, as we explore further in Appendix\ \ref{noptasstats}.  Past works achieving such a PTAS-style result for related models (e.g. \cite{duttingetal2023,agrawaletal2014,goldberg2018polynomial}) 
typically exploit special structure and assumptions that do not hold in our model, as we also explore in Appendix\ \ref{noptasstats}.
The reason $T$ appears in our error bounds is that certain quantities arising in convex optimization (e.g. the norms of solutions and gradients, and smoothness parameters) seem to have a fundamental dependence on $T$.  Similarly, the $\textrm{poly}(\frac{1}{\epsilon})$ appearing in the exponent in our complexity analysis is inherited from the number of iterations required by SGA (which is exponentiated due to the recursive nature of our algorithm).  We leave as an open question whether these dependencies are fundamental.
\\\indent Finally, let us emphasize several limitations of our results.  First, ${\mathcal C}$ will in general be at least polynomial in $T$ and $D$, and thus our algorithms' per-decision complexity is not truly independent of $T$ and $D$.  Second, supposing that ${\mathcal C}$ is polynomial in $T$ and $D$, our algorithms' runtimes will scale roughly as $\textrm{poly}(T,D) \times e^{c_0 \epsilon^{-1}}.$  Alternatively, many past methods have a complexity which is exponential in $T$ (but polynomial in $\epsilon^{-1}$ for fixed $T$).  Thus for $T \leq \frac{c_0}{\epsilon}$, those methods may be both faster and more accurate.  If the problem has additional structure (e.g. satisfies stronger notions of continuity), then past methods can achieve a superior polynomial complexity broadly.  Thus our algorithms would be most relevant when : (1) $\epsilon$ is small but fixed independent of $T$ (so high accuracy is not required), (2) the horizon is very long, (3) the optimal value scales linearly in $T$, and (4) the stochasticity is unstructured and discontinuous.  Furthermore, given the high computational cost of implementing such policies, in many cases either algorithms coming from worst-case analysis, or heuristics based on intuitions from the independent setting, approximate dynamic programming, or neural networks may be preferred.  In this light, our results should be viewed as a proof-of-concept that o.s.p. with general correlations is theoretically (although not practically) tractable in certain settings where past methods either required exponential time or did not have theoretical guarantees.

\subsection{Proof outline}\label{proofoutlinesec} First, in Section\ \ref{pensec1}, we define a smoothed-penalty relaxation \ref{pen-smooth} of \textsc{lp}.  We derive an on-the-fly algorithm $\textrm{A}_{\textsc{pen}^{\theta}}$ for \ref{pen-smooth}, proving a result analogous to Theorem\ \ref{NRMAIN1}, as follows.  First, we prove (using standard results in convex optimization) that a family of SGA algorithms (including accelerated and unaccelerated variants with different types of sampling), run on the massive problem \ref{pen-smooth}, yields an $\epsilon T$-approximately optimal solution (in expectation) in an appropriately bounded number of iterations.  Second, we prove that a simple subroutine $\textrm{R}$ can efficiently compute the value of any one variable in the above SGA methods after any given number of iterations, and combine with some additional analysis to complete the proof.
\\\indent Second, in Section\ \ref{penpensec1}, we define a non-smooth penalty relaxation \ref{pen} of \textsc{lp} to act as a bridge between \ref{pen-smooth} and \textsc{lp} (as it is easier to control the resource constraint violation).  We derive an on-the-fly algorithm $\textrm{A}_{{\textsc{pen}}}$ for \ref{pen}, proving a result analogous to Theorem\ \ref{NRMAIN1}, by leveraging our previous result for \ref{pen-smooth} and showing that the solution returned by $\textrm{A}_{\textsc{pen}^{\theta}}$ is also nearly optimal for \textsc{pen} (and nearly feasible for \textsc{lp}) for appropriate $\theta$.
\\\indent Third, in Section\ \ref{lpsec}, we derive an on-the-fly algorithm $\textrm{A}_{{\textsc{lp}}}$ for $\textsc{lp}$, proving a result analogous to Theorem\ \ref{NRMAIN1} by patching the infeasibility in our algorithm for \ref{pen} (formalized through a function $\textrm{F}$).
\\\indent Fourth, in Section\ \ref{nrmain1sec1}, we prove our main result for $\textrm{NRM}$ Theorem\ \ref{NRMAIN1} (using on-the-fly algorithm $\textrm{A}_{\textsc{nrm}}$) by combining the general logic of our algorithm and proof for $\textsc{lp}$ with a randomized rounding scheme (independent across time, dependent across options, formalized with functions $\textrm{H}$ and $\textrm{J}$).  
\\\indent Finally, in Sections\ \ref{issec1} and \ref{matchsec1}, we prove our main results for \textsc{is}, \textsc{mwmlp}, \textsc{mmo} (Theorem \ref{NRMAINCOMBOPT}, using on-the-fly algorithms $\textrm{A}_{\textsc{is}}$, $\textrm{A}_{\textsc{mwmlp}}$, $\textrm{A}_{\textsc{mmo}}$) by combining $\textrm{A}_{{\textsc{lp}}}$ with additional rounding schemes.  For \textsc{is}, in Section\ \ref{issec1} we use a dependent rounding scheme which compares the fractional values returned by our algorithms to a single common uniform r.v.  For \textsc{mmo}, in Section\ \ref{matchsec1} we plug the values returned by our algorithms into the online rounding algorithm of \cite{harrisetal2026}.

\section{On-the-fly algorithm for smoothed-penalty formulation}\label{pensec1}
\subsection{Additional definitions and notations}
\begin{itemize}
\item $a^{r,+}(E) \stackrel{\Delta}{=} \{i : a^r_i(E) > 0\}$ : requested resources for option $r$ on $E \in {\mathcal E}$
\item $a^+(E) \stackrel{\Delta}{=} \bigcup_{r=1}^{q-1} a^{r,+}(E)$ : requested resources across all options on $E \in {\mathcal E}$
\item $\mathcal{T}_i(S) \stackrel{\Delta}{=} \bigcup_{r=1}^{q-1}\{t : a^r_i(S_{[t]}) > 0\}$ : times $\textrm{RES}_i$ is requested (by at least one option) on $S \in {\mathcal S}$
\item $\mathcal{RT}_i(S) \stackrel{\Delta}{=} \lbrace (r,t) : a^r_i(S_{[t]}) > 0 \rbrace$ : $(r,t)$ s.t. $\textrm{RES}_i$ is requested by option $r$ at time $t$ on $S \in {\mathcal S}$
\item $U \stackrel{\Delta}{=} \max_{S \in \mathcal{S}, i \in \{1,\ldots,m\}} |\mathcal{T}_i(S)|$ : bound on $|\mathcal{T}_i(S)|$, where we note that $U \leq T$
\item $U_r \stackrel{\Delta}{=} \max_{S \in \mathcal{S}, i \in \{1,\ldots,m\}} |\mathcal{RT}_i(S)|$ : bound on $|\mathcal{RT}_i(S)|$, where we note that $U_r \leq \kappa T$
\item For $E \in {\mathcal S}^t, E' \in {\mathcal S}^{t'}$ with $t' \leq t$, we say $E' \subseteq E$ if $E' = E_{[t']}$, i.e. $E'$ is a prefix of $E$ 
\item $W \stackrel{\Delta}{=} \max_{S \in \mathcal{S}, S' \subseteq S, r \in \lbrace 0,\ldots,q-1 \rbrace} \sum_{t=1}^T \sum_{l = 1}^{q-1} |a^{l,+}(S_{[t]}) \cap a^{r,+}(S')|$ : bound on the total resource overlap between any one option for any one arriving item and all the r.c.v.s associated with all options of all arrivals on any one trajectory, where we note that $W \leq L U_r$
\item ${\mathcal Q}^0 \stackrel{\Delta}{=} \bigg\lbrace X^r, r = 0,\ldots,q-1\ \textrm{s.t.}\ X^r \in [0,1]\ \forall\ r; \sum_{r=0}^{q-1} X^r = 1\bigg\rbrace$ : $q$-dimensional simplex 
\item ${\mathcal Q}^1 \stackrel{\Delta}{=} \bigg\lbrace X^r_t, r = 0,\ldots,q-1, t = 1,\ldots,T\ \textrm{s.t.}\ X^r_t \in [0,1]\ \forall\ r,t ; \sum_{r=0}^{q-1} X^r_t = 1\ \forall t \bigg\rbrace$ : $q \times T$-dimensional polytope corresponding to the set of (possibly non-admissible) static policies
\item ${\mathcal Q}^2 \stackrel{\Delta}{=} \bigg\lbrace X^r(E), r = 0,\ldots,q-1, E \in {\mathcal E}\ \textrm{s.t.}\ X^r(E) \in [0,1]\ \forall\ r,E ; \sum_{r=0}^{q-1} X^r(E) = 1\ \forall E \bigg\rbrace$ : $q \times |{\mathcal E}|$-dimensional polytope corresponding to the set of all (possibly non-admissible) policies
\item For polytope $Q$ and $\tau \in {\mathcal R}^+$, $\tau Q \stackrel{\Delta}{=} \bigcup_{X \in Q} \tau X$ denotes that scaling of $Q$ (also a polytope).  
\item For polytope $Q$ and vector $X$ of the same dimension, let $\Pi_Q(X)$ denote the Euclidean projection of $X$ onto $Q$, where we note that $\Pi_{\tau {\mathcal Q}^0}(X)$ can be computed exactly in time $O\big((q \log(q)\big)$ in our computational model (independent of the value of $\tau$ and $X$, see \cite{chenetal2011})
\item $V \stackrel{\Delta}{=} \max_{S \in \mathcal{S} , \lbrace X^r_t \rbrace \in {\mathcal Q}^1} \left| \left\{ i : \sum_{t=1}^T \sum_{r=1}^{q-1} X^r_t a^r_i(S_{[t]}) \geq \frac{1}{2} b_i \right\} \right|$ : bound on the number of resources that can have their budget at least half-saturated on any one trajectory, where we note that $V \leq m$
\item $I(\cdot)$ : indicator of the corresponding event
\item $\overline{e}^q$ : length-q vector with component 0 set to 1 and all $q-1$ other components set to 0
\end{itemize}

\subsection{Smoothed-penalty formulation} For $\theta \in (0,T]$, let $\phi_{\theta} : {\mathcal R} \rightarrow {\mathcal R}$ (with derivative $\phi'_{\theta}$) be the following one-sided Huber loss : 
\[
\phi_{\theta}(x) = \begin{cases}
0 & \text{if } x \leq 0,\\
\frac{1}{2 \theta} x^2  & \text{if } x \in [0, \theta],\\
x - \frac{1}{2} \theta  & \text{if } x > \theta.
\end{cases}
\]
Let $f^{\theta} : {\mathcal R}^{q \times |{\mathcal E}|} \rightarrow {\mathcal R}$ denote the mapping
$$f^{\theta}(X) \stackrel{\Delta}{=} \sum_{E \in {\mathcal E}} \mu(E) \sum_{r=1}^{q-1} Z^r(E) X^r(E) - 2 \iota^{-1} \sum_{S \in {\mathcal S}} \mu(S) \sum_{i=1}^m \phi_{\theta}\big( \sum_{t=1}^T \sum_{r=1}^{q-1} a^r_i(S_{[t]}) X^r(S_{[t]}) - b_i \big).$$  Let $\textsc{pen}^{\theta}$ denote the following concave program, with optimal value $\textrm{OPT}_{\textsc{pen}^{\theta}}$ and solution $X_{*,\textsc{pen}^{\theta}}.$
\vspace{-.15in}
\begin{align}
    \tag{$\textsc{pen}^{\theta}$} \label{pen-smooth}
    \max f^{\theta}(X) \quad s.t. \ X \in {\mathcal Q}^2
\end{align}
Our main result for \ref{pen-smooth} (below) utilizes policy $\textrm{A}_{\textsc{pen}^{\theta}}$ which is the formalization and generalization (from multiple stopping to \textsc{nrm}) of the on-the-fly SGA policy sketched in Section\ \ref{introsec1}, and which we formally define at the end of this section.

\begin{theorem}\label{convergence1}
$\forall\ \epsilon \in (0, 1)$, on any trajectory $S \in \mathcal{S}$ one can compute decisions $\textrm{A}_{\textsc{pen}^{\theta}}(S_{[t]})\ \forall\ t = 1, \dots, T$ on-the-fly, in per-decision time 
$$
{\mathcal C} \times  \min \Bigg( \big( c_2 \frac{L T q}{\iota \theta \epsilon} \big)^{c_2 (\frac{L}{\iota \epsilon})^{2} q} , \big( c_2 \frac{L T q}{\iota \theta \epsilon} \big)^{c_2 \lceil \frac{(U_r L W)^{\frac{1}{4}}}{\sqrt{\iota \theta}}\rceil \epsilon^{-\frac{1}{2}}} , \big( c_2 \frac{L U q}{\iota \epsilon} \big)^{c_2 \lceil \frac{(U_r L W)^{\frac{1}{4}}}{\sqrt{\iota \theta}}\rceil \epsilon^{-\frac{1}{2}}} \Bigg),
$$
where $c_2$ is an absolute constant. Also, $\mathbb{E}\big[ \sum_{E \in {\mathcal E}} \mu(E) \sum_{r=1}^{q-1} Z^r(E) \textrm{A}_{\textsc{pen}^{\theta}}^r(E) \big] \geq \textrm{OPT}_{\textsc{pen}^{\theta}} - \epsilon T.$ 
\end{theorem}
To prove Theorem\ \ref{convergence1}, we first analyze the aforementioned family of gradient methods for the massive problem \ref{pen-smooth}, parametrized by a step-size $\alpha$, a set of non-negative ``momentum constants" $\lbrace \beta_k, k \geq 0 \rbrace$ (allowing us to consider both accelerated and unaccelerated methods), and sampling parameters $\eta_1,\eta_2 \in Z^+.$  For $k \geq 0$, we define a random vector-valued function $\hat{G}^k: \mathcal{R}^{q \times |\mathcal{E}|} \to \mathcal{R}^{q \times |\mathcal{E}|}$, which will act as a (biased) stochastic gradient (in iteration $k$ of the gradient method).  For any $X \in \mathcal{R}^{q \times |\mathcal{E}|}$ and $E \in \mathcal{E}, l \in \lbrace 0,\ldots,q-1 \rbrace$, the $(E,l)$-th coordinate of $\hat{G}^k$ acting on $X$ is 
\begin{equation}\label{ghatkdef}
\hat{G}^k(X)_{E,l} = 
Z^l(E) - 2 \iota^{-1} \sum_{i =1}^m a^l_{i}(E) \times \eta^{-1}_1 \sum_{S' \in {\mathcal S}^{E,k}}\phi'_{\theta}\bigg( \frac{T}{\eta_2} \sum_{t \in \aleph^k} \sum_{r=1}^{q-1} a^r_{i}({S'}_{[t]}) X^r({S'}_{[t]}) - b_i \bigg),
\end{equation}
where ${\mathcal S}^{E,k}$ is a multi-set of $\eta_1$ independent draws from \textrm{SIM}(E) ($\eta_1$ complete trajectories conditional on $E$, independent across $E$, common across $X$), and $\aleph^k$ is a set of $\eta_2$ indices selected uniformly at random without replacement from $\{1, \dots, T\}$ (with the same set used across all $E$ and $X$, where we note that sampling with replacement would also work), with these random sets independent across $k$.  For $E \in {\mathcal E}$, let $\hat{G}^k(X)_{E}$ denote the $q$-dimensional vector with component $l$ equal to $\hat{G}^k(X)_{E,l}$.  
\\\indent We now define our family of SGA methods Algorithm \ref{alg:meta-gd}, where (to be clear) ${\mathcal X}^k(E)$ is a q-dimensional vector for $k \geq -1$ and $E \in {\mathcal E}$, and ${\mathcal X}^k$ is the corresponding $q \times |{\mathcal E}|$-dimensional vector.
\vspace{-.45in}
\LinesNotNumbered
\begin{algorithm}[ht]
\caption{\textup{\textsf{Stochastic Gradient for \ref{pen-smooth}}}}
\label{alg:meta-gd}
\vspace{2mm}
   Initialize ${\mathcal X}^{-1}(E), {\mathcal X}^{0}(E) \gets \overline{e}^q\ \forall\ E \in \mathcal{E}$ 
\\For $k\gets\ 0\ \KwTo\ K\ \textrm{and}\ E \in \mathcal{E}$
\\\ \ \ ${\mathcal X}^{k+1}(E) = \Pi_{{\mathcal Q}^0}\Bigg( (1 + \beta_k) {\mathcal X}^{k}(E) - \beta_k {\mathcal X}^{k-1}(E) + \alpha \hat{G}^{k}\big( (1 + \beta_k) {\mathcal X}^{k} - \beta_k {\mathcal X}^{k-1} \big)_E \Bigg)$ 
\end{algorithm}
\par \vspace{-.35in}
Using standard results and arguments from the literature on gradient methods in convex optimization, in particular \cite{schmidtetal2011} and \cite{nemirovski2012}, we derive the following convergence result for Algorithm \ref{alg:meta-gd}.  We defer the proof to Appendix\ \ref{appconvergence1sec}.
\begin{theorem}\label{convclaim1cor}
Suppose $\beta_k = 0$ for all $k \geq 0$, in which case $\textrm{Algorithm \ref{alg:meta-gd}}$ corresponds to projected SGA (as analyzed in \cite{nemirovski2012}).  If $\alpha = \frac{\iota^2 \epsilon}{24 L^2 q}, K = \lceil \frac{288 L^2 q}{\epsilon^2 \iota^2} \rceil, \eta_1 \geq \frac{2304 L^2 q}{\iota^2 \epsilon^2}, \eta_2 \geq \min\big(\frac{20736 L^2 T^2 q}{\iota^2 \theta^2 \epsilon^2} , T \big)$, then $\textrm{OPT}_{\textsc{pen}^{\theta}} - \mathbb{E}\big[f^{\theta}(K^{-1} \sum_{j=1}^K {\mathcal X}^{j})\big] \leq \epsilon T$.  
\\\indent Alternatively, suppose $\beta_0 = 0$, and $\beta_k = \frac{k-1}{k+2}$ for all $k \geq 1$, in which case $\textrm{Algorithm \ref{alg:meta-gd}}$ corresponds to the accelerated proximal-gradient method of \cite{schmidtetal2011}.  If $\alpha = \frac{1}{4} \lceil \frac{ (U_r L W)^{\frac{1}{4}}}{\sqrt{\iota \theta}} \rceil^{-2}, K = 14 \lceil \frac{ (U_r L W)^{\frac{1}{4}}}{\sqrt{\iota \theta}} \rceil \lceil \epsilon^{-\frac{1}{2}} \rceil, \eta_1 \geq 142800 \frac{L^2 q}{\iota^2 \epsilon^2}, \eta_2 \geq \min\big( 691200 \frac{L^2 T^2 q}{\iota^2 \theta^2 \epsilon^2}, T\big)$, then $\textrm{OPT}_{\textsc{pen}^{\theta}} - \mathbb{E}\big[f^{\theta}({\mathcal X}^{K})\big] \leq \epsilon T$.  
\end{theorem}
We now define a subroutine $\textrm{R}$ to compute ${\mathcal X}^k(E)$ more efficiently.  $\textrm{R}$ generates the multiset ${\mathcal S}^{E,k}$ ``on-the-fly", and utilizes a memoization table (i.e. a matrix $\Upsilon$ with $|\mathcal{E}|$ rows and an infinite number of columns) to consistently implement the required calculations.  $\Upsilon(E, 0), \Upsilon(E, -1)$ are initialized to $\overline{e}^q\ \forall\ E$ and all other entries are initialized to a dummy value $\emptyset.$  $\Upsilon(E, k)$ stores the (q-dimensional vector) value of ${\mathcal X}^k(E)$ (if computed).  $\textrm{R}$ can be used by the $\textrm{DM}$ to solve \ref{pen-smooth} on-the-fly by setting a desired number of iterations $K$, then calling $\textrm{R}(M_{[t]},K)$ and following the decision $K^{-1} \sum_{j=1}^K \Upsilon(M_{[t]},j)$ (in the unaccelerated setting) or $\Upsilon(M_{[t]},K)$ (in the accelerated setting) at each time $t$.  We may view $\Upsilon(\cdot, k)$ as a mapping from ${\mathcal E}$ to ${\mathcal R}^q \bigcup \emptyset$, and denote this mapping (in vector form) by $\Upsilon^k$.  We formally define $\textrm{R}$ in Algorithm \ref{alg:gd-implement}, and note that for $E \in {\mathcal E}$ and $S' \in {\mathcal S}^{E,k-1}$, $\aleph^{k-1} \bigcap \bigcup_{i \in a^+(E)} \mathcal{T}_i(S')$ represents the set of times $t$ for which $S'_{[t]}$ directly manifests in the calculation of $\hat{G}^{k-1}\big( (1 + \beta_{k-1}) \Upsilon^{k-1} - \beta_{k-1} \Upsilon^{k-2})_E$.  Thus $\bigcup_{S' \in {\mathcal S}^{E,k-1}} \big( \aleph^{k-1} \cap \bigcup_{i \in a^+(E)} \mathcal{T}_i(S') \big)$ corresponds to the set of necessary direct recursive calls in the calculation of $\Upsilon(E,k)$ by $\textrm{R}$ (modulo one more call to compute $\Upsilon(E,k-1)$ and $\Upsilon(E,k-2)$ for the given $E$, where we note that a call to $\textrm{R}(E,k-1)$ actually results in the computation of $\Upsilon(E,j)\ \forall\ j \leq k - 1$ as we show in Claim\ \ref{algclaim1} ).
\LinesNotNumbered
\vspace{-.25in}
\begin{algorithm}[ht]
\caption{$\textrm{R}(E,k)$}
\label{alg:gd-implement}
\DontPrintSemicolon
 If $\Upsilon(E, k) = \emptyset$
\\\ \ \ \ Call $\textrm{SIM}(E)$ $\eta_1$ times, store output trajectories in $\mathcal{S}^{E,k-1}$
\\\ \ \ \ If $\Upsilon(E, k-1) = \emptyset$\ then\ Call $\textrm{R}(E,k-1)$
\\\ \ \ \ For $S' \in \mathcal{S}^{E,k-1}\ \textrm{and}\ t \in \aleph^{k-1} \cap \bigcup_{i \in a^+(E)} \mathcal{T}_i(S')$
\\\ \ \ \ \ \ \ \ If $\Upsilon(S'_{[t]}, k-1) = \emptyset$\ then\ Call $\textrm{R}(S'_{[t]},k-1)$
{\small
\\\ \ \ \ $\Upsilon(E,k) = \Pi_{{\mathcal Q}^0}\bigg( (1 + \beta_{k-1}) \Upsilon(E,k-1) - \beta_{k-1} \Upsilon(E,k-2) + \alpha \hat{G}^{k-1}\big( (1 + \beta_{k-1}) \Upsilon^{k-1} - \beta_{k-1} \Upsilon^{k-2}\big)_E \bigg)$}
\end{algorithm}
\par\vspace{-.35in}
We now state the fact that a call to $\textrm{R}(E,k)$ results in $\Upsilon(E,j)$ being assigned value ${\mathcal X}^j(E)\ \forall\ j \in \lbrace -1,\ldots,k \rbrace.$  We defer the proof to Appendix \ref{APPalgclaim1sec}. 

\begin{claim}\label{algclaim1}
Every call to $\textrm{R}(E,k)$ terminates in finite time, and upon termination $\Upsilon(E,j)$ will have permanently been assigned value ${\mathcal X}^j(E)$ for all $j \in \lbrace -1,\ldots,k \rbrace$.
\end{claim}

We now analyze the runtime of $\textrm{R}(E,k)$.  Let $\textrm{COMPLEXITY}(k)$ denote the supremum, over all $E \in {\mathcal E}$, of the computational complexity (including the time for all simulations and recursive calls) to execute a call to $\textrm{R}(E,k)$.  We prove the following, deferring the proof to Appendix\ \ref{APPcomplexity1bsec}.  
\begin{lemma}\label{complexity1b}
$\textrm{COMPLEXITY}(k) \leq {\mathcal C} L q^2 (c_{comp} \eta_1 \eta_2)^{k+2}$; and if $\eta_2 = T$ (i.e. there is no subsampling) then $\textrm{COMPLEXITY}(k) \leq {\mathcal C} (c_{comp} \eta_1 U L q)^{k+2}$, for $c_{\textrm{comp}}$ an absolute constant.
\end{lemma}

We now complete the proof of Theorem\ \ref{convergence1}.  We define $\textrm{A}_{\textsc{pen}^{\theta}}$ to be the policy that (at each time $t$) calls $R(S_{[t]},K)$ (with the appropriate $K$ and other parameters as dictated by Theorem\ \ref{convclaim1cor}), and returns either $K^{-1}\sum_{j=1}^K \Upsilon(E,j)$ (in the unaccelerated case) or $\Upsilon(E,K)$ (in the accelerated case).  These cases (accelerated vs. unaccelerated) are delineated by the momentum and sampling parameters used to implement $\textrm{A}_{\textsc{pen}^{\theta}}$, and result in the three-piece complexity term appearing in Theorem\ \ref{convergence1} (each piece coming from a different way to set these parameters).

\begin{proof}{Proof of Theorem\ \ref{convergence1} :}
The result follows by combining Theorem\ \ref{convclaim1cor} with Claim\ \ref{algclaim1} and Lemma\ \ref{complexity1b} in the natural manner, and we omit the details.  $\qed$
\end{proof}

\section{On-the-fly algorithm for non-smooth penalty formulation}\label{penpensec1}
Let $f : {\mathcal R}^{q \times |{\mathcal E}|} \rightarrow {\mathcal R}$ denote the mapping
$$f(X) \stackrel{\Delta}{=} \sum_{E \in {\mathcal E}} \mu(E) \sum_{r=1}^{q-1} Z^r(E) X^r(E) - 2 \iota^{-1} \sum_{S \in {\mathcal S}} \mu(S) \sum_{i=1}^m \big( \sum_{t=1}^T \sum_{r=1}^{q-1} a^r_i(S_{[t]}) X^r(S_{[t]}) - b_i \big)^+.$$  Let $\textsc{pen}$ denote the following concave program, with optimal value $\textrm{OPT}_{\textsc{pen}}$ and solution $X_{*,\textsc{pen}}$.
\vspace{-.15in}
\begin{align}
    \tag{$\textsc{pen}$} \label{pen}
    \max f(X) \quad s.t. \ X \in {\mathcal Q}^2
\end{align}
\par \vspace{-.15in}
We use our results for \ref{pen-smooth} to prove an analogous result for \ref{pen}, in which we also control the feasibility violation (relative to $\textsc{lp}$ in which the constraints are enforced w.p.1).  Our main result for \ref{pen} (below) utilizes policy $\textrm{A}_{\textsc{pen}}$ which is simply $\textrm{A}_{\textsc{pen}^{\theta}}$ run with parameters $\theta = \frac{\epsilon T \iota}{4 V}$ and $\epsilon' = \frac{\epsilon}{2}$.

\begin{theorem}\label{NRMAIN3}
$\forall\ \epsilon \in (0, 1)$, on any trajectory $S \in \mathcal{S}$ one can compute decisions $\textrm{A}_{{\textsc{pen}}}(S_{[t]})\ \forall\ t = 1, \dots, T$ on-the-fly, in per-decision time
$$
{\mathcal C} \times \min\Bigg( \big( c_3 \frac{L V q}{\iota \epsilon} \big)^{c_3 (\frac{L}{\iota \epsilon})^2 q} , \big( c_3 \frac{L V q}{\iota \epsilon} \big)^{c_3 \sqrt{\frac{L \kappa V}{\iota^2 \epsilon^2}}} , \big( c_3 \frac{L U q}{\iota \epsilon} \big)^{c_3 \lceil \frac{ (U_r L W)^{\frac{1}{4}} \sqrt{V}}{\iota \sqrt{\epsilon T}} \rceil \epsilon^{-\frac{1}{2}}} \Bigg),
$$
where $c_3$ is an absolute constant.  Also, $\mathbb{E}\big[ \sum_{E \in {\mathcal E}} \mu(E) \sum_{r=1}^{q-1} Z^r(E) \textrm{A}_{\textsc{pen}}^r(E) \big] \geq \textrm{OPT}_{\textsc{pen}} - \epsilon T;$
and  $\mathbb{E}\big[ \sum_{S \in {\mathcal S}} \mu(S) \sum_{i=1}^m \big( \sum_{t=1}^T \sum_{r=1}^{q-1} a^r_i(S_{[t]}) \textrm{A}_{\textsc{pen}}^r(S_{[t]}) - b_i \big)^+ \big] \leq \iota \epsilon T.$
\end{theorem} 
To prove Theorem\ \ref{NRMAIN3}, we proceed as follows.  First, we show that any approximately optimal solution to \ref{pen-smooth} is also approximately optimal to \ref{pen}, and defer the proof to Appendix\ \ref{NRMAIN3citenoweq1proofsec}.
\begin{claim}\label{NRMAIN3citenoweq1}
$\forall\ X \in {\mathcal Q}^2$, $\textrm{OPT}_{\textsc{pen}}  - f(X) \leq \textrm{OPT}_{\textsc{pen}^{\theta}} - f^{\theta}(X) + 2 \iota^{-1} V \theta.$
\end{claim}
Second, we prove that any solution (of \ref{pen}) with large aggregate feasibility violation must be suboptimal in value by a large margin.  
\begin{lemma}\label{NRMAIN3proof200}
$\forall X \in {\mathcal Q}^2, \sum_{S \in {\mathcal S}} \mu(S) \sum_{i=1}^m \bigg( \sum_{t=1}^T \sum_{r=1}^{q-1} a^r_i(S_{[t]}) X^r(S_{[t]}) - b_i \bigg)^+ \leq \iota\big(\textrm{OPT}_{\textsc{pen}} - f(X)\big).$
\end{lemma}
\ \indent We prove Lemma\ \ref{NRMAIN3proof200} by showing that given any $X \in {\mathcal Q}^2$ for \ref{pen}, one can construct a significantly better solution whenever $X$ has large aggregate feasibility violation.  We will conclude that any approximately optimal solution for \ref{pen} cannot have large aggregate feasibility violation.  More formally, for $X \in {\mathcal Q}^2$, we now define a mapping (from ${\mathcal Q}^2$ to ${\mathcal Q}^2$) $\textrm{F}(X)$ which is admissible (i.e. feasible) for $\textsc{lp}$, and corresponds to the policy which implements $X$, except whenever a constraint would be violated the corresponding values are ``pushed" to option 0 to maintain feasibility.  We define 
$\textrm{F}(X)$ using forward induction. For $E$ a $D \times 1$ matrix, $\textrm{F}(X)^1(E) = \min\bigg( X^1(E), \min_{i \in a^{1,+}(E)} \frac{b_i}{a^1_i(E)} \bigg).$  Supposing (for some $t \geq 1$ and $r \in \lbrace 1,\ldots q-1 \rbrace$) we have defined $\textrm{F}(X)^l(E)$ for 
all $l \in \lbrace 1,\ldots,q-1 \rbrace$ when $E$ is a $D \times u$ matrix with $u \in \lbrace 1,\ldots,t-1 \rbrace$, and for all $l \in \lbrace 1,\ldots,r-1 \rbrace$ when $E$ is a $D \times t$ matrix, we define
$\textrm{F}(X)^r(E)$ for $E$ a $D \times t$ matrix as
$$\textrm{F}(X)^r(E) = \min\bigg( X^r(E), \min_{i \in a^{r,+}(E)} \frac{b_i - \sum_{u=1}^{t-1} \sum_{l=1}^{q-1} a^l_i(E_{[u]}) \textrm{F}(X)^l(E_{[u]}) - \sum_{l=1}^{r - 1} a^l_i(E) \textrm{F}(X)^l(E)}{a^r_i(E)} \bigg).$$
We note that when $t = 1$, $\sum_{u=1}^{t-1} \sum_{l=1}^{q-1} a^l_i(E_{[u]}) \textrm{F}(X)^l(E_{[u]}) = 0$; and when $r = 1$, $\sum_{l=1}^{r - 1} a^l_i(E) \textrm{F}(X)^l(E) = 0$.  If $|a^{r,+}(E)| = 0$, we instead set $\textrm{F}(X)^r(E) = X^r(E)$.  We set $\textrm{F}(X)^0(E) = 1 - \sum_{r=1}^{q-1} \textrm{F}(X)^r(E)$.  Note that $\textrm{F}(X)$ is feasible for $\textsc{lp}$.  We now quantify the reward lost by $\textrm{F}(X)$ (relative to $X$), deferring the proof to Appendix \ \ref{APPcomparefeasible1sec}.
\begin{lemma}\label{comparefeasible1}
$\forall\ X \in {\mathcal Q}^2, \sum_{E \in {\mathcal E}} \mu(E) \sum_{r=1}^{q-1} Z^r(E) \textrm{F}(X)^r(E)$ is at least  
$$\sum_{E \in {\mathcal E}} \mu(E) \sum_{r=1}^{q-1} Z^r(E) X^r(E) -  \iota^{-1} \sum_{S \in {\mathcal S}} \mu(S) \sum_{i=1}^m \big( \sum_{t=1}^T \sum_{r=1}^{q-1} a^r_i(S_{[t]}) X^r(S_{[t]}) - b_i \big)^+.$$
\end{lemma}
\ \indent We now combine Lemma\ \ref{comparefeasible1} with the fact that reducing the feasibility violation significantly increases the objective of \ref{pen} (due to the $2 \iota^{-1}$ multiplier in $f$) to complete the proof of Lemma\ \ref{NRMAIN3proof200}. 
\begin{proof}{Proof of Lemma\ \ref{NRMAIN3proof200} :}Combining Lemma\ \ref{comparefeasible1} with the feasibility (for $\textsc{lp}$) of $\textrm{F}(X)$ and definition of $f$ implies $f\big( \textrm{F}(X) \big) - f(X)$ equals 
{\small
$$\sum_{E \in {\mathcal E}} \mu(E) \sum_{r=1}^{q-1} Z^r(E) \textrm{F}(X)^r(E) - \bigg( \sum_{E \in {\mathcal E}} \mu(E) \sum_{r=1}^{q-1} Z^r(E) X^r(E) - 2 \iota^{-1} 
\sum_{S \in {\mathcal S}} \mu(S) \sum_{i=1}^m \big( \sum_{t=1}^T \sum_{r=1}^{q-1} a^r_i(S_{[t]}) X^r(S_{[t]}) - b_i \big)^+ \bigg),$$} itself at least 
$$\iota^{-1} 
\sum_{S \in {\mathcal S}} \mu(S) \sum_{i=1}^m \big( \sum_{t=1}^T \sum_{r=1}^{q-1} a^r_i(S_{[t]}) X^r(S_{[t]}) - b_i \big)^+.$$  As $\textrm{F}(X)$ is also feasible for \ref{pen}, it follows that $$\textrm{OPT}_{\textsc{pen}} - f(X) \geq f\big( \textrm{F}(X) \big) - f(X) \geq 
\iota^{-1} 
\sum_{S \in {\mathcal S}} \mu(S) \sum_{i=1}^m \big( \sum_{t=1}^T \sum_{r=1}^{q-1} a^r_i(S_{[t]}) X^r(S_{[t]}) - b_i \big)^+.\ \ \qed$$
\end{proof}
\ \indent Finally, we combine all of the above to complete the proof of Theorem\ \ref{NRMAIN3}.  

\begin{proof}{Proof of Theorem\ \ref{NRMAIN3} :}
The first part of the proof follows by applying Claim\ \ref{NRMAIN3citenoweq1} with $X = k^{-1} \sum_{j=1}^k {\mathcal X}^j$ (in the unaccelerated case) or $X = {\mathcal X}^k$ (in the accelerated case) for appropriate $k$, taking expectations, and combining (in the natural manner) with Theorem\ \ref{convergence1} (applied with $\theta = \frac{\epsilon T \iota}{4 V}$ and replacing $\epsilon$ with $\frac{\epsilon}{2}$).  For the second of the three complexity terms appearing in the minimum, we also use the bounds $U_r \leq \kappa T, W \leq L \kappa T.$  The second part of the proof then follows directly from Lemma\ \ref{NRMAIN3proof200} after taking expectations on both sides.  $\qed$
\end{proof}

\section{On-the-fly algorithm for $\textsc{lp}$}\label{lpsec}
We now state and prove our main result for $\textsc{lp}$, which utilizes policy $\textrm{A}_{\textsc{lp}}$ (which is simply the feasibility-patching procedure $\textrm{F}$ composed with $\textrm{A}_{\textsc{pen}}$ run with $\epsilon' = \frac{1}{2} \epsilon$).

\begin{theorem}\label{NRMAIN2}
$\forall\ \epsilon \in (0, 1)$, on any trajectory $S \in \mathcal{S}$ one can compute decisions $\textrm{A}_{\textsc{lp}}(S_{[t]})\ \forall\ t = 1, \dots, T$ on-the-fly, in per-decision time $$
{\mathcal C} \times \min\Bigg( \big( c_1 \frac{L V q}{\iota \epsilon} \big)^{c_1 (\frac{L}{\iota \epsilon})^2 q} , \big( c_1 \frac{L V q}{\iota \epsilon} \big)^{c_1 \sqrt{\frac{L \kappa V}{\iota^2 \epsilon^2}}} , \big( c_1 \frac{L U q}{\iota \epsilon} \big)^{c_1 \lceil \frac{ (U_r L W)^{\frac{1}{4}} \sqrt{V}}{\iota \sqrt{\epsilon T}} \rceil \epsilon^{-\frac{1}{2}}} \Bigg),
$$
where $c_{1}$ is an absolute constant. Also, $\mathbb{E}\big[ \sum_{E \in {\mathcal E}} \mu(E) \sum_{r=1}^{q-1} Z^r(E) \textrm{A}_{\textsc{lp}}^r(E) \big] \geq \textrm{OPT}_{\textsc{lp}} - \epsilon T.$
\end{theorem}
We proceed by patching the infeasibility in our solution to \ref{pen} (combining $\textrm{F}$ and Theorem\ \ref{NRMAIN3}).
\begin{proof}{Proof of Theorem\ \ref{NRMAIN2} :} Consider $\textrm{F}(\textrm{A}_{\textsc{pen}})$ implemented with $\epsilon' = \frac{1}{2} \epsilon$.  Theorem\ \ref{NRMAIN3} implies
{\small 
$$\mathbb{E}\big[ \sum_{E \in {\mathcal E}} \mu(E) \sum_{r=1}^{q-1} Z^r(E) \textrm{A}_{\textsc{pen}}^r(E) \big] \geq \textrm{OPT}_{\textsc{pen}} - \frac{1}{2} \epsilon T\ \ ,\ \ \textrm{and}\ \ \mathbb{E}\big[ \sum_{S \in {\mathcal S}} \mu(S) \sum_{i=1}^m \big( \sum_{t=1}^T \sum_{r=1}^{q-1} a^r_i(S_{[t]}) \textrm{A}_{\textsc{pen}}^r(S_{[t]}) - b_i \big)^+ \big] \leq \frac{1}{2} \iota \epsilon T.$$
}
Combining with Lemma\ \ref{comparefeasible1}, the feasibility (for $\textsc{lp}$) of $\textrm{F}(\textrm{A}_{\textsc{pen}})$, and the fact that $\textrm{OPT}_{\textsc{pen}} \geq \textrm{OPT}_{\textsc{lp}}$ then yields the desired performance guarantee.  As we prove in Appendix\ \ref{compfsec} that $\textrm{F}$ may be implemented on-the-fly in $O(L \times q)$ time, combining with Theorem\ \ref{NRMAIN3} completes the proof.  $\qed$
\end{proof}

\section{Proof of main result for $\textrm{NRM}$ Theorem\ \ref{NRMAIN1}}\label{nrmain1sec1}
We now complete the proof of Theorem\ \ref{NRMAIN1}.  Given $X \in {\mathcal Q}^2$, we define a (random) mapping $\textrm{H}(X)$ from ${\mathcal Q}^2$ to ${\mathcal Q}^2$, corresponding to the natural dependent rounding scheme where for each $E$ we select exactly one r $\in \lbrace 0,\ldots,q-1 \rbrace$ as the chosen option to implement.  More formally, independently for each $E$, an $\hat{r}(X,E) \in \lbrace 0,\ldots,q-1 \rbrace$ is selected randomly s.t. $\mathbb{P}\big(\hat{r}(X,E) = l\big) = X^l(E)$ (for $l = 0,\ldots,q-1$), and $\textrm{H}(X)^r(E)$ is set equal to $I\big(r = \hat{r}(X,E)\big)$ for $r \in \lbrace 0,\ldots,q-1 \rbrace$.
We now prove that $\textrm{H}(X)$ will achieve the same expected reward as $X$, but not have too much more (expected) aggregate inequality violation, deferring the proof to Appendix\ \ref{APPnrmroundlemma1sec}.

\begin{lemma}\label{nrmroundlemma1}
$\forall\ X \in {\mathcal Q}^2$, 
$\mathbb{E}\big[ \sum_{E \in {\mathcal E}} \mu(E) \sum_{r=1}^{q-1} Z^r(E) \textrm{H}(X)^r(E) \big] = \sum_{E \in {\mathcal E}} \mu(E) \sum_{r=1}^{q-1} Z^r(E) X^r(E)$, and
$\mathbb{E}\big[ \sum_{S \in {\mathcal S}} \mu(S) \sum_{i=1}^m \big( \sum_{t=1}^T \sum_{r=1}^{q-1} a^r_i(S_{[t]}) \textrm{H}(X)^r(S_{[t]}) - b_i \big)^+ \big]$ is at most  
$$\sum_{S \in {\mathcal S}} \mu(S) \sum_{i=1}^m \big( \sum_{t=1}^T \sum_{r=1}^{q-1} a^r_i(S_{[t]}) X^r(S_{[t]}) - b_i \big)^+ +  \sqrt{\frac{\pi}{2}} \min\big( \sqrt{L q m} \sqrt{T} , V \sqrt{T} + m T \exp( - \frac{\nu^2}{2} T ) \big).
$$
\end{lemma}

Next we prove that although $\textrm{F}\big( \textrm{H}( \textrm{A}_{\textsc{pen}} ) \big)$ may be non-integral even though $\textrm{H}( \textrm{A}_{\textsc{pen}})$ is integral, the effect is mild, and also that $V \leq \lambda \stackrel{\Delta}{=} \min(m, \frac{2L}{\nu})$, deferring the proof to Appendix\ \ref{APPnrmroundlemma2sec}.

\begin{lemma}\label{nrmroundlemma2}
Given $X \in {\mathcal Q}^2$ s.t. $X^r(E) \in \lbrace 0,1 \rbrace$ for all $r$ and $E$, it holds that for any $S \in {\mathcal S}$, $\lbrace \textrm{F}(X)^r(S_{[t]}), r = 0,\ldots,q-1; t = 1,\ldots,T \rbrace$ has at most $V$ non-integer values.  Also, $V \leq \lambda$.
\end{lemma}

We now complete the proof of Theorem\ \ref{NRMAIN1}.  Given $X \in {\mathcal Q}^2$, let $\textrm{J}(X)$ denote the policy that sets $\textrm{J}(X)^r(E) = X^r(E)\ \forall\ r$ if $X^r(E) \in \lbrace 0,1 \rbrace\ \forall\ r$, and otherwise sets $X^0(E) = 1, X^r(E) = 0$ for $r \neq 0$ (a ``round down" operation, to handle the mild non-integrality formalized in Lemma\ \ref{nrmroundlemma2}).

We define $\textrm{A}_{\textsc{nrm}}$ to equal $\textrm{J}\bigg( \textrm{F} \big( \textrm{H}( \textrm{A}_{\textsc{pen}}) \big) \bigg)$ implemented with $\epsilon' = .45 \epsilon$, and now complete the proof of Theorem\ \ref{NRMAIN1}.  We actually prove the below stronger result (as alluded to in Section\ \ref{mainresultssec}), from which Theorem\ \ref{NRMAIN1} follows after some straightforward algebra.  Here $c'_0$ is an absolute constant.

\begin{lemma}\label{NRMAIN1lemlem}
$\forall\ \epsilon \in (0, 1)$, $\textrm{A}_{\textsc{nrm}}$ is admissible for $\textsc{nrm}$, and on any trajectory $S \in \mathcal{S}$ one can compute decisions $\textrm{A}_{\textsc{nrm}}(S_{[t]})\ \forall\ t = 1, \dots, T$ on-the-fly, in per-decision time ${\mathcal C} \times \min\bigg(\big( c'_0 \frac{ L \lambda q}{\iota \epsilon} \big)^{c'_0 \frac{L^2 q}{\iota^2} \epsilon^{-2} } , \big( c'_0 \frac{L \lambda q}{\iota \epsilon} \big)^{c'_0 \sqrt{\frac{L \kappa \lambda}{\iota^2}} \epsilon^{-1}} \bigg).$  Also, $\mathbb{E}\big[ \sum_{E \in {\mathcal E}} \mu(E) \sum_{r=1}^{q-1} Z^r(E) \textrm{A}_{\textsc{nrm}}^r(E) \big] \geq \textrm{OPT}_{\textsc{nrm}} - \epsilon T,$
so long as $T \geq c'_0 \min\big( \frac{ m L q}{\iota^2 \epsilon^2} , \frac{\lambda^2}{\iota^2 \epsilon^2} + \frac{\log(1 + \frac{m}{\iota \epsilon})}{\nu^2} \big).$
\end{lemma}

\begin{proof}{Proof of Lemma\ \ref{NRMAIN1lemlem} and Theorem\ \ref{NRMAIN1} :}
Consider $\textrm{A}_{\textsc{pen}}$ with $\epsilon' = .45 \epsilon$.  Theorem\ \ref{NRMAIN3} implies that 
{\small
$$\mathbb{E}\big[ \sum_{E \in {\mathcal E}} \mu(E) \sum_{r=1}^{q-1} Z^r(E) \textrm{A}_{\textsc{pen}}^r(E) \big] \geq \textrm{OPT}_{\textsc{pen}} - .45 \epsilon T\ ,\ 
\mathbb{E}\bigg[ \sum_{S \in {\mathcal S}} \mu(S) \sum_{i=1}^m \big( \sum_{t=1}^T \sum_{r=1}^{q-1} a^r_i(S_{[t]}) \textrm{A}_{\textsc{pen}}^r(S_{[t]}) - b_i \big)^+ \bigg] \leq .45 \iota \epsilon T.$$
}
Combining with Lemma\ \ref{nrmroundlemma1} and our assumptions on $T$ (by taking $c_0$ sufficiently large) implies
{\small
$$\mathbb{E}\big[ \sum_{E \in {\mathcal E}} \mu(E) \sum_{r=1}^{q-1} Z^r(E) \textrm{H}\big(\textrm{A}_{\textsc{pen}}\big)^r(E) \big] \geq \textrm{OPT}_{\textsc{pen}} - .45 \epsilon T\ ,\ 
\mathbb{E}\bigg[ \sum_{S \in {\mathcal S}} \mu(S) \sum_{i=1}^m \big( \sum_{t=1}^T \sum_{r=1}^{q-1} a^r_i(S_{[t]}) \textrm{H}\big(\textrm{A}_{\textsc{pen}}\big)^r(S_{[t]}) - b_i \big)^+ \bigg] \leq .5 \iota \epsilon T.$$}
Combining with Lemma\ \ref{comparefeasible1}, the fact that $\textrm{OPT}_{\textsc{pen}} \geq \textrm{OPT}_{\textsc{lp}}$, and Lemma\ \ref{nrmroundlemma2} implies
$$\mathbb{E}\big[ \sum_{E \in {\mathcal E}} \mu(E) Z(E) \textrm{J}\big(\textrm{F}\big( \textrm{H}(\textrm{A}_{\textsc{pen}}) \big)\big) (E) \big] \geq \textrm{OPT}_{\textsc{lp}} - .95 \epsilon T - \lambda.$$  As our assumptions on $T$ also imply $\lambda \leq .05 \epsilon T$, we conclude the performance guarantee.  Our proof of Theorem\ \ref{NRMAIN2} implies $\textrm{F}$ can be implemented in $O(L \times q)$ time, while $\textrm{H}$ and $\textrm{J}$ may be implemented in $O(q)$ time.  Combining with Theorem\ \ref{NRMAIN3} completes the proof.  $\qed$
\end{proof}

\section{Proof of main result for $\textsc{is}$ (part 1 of Theorem\ \ref{NRMAINCOMBOPT})}\label{issec1}
As explained in Section\ \ref{applicationssec}, $\textsc{is}$ can be put in the $\ref{pack}$ framework, with $q = 2, T = n, m = \lfloor \frac{1}{2} \Delta n \rfloor$; and it is easily verified that one can take $L = \Delta, U = 2, U_r = 2, W = 2 \Delta, V = \lfloor \frac{1}{2} \Delta n \rfloor, \iota = 1.$  Let us denote the LP relaxation of $\textsc{is}$ (in which all constraints of the form $X^r(E) \in \{0,1\}$ are replaced by $X^r(E) \in [0,1]$) by \textsc{lp-is} (with optimal value $\textrm{OPT}_{\textsc{lp-is}}$).  Then Theorem\ \ref{NRMAIN2} directly implies the natural analogue of part 1 of Theorem\ \ref{NRMAINCOMBOPT} for \textsc{lp-is} (with the same per-decision complexity and guarantee, and $\textrm{OPT}_{\textsc{is}}$ replaced by $\textrm{OPT}_{\textsc{lp-is}}$).  We denote the corresponding policy by $\textrm{A}_{\textsc{lp-is}}$.  We now use this fact to complete the proof of Theorem\ \ref{NRMAINCOMBOPT} for $\textsc{is}$ by combining with an appropriate rounding scheme.
\begin{proof}{Proof of Theorem\ \ref{NRMAINCOMBOPT} for $\textsc{is}$ :}
We first define $\textrm{A}_{\textsc{is}}$ using $\textrm{A}_{\textsc{lp-is}}$ and a rounding scheme.  Let ${\mathcal U}$ be a fixed $\textrm{Uniform}[0,1]$ r.v., generated once at time zero.  Then we define $\textrm{A}^1_{\textsc{is}}(E)$ to equal $I\big( \textrm{A}^1_{\textsc{lp-is}}(E) > {\mathcal U} \big)$ if $E$ is a $D \times t$ matrix where node $t$ is in partite ${\mathcal L}$; and $\textrm{A}^1_{\textsc{is}}(E) = I\big( \textrm{A}^1_{\textsc{lp-is}}(E) > 1 - {\mathcal U} \big)$ if $E$ is a $D \times t$ matrix where node $t$ is in partite ${\mathcal R}$.  That $\mathbb{E}\big[ \sum_{E \in {\mathcal E}} \mu(E) Z(E) \textrm{A}^1_{\textsc{is}}(E) \big] \geq \textrm{OPT}_{\textsc{is}} - \epsilon n$ follows from our previous analogue of part 1 of Theorem\ \ref{NRMAINCOMBOPT} for \textsc{lp-is}.  Thus we need only verify that $\lbrace \textrm{A}_{\textsc{is}}(E), E \in {\mathcal E} \rbrace$ is w.p.1 feasible for $\textsc{is}$.  Suppose for contradiction it is not.  Then there must exist $S \in {\mathcal S}$, $t_L$ corresponding to a node in partite ${\mathcal L}$, and $t_R$ corresponding to a node in partite ${\mathcal R}$, s.t. for some potential edge $i$ it holds that $a^1_i(S_{[t_L]}) = 1$ and $a^1_i(S_{[t_R]}) = 1$, but also $\textrm{A}^1_{\textsc{is}}(S_{[t_L]}) = 1$ and $\textrm{A}^1_{\textsc{is}}(S_{[t_R]}) = 1$.  $\textrm{A}^1_{\textsc{is}}(S_{[t_L]}) = 1$ and $\textrm{A}^1_{\textsc{is}}(S_{[t_R]}) = 1$ implies $\textrm{A}^1_{\textsc{lp-is}}(S_{[t_L]}) > {\mathcal U}$ and $\textrm{A}^1_{\textsc{lp-is}}(S_{[t_R]}) > 1 - {\mathcal U},$ implying $\textrm{A}^1_{\textsc{lp-is}}(S_{[t_L]}) + \textrm{A}^1_{\textsc{lp-is}}(S_{[t_R]}) > 1$.  But $a^1_i(S_{[t_L]}) = 1, a^1_i(S_{[t_R]}) = 1$, and the feasibility of $\textrm{A}^1_{\textsc{lp-is}}$ for $\textsc{lp-is}$, implies $\textrm{A}^1_{\textsc{lp-is}}(S_{[t_L]}) + \textrm{A}^1_{\textsc{lp-is}}(S_{[t_R]}) \leq 1$, yielding a contradiction.  $\qed$
\end{proof}

\section{Proof of main result for $\textsc{mwmlp}$ and $\textsc{mmo}$ (parts 2-3 of Theorem\ \ref{NRMAINCOMBOPT})}\label{matchsec1}
Our proof of Theorem\ \ref{NRMAINCOMBOPT} for \textsc{mwmlp} follows from Theorem\ \ref{NRMAIN2}, and we defer the proof to Appendix\ \ref{MWMLPproofsec}.  Our proof of Theorem\ \ref{NRMAINCOMBOPT} for \textsc{mmo} combines Theorem\ \ref{NRMAIN2} with the result for rounding fractional matchings in \cite{harrisetal2026}.  First, we review the results of \cite{harrisetal2026}.  Suppose there is an unknown bipartite graph $G$ with $n$ nodes and maximum degree $\Delta$, of which partite ${\mathcal L}$ are offline nodes present at time 0, and partite ${\mathcal R}$ are online nodes which arrive sequentially (over $|{\mathcal R}|$ periods).  Upon arrival of a given node, all of its incident edges are revealed, along with (possibly fractional) values in $[0,1]$ (one for each incident edge).  For each edge $e$ in $G$, let $f_e$ denote the fractional value revealed for that edge.  For a given node $v$ in $G$ (offline or online), let $E_v$ denote the set of edges incident to $v$ in unknown graph $G$.  It is promised that $\sum_{e \in E_v} f_e \leq 1$ for all nodes $v$ in $G$, i.e. the fractional values revealed online (along with the graph) constitute a feasible fractional matching.  Then Theorem 1.2 of \cite{harrisetal2026} proves the following.  Although those authors do not explicitly prove a $O(\Delta)$ per-decision runtime, it is easily verified.
\begin{theorem}[\cite{harrisetal2026} Theorem 1.2]\label{naorresult}
$\forall\ \epsilon \in (0,1)$, $\exists$\ a (randomized) online algorithm $A_{\textrm{H}}$ that at each time randomly selects at most one of the edges incident to the online node revealed, s.t. : (1) the set $M$ of selected edges is a feasible matching in $G$ w.p.1; (2) \ $\forall\ $\ edges $e$ in $G$ it holds that $\mathbb{P}(e \in M) \geq .68 f_e;$ and (3) in each period, $A_{\textrm{H}}$ can be implemented in time $O(\Delta)$.
\end{theorem}
\begin{proof}{Proof of Theorem\ \ref{NRMAINCOMBOPT} for \textsc{mmo} :} The proof follows by : (1) observing that the LP relaxation of $\textsc{MMO}$ can be put in the framework of $\textsc{LP}$ (as outlined in Section\ \ref{applicationssec}) and applying Theorem\ \ref{NRMAIN2} with $T = |{\mathcal R}|, m = |{\mathcal L}|, q = \Delta + 1, L = 1, U = \Delta, U_r = \Delta, W = \Delta, V = |{\mathcal L}|, \iota = 1;$ (2) setting $f_e$ equal to $\textrm{A}^r_{\textsc{lp}}(E)$ for the $E$ and $r$ corresponding to the revealed edges; and (3) using Theorem\ \ref{naorresult} to round these fractional values online.  To complete the proof, we observe that (since $T = |{\mathcal R}|$, not $n$), to achieve error $\epsilon n$ (and as $n = |{\mathcal L}| + |{\mathcal R}|$) we apply Theorem\ \ref{NRMAIN2} with $\epsilon' = \min\big( \frac{ |{\mathcal L}| + |{\mathcal R}|}{|{\mathcal R}|} \epsilon, 1 \big)$.  Combining with the fact that $\frac{\sqrt{|{\mathcal L}| |{\mathcal R}|} }{|{\mathcal L}| + |{\mathcal R}|} \leq \frac{1}{2}$ (since $\frac{\sqrt{x y}}{x + y} \leq \frac{1}{2}\ \forall\ x,y > 0$) completes the proof.  We also note that if $\epsilon' = 1$, the $\textrm{DM}$ can simply return the empty matching.  $\qed$
\end{proof}

\section{Conclusion}\label{concsec1}
In this work, we derived algorithms for o.s.p./c.a. with general correlations (also network revenue management and online bipartite independent set and matching) whose runtime scales as the time to simulate a single sample path of the underlying stochastic process, multiplied by a constant independent of the time horizon or number of states of the o.i.p.  At the heart of our algorithms is a new way to conceptualize and implement SGA methods in a completely on-the-fly/recursive manner for the associated massive deterministic-equivalent LP, and a recognition that to solve such online problems one need only compute the values of very few variables in this exponentially large LP.  Our work leaves many directions for future work, including an exploration of : (1) the full range of problems to which our approach may be applied; (2) whether our approach can be made practical by combining with techniques from optimization and simulation; (3) lower bounds; (4) connections to other computational models such as parallel, distributed, local, and quantum; (5) the implications of such an ``$\textrm{efficient simulator}\ \rightarrow\ \textrm{efficient algorithm}$" result; and (6) the relationship between the online optimization and multistage stochastic programming literatures broadly.


\newpage
\section{Electronic Compendium : Technical Appendix}\label{techsec1}
\subsection{Proof of Theorem\ \ref{convclaim1cor}}\label{appconvergence1sec} 
We now use results from the literature on (accelerated) inexact gradient methods to analyze the convergence of Algorithm \ref{alg:meta-gd} and prove Theorem\ \ref{convclaim1cor}.  Recall from our discussion of gradients in Section\ \ref{introsec1} that component $(r,E)$ of the gradient (of $f^{\theta}$) is scaled by $\mu(E)$, while variable $X^r(E)$ itself is not.  There we claimed (loosely to build intuition) that we could implicitly account for this ``mismatch" using a weighted Euclidean norm.  Indeed, a natural framework for handling this mismatch and analyzing Algorithm \ref{alg:meta-gd} would be that of mirror descent with inexact gradients, with weighted Euclidean norm $||X||_W = \sqrt{ \sum_{E \in {\mathcal E}} \mu(E) \sum_{r=0}^{q-1} \big(X^r(E)\big)^2 }$ (for $X \in {\mathcal R}^{q \times |{\mathcal E}|}$).  However, to the best of our knowledge the precise type of inexactness required for our analysis in the accelerated case is not available in the literature at the full generality of mirror descent (although for very recent independent and concurrent progress we refer the reader to \cite{zhangjaillet2025}).  Such a result is, however, available for inexact accelerated gradient descent under the Euclidean norm (\cite{schmidtetal2011}).  Thankfully, the weighted Euclidean norm $||X||_W$ is sufficiently similar to the Euclidean norm that one can directly reduce the desired ``inexact mirror descent analysis" to a very similar analysis under the Euclidean norm (by implementing a few linear transformations) and then apply the results of \cite{schmidtetal2011}.  This is how we handle the mismatch (in scaling between gradients and variables) and prove Theorem\ \ref{convclaim1cor}.  For completeness, we make the relevant statement precise and include a proof.  In the unaccelerated case we directly apply a more general result of \cite{nemirovski2012}, and note that although a tighter result can be proven using an argument similar to that of \cite{hamoudietal2019}, for simplicity we use the results of \cite{nemirovski2012}.
\\\indent In summary, our proof of Theorem\ \ref{convclaim1cor} proceeds as follows.
\begin{itemize}
\item First, we define a linearly transformed (by probabilities $\lbrace \mu(E), E \in {\mathcal E} \rbrace$) problem and algorithm.
\item Second, we prove that the rate of convergence to optimality of the transformed algorithm on the transformed problem (in which there is no longer a mismatch between the scaling of the underlying variables and the scaling of the gradients) is identical to the rate of convergence to optimality of \textrm{Algorithm} \ref{alg:meta-gd} on problem \ref{pen-smooth} (all in the Euclidean norm).
\item Third, we combine results of \cite{nemirovski2012} (in the unaccelerated case) and \cite{schmidtetal2011} (in the accelerated case), with some additional analysis (of smoothness parameters etc.), to analyze the convergence of the transformed algorithm to optimality in the transformed problem, and then transfer the result to the convergence of \textrm{Algorithm} \ref{alg:meta-gd} on \ref{pen-smooth}.
\end{itemize}
\subsubsection{Linearly transformed problem and algorithm}
We begin by defining the aforementioned linearly transformed problem and algorithm, which are essentially the same problem and algorithm but working with transformed variables $X'^r(E) = X^r(E) \sqrt{\mu(E)}$ (which intuitively ensures the underlying variables and gradients have the same $\sqrt{\mu(E)}$ scaling).  Let ${\mathcal Q}^3$ denote the rescaling of ${\mathcal Q}^2$ in which the constraints corresponding to $X^r(E)$ are rescaled by $\sqrt{\mu(E)}$ : 
$${\mathcal Q}^3 \stackrel{\Delta}{=} \bigg\lbrace X^r(E), r = 0,\ldots,q-1, E \in {\mathcal E}\ \textrm{s.t.}\ X^r(E) \in [0,\sqrt{\mu(E)}]\ \forall\ r,E ; \sum_{r=0}^{q-1} X^r(E) = \sqrt{\mu(E)}\ \forall E \bigg\rbrace.$$  Given $a \leq 0, b \geq 0$, we also define a family of polytopes which instead enforce that a certain rescaling (again by $\sqrt{\mu(E)}$) of the amount of each resource consumed (collectively by all options at any given time) lies in an interval defined by $a$ and $b$ : 
$${\mathcal Q}^4(a,b) \stackrel{\Delta}{=} \bigg\lbrace X^r(E), r = 0,\ldots,q-1, E \in {\mathcal E}\ \textrm{s.t.}\ \sum_{r=1}^{q-1} a^r_i(E) X^r(E) \in [a \sqrt{\mu(E)}, b \sqrt{\mu(E)}]\ \forall\ i,E \bigg\rbrace.$$  For $X \in {\mathcal R}^{q \times |{\mathcal E}|}$, let $f^{\mu,\theta} : {\mathcal R}^{q \times |{\mathcal E}|} \rightarrow {\mathcal R}$ denote the mapping
$$f^{\mu,\theta}(X) \stackrel{\Delta}{=} \sum_{E \in {\mathcal E}} \mu(E) \sum_{r=1}^{q-1} Z^r(E) \frac{X^r(E)}{\sqrt{\mu(E)}} - 2 \iota^{-1} \sum_{S \in {\mathcal S}} \mu(S) \sum_{i=1}^m \phi_{\theta}\big( \sum_{t=1}^T \sum_{r=1}^{q-1} a^r_i(S_{[t]}) \frac{X^r(S_{[t]})}{\sqrt{\mu(S_{[t]})}} - b_i \big).$$  The transformed problem \ref{pen-smooth-mu} (whose optimal value we denote $\textrm{OPT}_{\textsc{pen}^{\mu,\theta}}$ and optimal solution we denote $X_{*,\textsc{pen}^{\mu,\theta}}$) is defined as follows.  
\begin{align}
    \tag{$\textsc{pen}^{\mu, \theta}$} \label{pen-smooth-mu}
    \max f^{\mu, \theta}(X) \quad s.t. \ X \in {\mathcal Q}^3
\end{align}
As \ref{pen-smooth-mu} is a simple rescaling of \ref{pen-smooth}, one may easily verify the following claim.  
\begin{claim}\label{sameopt}
\ref{pen-smooth} and \ref{pen-smooth-mu} have the same optimal value, i.e. $\textrm{OPT}_{\textsc{pen}^{\mu,\theta}} = \textrm{OPT}_{\textsc{pen}^{\theta}}$.  The vector $X$ s.t. $X^r(E) = \sqrt{\mu(E)} X^r_{*,\textsc{pen}^{\theta}}(E)$ is optimal for \ref{pen-smooth-mu}, and we may w.l.o.g. assume $X^r_{*,\textsc{pen}^{\mu,\theta}}(E) = \sqrt{\mu(E)} X^r_{*,\textsc{pen}^{\theta}}(E)$ for all $E \in {\mathcal E}$ and $r \in \lbrace 0,\ldots,q-1 \rbrace$ (which we do in the remainder of the paper).
\end{claim}
For $X \in {\mathcal R}^{q \times |{\mathcal E}|}$, let $\hat{G}^{\mu,k}(X)$ denote the $q \times |{\mathcal E}|$-dimensional vector s.t.
\begin{equation}\label{ghatmukdef}
\hat{G}^{\mu,k}(X)_{E,l} = 
\sqrt{\mu(E)}\Bigg( Z^l(E) - 2 \iota^{-1} \sum_{i =1}^m a^l_{i}(E) \times \eta^{-1}_1 \sum_{S' \in {\mathcal S}^{E,k}}\phi'_{\theta}\bigg( \frac{T}{\eta_2} \sum_{t \in \aleph^k} \sum_{r=1}^{q-1} a^r_{i}({S'}_{[t]}) \frac{X^r({S'}_{[t]})}{\sqrt{\mu(S'_{[t]})}} - b_i \bigg) \Bigg).
\end{equation}
We let $\hat{G}^{\mu,k}(X)_{E}$ denote the $q$-dimensional vector with component $l$ equalling $\hat{G}^{\mu,k}(X)_{E,l}$.
\\\indent We now define our family of linearly transformed SGA methods Algorithm \ref{alg:meta-gd-mu}.  Here ${\mathcal X}^{\mu,k}(E)$ is a q-dimensional vector for each $k \geq -1$ and $E \in {\mathcal E}$, and ${\mathcal X}^{\mu,k}$ is the corresponding $q \times |{\mathcal E}|$-dimensional vector.
\begin{algorithm}[ht]
\caption{\textup{\textsf{Stochastic Gradient for \ref{pen-smooth-mu}}}}
\label{alg:meta-gd-mu}
\vspace{2mm}
   Initialize ${\mathcal X}^{\mu,-1}(E), {\mathcal X}^{\mu,0}(E) \gets \sqrt{\mu(E)} \overline{e}^q\ \forall E \in \mathcal{E}$
\\For $k\gets\ 0\ \KwTo\ K\ \textrm{and}\ E \in \mathcal{E}$
\\\ \ \ {\small ${\mathcal X}^{\mu, k+1}(E) = \Pi_{\sqrt{\mu(E)} {\mathcal Q}^0}\Bigg( (1 + \beta_k) {\mathcal X}^{\mu, k}(E) - \beta_k {\mathcal X}^{\mu, k-1}(E) + \alpha \hat{G}^{\mu, k}\big( (1 + \beta_k) {\mathcal X}^{\mu, k} - \beta_k {\mathcal X}^{\mu, k-1} \big)_E \Bigg).$}
\end{algorithm}
Then we have the following equivalence between the iterates $\lbrace {\mathcal X}^k, k \geq -1 \rbrace$ and (scalings of) $\lbrace {\mathcal X}^{\mu,k}, k \geq -1 \rbrace$.  The result follows from a straightforward induction and the fact that for any polytope $Q$, vector $Y$ in the same dimension, and $\tau \in {\mathcal R}^+$, $\Pi_Q(Y) = \frac{1}{\tau} \Pi_{\tau Q}(\tau Y)$ (a basic and well-known property of Euclidean projection) applied with $Q = Q^0, \tau = \sqrt{\mu(E)}$ and $Y = (1 + \beta_k) {\mathcal X}^{k}(E) - \beta_k {\mathcal X}^{k-1}(E) + \alpha \hat{G}^{k}\big( (1 + \beta_k) {\mathcal X}^{k} - \beta_k {\mathcal X}^{k-1} \big)_E$, and we omit the details.

\begin{claim}\label{sameiterates}
For all $k \geq -1$ and $E \in {\mathcal E}$, ${\mathcal X}^k(E) = \frac{{\mathcal X}^{\mu,k}(E)}{\sqrt{\mu(E)}}$ (i.e. the vectors are equivalent).
\end{claim}
Combining with Claim\ \ref{sameopt} and definitions, we conclude the following.
\begin{corollary}\label{samegap}
For all $k \geq 1$, $\textrm{OPT}_{\textsc{pen}^{\theta}} - f^{\theta}(k^{-1} \sum_{j=1}^k {\mathcal X}^{j}) = 
\textrm{OPT}_{\textsc{pen}^{\mu,\theta}} - f^{\mu,\theta}(k^{-1} \sum_{j=1}^k {\mathcal X}^{\mu,j})$, and
$\textrm{OPT}_{\textsc{pen}^{\theta}} - f^{\theta}({\mathcal X}^{k}) = 
\textrm{OPT}_{\textsc{pen}^{\mu,\theta}} - f^{\mu,\theta}({\mathcal X}^{\mu,k}).$
\end{corollary}
Next, let us explicitly describe the gradients of $f^{\theta}$ and $f^{\mu,\theta}$ for use in later arguments.  First, we state some well-known properties of $\phi_{\theta}$ (\cite{tatarenkoetal2021}).
\begin{claim}\label{huberprops}
$\phi_{\theta}$ is convex and continuously differentiable on ${\mathcal R}$ (with derivative $\phi'_{\theta}$).  Furthermore : 
\begin{itemize}
\item $\phi'_{\theta}(x) = \min(\frac{x^+}{\theta},1)$ for all $x \in {\mathcal R}$, with $x^+ = \max(0,x)$;
\item $\big| \phi'_{\theta}(x) - \phi'_{\theta}(y)\big| \leq \theta^{-1} |x - y|$ for all $x,y \in {\mathcal R}$ (i.e. $\phi'_{\theta}$ is $\theta^{-1}$-Lipschitz);
\item $\phi_{\theta}(x) \leq x^+ \leq \phi_{\theta}(x) + \frac{1}{2} \theta$ for all $x \in {\mathcal R}$.  
\end{itemize}
\end{claim}
Then the following gradient descriptions follow directly from Claim\ \ref{huberprops}, and we omit the details.  
\begin{claim}\label{whatisderiv1}
For all $X \in {\mathcal R}^{q \times |{\mathcal E}|}$, $l \in \lbrace 0,\ldots,q-1 \rbrace$, and $E \in {\mathcal E}$, $$
\nabla f^{\theta}(X)_{E,l} = \mu(E) \bigg( Z^l(E) - 2 \iota^{-1} \sum_{S' \in {\mathcal S} : E \subseteq S'} \frac{\mu(S')}{\mu(E)} \sum_{i=1}^m a^l_i(E) \phi'_{\theta}\big(\sum_{t=1}^T \sum_{r=1}^{q-1} a^r_i(S'_{[t]}) X^r(S'_{[t]}) - b_i\big) \bigg),$$ and
$$
\nabla f^{\mu,\theta}(X)_{E,l} = \sqrt{\mu(E)} \Bigg( Z^l(E) - 2 \iota^{-1} \sum_{S' \in {\mathcal S} : E \subseteq S'} \frac{\mu(S')}{\mu(E)} \sum_{i=1}^m a^l_i(E) \phi'_{\theta}\big(\sum_{t=1}^T \sum_{r=1}^{q-1} a^r_i(S'_{[t]}) \frac{X^r(S'_{[t]})}{\sqrt{\mu(S'_{[t]})}} - b_i\big) \Bigg).$$  
\end{claim}
Next, let us state the relevant convergence results from the literature.  We first state the result in the unaccelerated setting, which follows directly from Theorem 1.1 of \cite{nemirovski2012} using the Euclidean norm, applied to our setting (after converting concave maximization \ref{pen-smooth-mu} to a corresponding convex minimization, and using the fact that (1) $\sum_{E \in {\mathcal E}} \mu(E) = T$, (2) $\hat{G}^{\mu,k}$ has the same distribution for all $k$, (3) $\sum_{r=0}^{q-1} x^2_r \leq \alpha^2$ for any $\alpha > 0$ and $x \in \alpha {\mathcal Q}^0$, (4) $\max_{X \in {\mathcal Q}^3} \sum_{E \in {\mathcal E}} \sum_{l=0}^{q-1} X^2_{E,l} \leq \sum_{E \in {\mathcal E}} (\sqrt{\mu(E)})^2 = T$, and
\begin{equation}\label{inQ4a}
X \in {\mathcal Q}^3\ \textrm{and}\ a^r_i(E) \in [0,1]\ \forall\ r,i,E\ \rightarrow\ X \in {\mathcal Q}^4(0,1).
\end{equation}
Indeed, we will in later proofs (to bound certain norms) repeatedly use this property that $\sum_{E \in {\mathcal E}} \mu(E) = T$, which follows from the fact that $\sum_{E \in {\mathcal E}} \mu(E) = \sum_{t=1}^T \sum_{E \in {\mathcal S}^t} \mu(E) = \sum_{t=1}^T 1 = T$.  Recall that as we are using the Euclidean norm, given a vector $X \in {\mathcal R}^{q \times |{\mathcal E}|}, ||X|| = \sqrt{ \sum_{E \in {\mathcal E}} \sum_{r=0}^{q-1} \big(X^r(E)\big)^2}$.   
\begin{claim}[\cite{nemirovski2012} Theorem 1.1]\label{convclaim1}
Suppose $\beta_k = 0$ for all $k \geq 1$ (i.e. \textrm{Algorithm} \ref{alg:meta-gd} is projected SGA).  Then for all $k \geq 1$, $\textrm{OPT}_{\textsc{pen}^{\mu,\theta}} - \mathbb{E}\big[f^{\mu,\theta}(k^{-1} \sum_{j=1}^k {\mathcal X}^{\mu,j})\big]$ is at most 
{\small
$$\frac{2 T}{k \alpha} + \alpha \sup_{X \in {\mathcal Q}^4(0,1)} ||\nabla f^{\mu,\theta}(X)||^2 +  1.5 \alpha \sup_{X \in {\mathcal Q}^4(0,1)} \mathbb{E}\bigg[ \big|\big|\nabla f^{\mu,\theta}(X) - \hat{G}^{\mu,1}(X)\big|\big|^2\bigg] + 2 \sqrt{2 T} \sup_{X \in {\mathcal Q}^4(0,1)} \bigg|\bigg| \mathbb{E}\big[\nabla f^{\mu,\theta}(X) - \hat{G}^{\mu,1}(X) \big] \bigg|\bigg|.$$}
\end{claim}
We now state the result in the accelerated setting, which follows directly from Proposition 2 of \cite{schmidtetal2011}.  Let 
$$L^{\mu,\theta} \stackrel{\Delta}{=} \inf\lbrace C \geq 0 : ||\nabla f^{\mu,\theta}(X) - \nabla f^{\mu,\theta}(Y)|| \leq C ||X - Y||\ \forall\ X,Y \in {\mathcal R}^{q \times |{\mathcal E}|}\rbrace$$
denote the smoothness of $f^{\mu,\theta}$, i.e. the Lipschitz continuity parameter of the gradient in the Euclidean norm.  Then combining \cite{schmidtetal2011} Proposition 2 with a straightforward conditioning argument, we conclude the following.  For completeness, we provide a more detailed proof of how our result follows from that of \cite{schmidtetal2011} Proposition 2 in Appendix\ \ref{convclaim2aproofsec}.
\begin{claim}[\cite{schmidtetal2011} Proposition 2]\label{convclaim2a}
Suppose $\beta_0 = 0, \beta_k = \frac{k-1}{k+2}$ for all $k \geq 1,$ and $\alpha \in (0,\frac{1}{L^{\mu,\theta}}]$.  Then for all $k \geq 1$, $\textrm{OPT}_{\textsc{pen}^{\mu,\theta}} - \mathbb{E}\big[f^{\mu,\theta}({\mathcal X}^{\mu,k})\big]$ is at most $\frac{16 T}{\alpha k^2} + 16 k^2 \alpha \sup_{X \in {\mathcal Q}^4(-1,2)} 
\mathbb{E}\big[ \big|\big| \nabla f^{\mu,\theta}(X) - 
\hat{G}^{\mu,1}(X)\big|\big|^2 \big].$
\end{claim}
By bounding the suprema appearing in Claims\ \ref{convclaim1} and \ref{convclaim2a} using Hoeffding's inequality, and separately bounding $L^{\mu,\theta}$, we prove the following, deferring the proof to Appendix\ \ref{convclaim2proofsec}.
\begin{corollary}\label{convclaim2}
Suppose $\beta_k = 0$ for all $k \geq 1$.  Then for all $k \geq 1, \alpha > 0$, and positive integers $\eta_1 \geq 1$  and $\eta_2 \in \lbrace 1,\ldots,T \rbrace$, $\textrm{OPT}_{\textsc{pen}^{\mu,\theta}} - \mathbb{E}\big[f^{\mu,\theta}(k^{-1} \sum_{j=1}^k \mathcal{X}^{\mu,j})\big]$ is at most 
$$
\frac{2  T}{k \alpha} + 4 \alpha \frac{L^2 T}{\iota^2} q + 1.5 \alpha \big( \frac{72 L^2 T^3}{\iota^2 \theta^2 \eta_2} + \frac{8 L^2 T}{\iota^2 \eta_1} \big) q + 2 \sqrt{2 T} \sqrt{\frac{72 L^2 T^3}{\iota^2 \theta^2 \eta_2} + \frac{8 L^2 T}{\iota^2 \eta_1}} \sqrt{q}.
$$
If in addition $\eta_2 = T$ (i.e. the relevant sum is computed exactly), then 
$$\textrm{OPT}_{\textsc{pen}^{\mu,\theta}} - \mathbb{E}\big[f^{\mu,\theta}(k^{-1} \sum_{j=1}^k \mathcal{X}^{\mu,j})\big] \leq 
\frac{2  T}{k \alpha} + 4 \alpha \frac{L^2 T}{\iota^2} q + 1.5 \alpha \times \frac{8 L^2 T}{\iota^2 \eta_1} q + 2 \sqrt{2 T} \sqrt{\frac{8 L^2 T}{\iota^2 \eta_1}} \sqrt{q}.
$$
Suppose $\beta_0 = 0$, $\beta_k = \frac{k-1}{k+2}$ for all $k \geq 1$ and $\alpha \in \big(0 , \frac{1}{2} \iota \theta (U_r L W)^{-\frac{1}{2}}\big]$.  Then for all $k \geq 1$ and positive integers $\eta_1 \geq 1$  and $\eta_2 \in \lbrace 1,\ldots,T \rbrace$,
$$\textrm{OPT}_{\textsc{pen}^{\mu,\theta}} - \mathbb{E}\big[f^{\mu,\theta}(\mathcal{X}^{\mu,k})\big]
\leq \frac{16 T}{\alpha k^2} + 16 k^2 \alpha \big( \frac{72 L^2 T^3}{\iota^2 \theta^2 \eta_2} + \frac{8 L^2 T}{\iota^2 \eta_1} \big) q.$$
If in addition $\eta_2 = T$ (i.e. the relevant sum is computed exactly), then
$$\textrm{OPT}_{\textsc{pen}^{\mu,\theta}} - \mathbb{E}\big[f^{\mu,\theta}(\mathcal{X}^{\mu,k})\big]
\leq \frac{16 T}{\alpha k^2} + 16 k^2 \alpha \times \frac{8 L^2 T}{\iota^2 \eta_1} q.$$
\end{corollary}

\subsubsection{Proof of Theorem\ \ref{convclaim1cor}}
We now complete the proof of Theorem\ \ref{convclaim1cor} by combining Corollary\ \ref{convclaim2} with some straightforward algebra.
\begin{proof}{Proof of Theorem\ \ref{convclaim1cor} :}
Let us first treat the case $\beta_k = 0$ for all $k$.  Combining Corollaries\ \ref{convclaim2} and\ \ref{samegap} with the fact that $\sqrt{x + y} \leq \sqrt{x} + \sqrt{y}$, it suffices to have $\frac{2T}{k \alpha} \leq \frac{1}{6} \epsilon T$ (equivalently $k \geq \frac{12}{\alpha \epsilon})$, $4 \alpha L^2 T \iota^{-2} q \leq \frac{1}{6} \epsilon T$ (equivalently $\alpha \leq \frac{1}{24} L^{-2} \iota^2 q^{-1} \epsilon)$, $\frac{1.5 \times 72 L^2 T^3}{\iota^2 \theta^2 \eta_2} \alpha q \leq \frac{1}{6} \epsilon T$ or $\eta_2 = T$ (equivalently $\eta_2 \geq \min\big(\frac{648 L^2 T^2}{\iota^2 \theta^2 \epsilon} \alpha q,T\big)$), $\frac{1.5 \times 8 L^2 T}{\iota^2 \eta_1} \alpha q \leq \frac{1}{6} \epsilon T$ (equivalently $\eta_1 \geq \frac{72 L^2}{\iota^2 \epsilon} \alpha q$), $2 \sqrt{2 T} \sqrt{\frac{72 L^2 T^3}{\iota^2 \theta^2 \eta_2}} \sqrt{q} \leq \frac{1}{6} \epsilon T$ or $\eta_2 = T$ (equivalently $\eta_2 \geq \min\big(\frac{20736 T^2 L^2}{\iota^2 \theta^2 \epsilon^2} q,T\big)$), and $2 \sqrt{2 T} \sqrt{\frac{8 L^2 T}{\iota^2 \eta_1}} \sqrt{q} \leq \frac{1}{6} \epsilon T$ (equivalently $\eta_1 \geq \frac{2304 L^2}{\iota^2 \epsilon^2} q$).  Combining the above completes the proof in this case.
\\\indent Next, let us treat the case $\beta_0 = 0, \beta_k = \frac{k-1}{k+2}$ for all $k \geq 1$.  Combining Corollaries\ \ref{convclaim2} and\ \ref{samegap}, it suffices to have $\alpha \leq \frac{\iota \theta}{2\sqrt{U_r L W}}$, $\frac{16 T}{\alpha k^2} \leq \frac{1}{3} \epsilon T$ (equivalently $\alpha k^2 \geq \frac{48}{\epsilon}$), $\frac{1152 \alpha k^2 L^2 T^3 q}{\iota^2 \theta^2 \eta_2} \leq \frac{1}{3} \epsilon T$ or $\eta_2 = T$ (equivalently $\eta_2 \geq \min\big( 3456 \times \alpha k^2 \times \frac{L^2 T^2}{\iota^2 \theta^2 \epsilon} \times q, T\big)$), and $\frac{128 \alpha k^2 L^2 T q}{\iota^2 \eta_1} \leq \frac{1}{3} \epsilon T$ (equivalently $\eta_1 \geq 384 \times \alpha k^2 \times \frac{L^2}{\iota^2 \epsilon} \times q$).  We will select $\alpha,k$ carefully to avoid having to use too many ceiling, min, or max operations.  Note that $\frac{1}{4} \lceil \frac{ (U_r L W)^{\frac{1}{4}}}{\sqrt{\iota \theta}} \rceil^{-2} \leq \frac{\iota \theta}{2\sqrt{U_r L W}}$, and thus we may take $\alpha = \frac{1}{4} \lceil \frac{ (U_r L W)^{\frac{1}{4}}}{\sqrt{\iota \theta}} \rceil^{-2}$.  Note that for this choice of $\alpha$, $\alpha \times \big( 14 \lceil \frac{ (U_r L W)^{\frac{1}{4}}}{\sqrt{\iota \theta}} \rceil \lceil \epsilon^{-\frac{1}{2}} \rceil  \big)^2 = 49  \lceil \epsilon^{-\frac{1}{2}} \rceil^2 \geq \frac{48}{\epsilon}$, and thus we may take $k = 14 \lceil \frac{ (U_r L W)^{\frac{1}{4}}}{\sqrt{\iota \theta}} \rceil \lceil \epsilon^{-\frac{1}{2}} \rceil.$  For this choice of $\alpha,k$, it holds that $\alpha k^2 = 49  \lceil \epsilon^{-\frac{1}{2}} \rceil^2 \leq 49 \big( 2 \epsilon^{-\frac{1}{2}} \big)^2 \leq 200 \epsilon^{-1}$.  Thus (plugging into our previous bounds and simplifying) it suffices to have 
$\eta_2 \geq \min\big( 691200 \frac{L^2 T^2}{\iota^2 \theta^2 \epsilon^2} q, T\big)$ and $\eta_1 \geq 142800 \frac{L^2}{\iota^2 \epsilon^2} q$.  Combining the above completes the proof.  $\qed$
\end{proof}

\subsection{Proof of Lemma\ \ref{comparefeasible1}}\label{APPcomparefeasible1sec}
It will simplify the exposition and intuition of our proof to begin with a more general result.  Consider a non-negative sequence $x_1,\ldots,x_N$ (for some $N \geq 1$).  Suppose we are given budgets $b_1,\ldots,b_m > 0$, and a general $m$ by $N$ non-negative matrix $A$ with all strictly positive components at least $\iota$.  In correspondence with our earlier definition of $\textrm{F}$, let us define $F(x_1) = \min\big( x_1, \min_{i : A_{i,1} > 0} \frac{b_i}{A_{i,1}} \big)$ if\ $\exists\ i : A_{i,1} > 0$, and $F(x_1) = x_1$ otherwise, and inductively
$F(x_n) = \min\big( x_n, \min_{i : A_{i,n} > 0} \frac{b_i - \sum_{u=1}^{n-1} A_{i,u} F(x_u)}{A_{i,n}} \big)$ if\ $\exists\ i : A_{i,n} > 0$, and $F(x_n) = x_n$ otherwise (for $n \geq 2$).  Note that $0 \leq F(x_n) \leq x_n$ for all $n$, and it is easily verified (by a straightforward induction the details of which we omit) that $b_i - \sum_{u=1}^{n-1} A_{i,u} F(x_u) \geq 0$ for all $n$ and $i$.  It thus follows (by a straightforward monotonicity argument) that so long as $|\lbrace i : A_{i,n} > 0 \rbrace| > 0$, one has
\begin{equation}\label{Fxnlb}
F(x_n) \geq \min\big( x_n, \min_{i : A_{i,n} > 0} \frac{\big(b_i - \sum_{u=1}^{n-1} A_{i,u} x_u\big)^+}{A_{i,n}} \big).
\end{equation}
Ultimately $x_1,\ldots,x_N$ will correspond to $X^r(S_{[t]}), r = 1,\ldots,q-1, t = 1,\ldots,T$ (with $N = (q-1) T$) and row $i$ of $A$ will correspond to $a^r_i(S_{[t]}), r = 1,\ldots,q-1, t = 1,\ldots,T$.  Then the key result (from which Lemma\ \ref{comparefeasible1} follows in a straightforward manner after summing over $n$) is the following.
\begin{lemma}\label{inductionstep1}
\begin{equation}\label{inductionstep1eq}
x_n - F(x_n) \leq \iota^{-1} \sum_{i=1}^m  \bigg( (\sum_{u=1}^n A_{i,u} x_u - b_i)^+ - (\sum_{u=1}^{n-1} A_{i,u} x_u - b_i)^+ \bigg)\ \forall\ n.
\end{equation}
\end{lemma}
\begin{proof}{Proof of Lemma\ \ref{inductionstep1}}
If $|\lbrace i : A_{i,n} > 0 \rbrace| = 0$ the bound holds trivially.  If $|\lbrace i : A_{i,n} > 0 \rbrace| \neq 0$, then
\begin{eqnarray*}
x_n - F(x_n) &\leq& x_n - \min\big( x_n, \min_{i : A_{i,n} > 0} \frac{\big(b_i - \sum_{u=1}^{n-1} A_{i,u} x_u\big)^+}{A_{i,n}} \big)\ \textrm{by}\ (\ref{Fxnlb})
\\&=& \max\bigg(0, x_n - \min_{i : A_{i,n} > 0} \frac{\big(b_i - \sum_{u=1}^{n-1} A_{i,u} x_u\big)^+}{A_{i,n}}\bigg)
\\&=& \max\bigg(0, \max_{i : A_{i,n} > 0} \big(x_n -  \frac{\big(b_i - \sum_{u=1}^{n-1} A_{i,u} x_u\big)^+}{A_{i,n}} \big) \bigg)
\\&=& \max\bigg(0, \max_{i : A_{i,n} > 0} \frac{1}{A_{i,n}} \big(A_{i,n} x_n - \big(b_i - \sum_{u=1}^{n-1} A_{i,u} x_u\big)^+ \big) \bigg)
\\&\leq& \iota^{-1} \max\bigg(0, \max_{i : A_{i,n} > 0} \big( A_{i,n} x_n - (b_i - \sum_{u=1}^{n-1} A_{i,u} x_u)^+ \big) \bigg)
\\&\leq& \iota^{-1} \sum_{i=1}^m \big( A_{i,n} x_n - (b_i - \sum_{u=1}^{n-1} A_{i,u} x_u)^+ \big)^+
\end{eqnarray*}
Alternatively, note that $(\sum_{u=1}^n A_{i,u} x_u - b_i)^+ - (\sum_{u=1}^{n-1} A_{i,u} x_u - b_i)^+$ equals $A_{i,n} x_n$ if $\sum_{u=1}^{n-1} A_{i,u} x_u - b_i \geq 0$ (equivalently $(b_i - \sum_{u=1}^{n-1} A_{i,u} x_u)^+ = 0$), and otherwise equals $(\sum_{u=1}^n A_{i,u} x_u - b_i)^+$.  Equivalently
$$(\sum_{u=1}^n A_{i,u} x_u - b_i)^+ - (\sum_{u=1}^{n-1} A_{i,u} x_u - b_i)^+ = \big( A_{i,n} x_n - (b_i - \sum_{u=1}^{n-1} A_{i,u} x_u)^+ \big)^+.$$
Combining the above completes the proof.  $\qed$
\end{proof}
\ \indent Recognizing that the right-hand-side of (\ref{inductionstep1eq}) telescopes (if summed over $n$), we find that Lemma\ \ref{inductionstep1} directly implies the following corollary.
\begin{corollary}\label{inductionstep1cor}
$\sum_{n=1}^N  \big( x_n - F(x_n) \big) \leq \iota^{-1} \sum_{i=1}^m  (\sum_{u=1}^N A_{i,u} x_u - b_i)^+.$
\end{corollary}
Corollary\ \ref{inductionstep1cor} directly implies the following result for our original problem, under the correspondence (for any fixed $S \in {\mathcal S}$) $N = (q-1) T$, $x_n = X^r(S_{[t]})$ and $A_{i,n} = a^r_i(S_{[t]})$ for $n = (t-1) \times (q-1) + r$, also using the fact that $Z^r(E) \in [0,1]$ for all $E$ and the definitional equivalence between the $F$ operator here and $\textrm{F}$ as previously defined.
\begin{corollary}\label{inductionstep1cor2}
{\small
$$\sum_{S \in {\mathcal S}} \mu(S) \big(  \sum_{t=1}^T \sum_{r=1}^{q-1} Z^r(S_{[t]}) X^r(S_{[t]}) - \sum_{t=1}^T \sum_{r=1}^{q-1} Z^r(S_{[t]}) \textrm{F}(X)^r(S_{[t]}) \big) \leq 
\iota^{-1} \sum_{S \in {\mathcal S}} \mu(S) \sum_{i=1}^m \big( \sum_{t=1}^T \sum_{r=1}^{q-1} a^r_i(S_{[t]}) X^r(S_{[t]}) - b_i \big)^+.$$}
\end{corollary}
\ \indent Recognizing that the left-hand-side corresponds to the difference of the expected reward earned under $X$ and $\textrm{F}(X)$ then completes the proof of Lemma\ \ref{comparefeasible1}.  $\qed$

\section{Electronic Compendium : Supplemental Appendix}\label{suppsec1}

\subsection{Analysis of Example 1 from Section\ \ref{introsec1}}\label{heuristicformalsec}
In this section we more formally define and analyze the heuristics from Section\ \ref{introsec1}, and compute the optimal policy (implying the claimed optimal solution for LP (\ref{LP000})).  Let us denote the multiple stopping problem of Example 1 by $\textrm{PROB}^1$.  We begin by proving that any feasible policy for $\textrm{PROB}^1$ that stops too many times in the first $\lfloor \frac{T}{3} \rfloor - 2$ periods (i.e. accepts too many rewards of value $.5$) in expectation must have a significant suboptimality gap relative to the optimal policy, which does not stop at all in the first $\lfloor \frac{T}{3} \rfloor - 2$ periods.  To accommodate the relevant heuristic policies we allow the policies to use exogenous randomization, and refer the reader to Section\ \ref{formulationsec} for further discussion of randomized policies in this context.  For a general feasible (possibly randomized) policy $\pi$, let $X^1_{\pi,t}$ denote the (random) value of $X^1_t$ under policy $\pi$, and $N_{\pi,.5} \stackrel{\Delta}{=} \sum_{t=1}^{\lfloor \frac{T}{3} \rfloor - 2} X^1_{\pi,t}$ denote the number of times $\pi$ stops at a reward of value $.5.$  Recall that $b = \lfloor \frac{T}{3} \rfloor - 2$.
\begin{lemma}\label{stoptoomuch}
An optimal policy for $\textrm{PROB}^1$ is for the DM to skip all the initial rewards of value .5, and then : (1) if the reward at time $\lfloor \frac{T}{3} \rfloor - 1$ is .001 the DM uses their entire budget on rewards of value $.45$; and (2) if instead the reward at time $\lfloor \frac{T}{3} \rfloor - 1$ is .01 the DM uses their entire budget on rewards of value $1$.  Furthermore, $\textrm{OPT}_T = .3 (\lfloor \frac{T}{3} \rfloor - 2) + .7 \times .45 \times (\lfloor \frac{T}{3} \rfloor - 2)$, and for any feasible policy $\pi$, $\textrm{OPT}_T - \textrm{VAL}_{\pi,T} \geq \big( .3 + .7 \times .45 - .5 \big) \mathbb{E}[N_{\pi,.5}]$.
\end{lemma}
\begin{proof}{Proof of Lemma\ \ref{stoptoomuch}}
For $t \in \lbrace 1,\ldots,T \rbrace, E \in {\mathcal S}^t$, and $B \in \lbrace 0,\ldots,b \rbrace$, let $V_t(E,B)$ denote the reward-to-go at time $t$ given partial reward sequence $E$, i.e. the greatest expected cumulative reward (in periods $t,\ldots,T$) achievable by any feasible policy given the partial reward sequence in periods $1,\ldots,t$ is $E$, where we refer the reader to \cite{puterman2014} for further details on reward-to-go functions.  It is easily verified that $V_{\lfloor\frac{T}{3}\rfloor - 1}\big( (.5,\ldots,.5,.01) , B \big) = B$ and $V_{\lfloor\frac{T}{3}\rfloor - 1}\big( (.5,\ldots,.5,.001) , B \big) = .45 B.$  It then follows from the fact that (1) the reward sequence in periods $1,\ldots,\lfloor\frac{T}{3}\rfloor - 2$ is deterministic, and (2) $V_t(E,B)$ is an upper bound on the cumulative expected reward (in periods $t,\ldots,T$) of any policy $\pi$ with remaining budget $B$ in period $t$ having encountered partial reward sequence $E$ in periods $1,\ldots,t$, that : (a) $\textrm{OPT}_T = \max_{k = 0,\ldots,b} \big( .5 k + .3 (b - k) + .7 \times .45 \times (b - k) \big)$, (b) for any feasible policy $\pi$ it holds that $\textrm{VAL}_{\pi,T} \leq .5 \mathbb{E}[N_{\pi,.5}] + .3 (b - \mathbb{E}[N_{\pi,.5}]) + .7 \times .45 \times (b - \mathbb{E}[N_{\pi,.5}])$, and (c) an optimal policy for $\textrm{PROB}^1$ is for the DM to accept exactly $argmax_{k = 0,\ldots,b} \big( .5 k + .3 (b - k) + .7 \times .45 \times (b - k) \big)$ rewards of value $.5$ and then if the reward at time $\lfloor \frac{T}{3} \rfloor - 1$ is .001 the DM uses their entire remaining budget on rewards of value $.45$, otherwise the DM uses their entire remaining budget on rewards of value $1$.  As $.5 k + .3 (b - k) + .7 \times .45 \times (b - k) = \big(.5 - .3  - .7 \times .45 \big) k + \big( .3 b + .7 \times .45 \times b \big)$ is a linear function of $k$, and $.5 - .3  - .7 \times .45 = -.115 < 0$, we immediately deduce the desired properties of the optimal policy and its expected reward.  To complete the proof, we observe that (by the above) for any feasible policy $\pi$, it holds that 
$\textrm{VAL}_{\pi,T} \leq .5 \mathbb{E}[N_{\pi,.5}] + .3 (b - \mathbb{E}[N_{\pi,.5}]) + .7 \times .45 \times (b - \mathbb{E}[N_{\pi,.5}])$, while $\textrm{OPT}_T = .3 b + .7 \times .45 \times b$, and the desired result follows.  $\qed$
\end{proof}
\ \indent We now prove that all three heuristic policies are suboptimal even on the fluid scale.  
\begin{lemma}\label{allaresubopt}
For $\pi \in \lbrace CE, Bayes, FBayes \rbrace$, it holds that $\textrm{OPT}_T - \textrm{VAL}_{\pi,T} \geq .023 \big( \lfloor \frac{T}{3} \rfloor - 2 \big).$  Thus for large $T$ (on the fluid scale)
$\textrm{OPT}_T - \textrm{VAL}_{\pi,T} \geq  .007 \times T.$
\end{lemma}
To prove Lemma\ \ref{allaresubopt}, we will apply Lemma\ \ref{stoptoomuch} to the 3 heuristics, showing they all stop too many times in the first $\lfloor \frac{T}{3} \rfloor - 2$ periods.  We begin by more formally defining these heuristics.
\\\textbf{C-E heuristic.}  The C-E heuristic is typically applied by reasoning about the set of customer types, with each type defined by its reward and r.c.v.  As in our setting all r.c.v.s are the scalar 1, the type is thus defined by the reward.  We thus define 6 types of customers : type 1 (with reward .5), type 2 (with reward $.45$), type 3 (with reward .01), type 4 (with reward .001), type 5 (with reward 1), and type 6 (with reward 0).  Let $A^i_{[t,T]}$ denote the expected number of type-i arrivals in $[t,T]$, and $r^i$ denote the reward of a type-i customer.  Then at each time $t$, the DM will solve the following LP : 
\begin{equation}\label{LPtosolve1}
\max \sum_{i=1}^6 r^i x_i\ \ s.t.\ \ \sum_{i=1}^6 x_i \leq b - \sum_{s=1}^{t-1} X^1_{CE,s} \ ,\ x_i \leq A^i_{[t,T]}\ ,\ x_i \geq 0\ \forall\ i.
\end{equation}
Supposing the arrival at time $t$ is of type $i$, the DM then sets $X^1_{CE,t} = 1$ w.p. $\frac{x_i}{A^i_{[t,T]}}$, and sets $X^1_{CE,t} = 0$ otherwise.  If setting $X^1_{CE,t} = 1$ turns out to be infeasible, the DM simply sets $X^1_{CE,t} = 0$ instead.
\\\textbf{Fluid Bayes Selector heuristic.}  The Fluid Bayes Selector heuristic is very similar to the C-E heuristic, except it uses hard-thresholding instead of randomization.  In particular, at each time $t$ the DM solves the same LP (\ref{LPtosolve1}) (with $X^1_{CE,s}$ replaced by $X^1_{FBayes,s}$), but now sets $X^1_{FBayes,t} = 1$ iff $\frac{x_i}{A^i_{[t,T]}} \geq \frac{1}{2}$, again setting $X^1_{FBayes,t} = 0$ if doing otherwise would lead to infeasibility.
\\\textbf{Bayes Selector heuristic.}  The Bayes Selector heuristic works a bit differently.  For $t = 1,\ldots,T$ and $\hat{S} \in \lbrace S^0,S^1 \rbrace$, let $IP_t(\hat{S})$ denote the following IP :     
\begin{equation}\label{Bayestosolve1}
\max \sum_{s=t}^T Z^1_s(\hat{S}_{[s]}) X_s\ \ s.t.\ \ \ \sum_{s=t}^T X_s \leq b - \sum_{s=1}^{t-1} X^1_{Bayes,s}\ ,\ X_s \in \lbrace 0,1 \rbrace\ \textrm{for}\ s = t,\ldots,T.
\end{equation}
For $\hat{S} \in \lbrace S^0,S^1 \rbrace$ and $i \in \lbrace 0,1 \rbrace$, let $A_t(i,\hat{S})$ denote the event that there exists an optimal solution to $IP_t(\hat{S})$ in which $X_t = i$.  Then under the Bayes Selector heuristic, the DM sets $X^1_{Bayes,t} = 1$ if $.3 \times I\big( A_t(1, S^1) \big) + .7 \times I\big( A_t(1, S^0) \big) \geq .3 \times I\big( A_t(0, S^1) \big) + .7 \times I\big( A_t(0, S^0) \big)$, and sets $X^1_{Bayes,t} = 0$ otherwise, and again setting $X^1_{Bayes,t} = 0$ if doing otherwise would lead to infeasibility.  The DM's action can be interpreted as trying to maximize the probability of agreement with the h.r., and we refer the reader to \cite{veraetal2021} for further details.
\begin{proof}{Proof of Lemma\ \ref{allaresubopt}}
We begin by proving a bound on $\mathbb{E}[N_{\pi,.5}]$ for $\pi \in \lbrace CE, FBayes \rbrace$.  For $\pi \in \lbrace CE, FBayes \rbrace$ and $t \in \lbrace 1,\ldots,\lfloor\frac{T}{3}\rfloor - 2 \rbrace$, LP (\ref{LPtosolve1}) becomes
\\$\max .5 x_1 + .45 x_2 + .01 x_3 + .001 x_4 + x_5$
\\s.t.
\\$\sum_{i=1}^6 x_i \leq \lfloor \frac{T}{3} \rfloor - 2 - \sum_{s=1}^{t-1} X^1_{\pi,s}$
\\$x_1 \leq \lfloor \frac{T}{3} \rfloor - 2 - (t-1) \ ,\ x_2 \leq T - 2 \lfloor \frac{T}{3} \rfloor + 3\ ,\ x_3 \leq .3\\\
x_4\ \leq\ .7\ ,\ x_5\ \leq\ .3 (\lfloor \frac{T}{3} \rfloor - 2)\ ,\ x_6 \leq .7(\lfloor \frac{T}{3} \rfloor - 2), x_i \geq 0\ \forall\ i$
\\It follows from the basic properties of (budgeted) fractional knapsack problems that so long as $\lfloor \frac{T}{3} \rfloor - 2 - \sum_{s=1}^{t-1} X^1_{\pi,s} \leq .3 (\lfloor \frac{T}{3} \rfloor - 2)$, any optimal solution to LP (\ref{LPtosolve1}) will set $x_5 = .3 (\lfloor \frac{T}{3} \rfloor - 2)$, and 
\begin{eqnarray*}
x_1 &=& \min\bigg( \lfloor \frac{T}{3} \rfloor - 2 - (t-1) , \lfloor \frac{T}{3} \rfloor - 2 - \sum_{s=1}^{t-1} X^1_{\pi,s} - .3 (\lfloor \frac{T}{3} \rfloor - 2) \bigg) 
\\&\geq& \min\bigg( \lfloor \frac{T}{3} \rfloor - 2 - (t-1) , \lfloor \frac{T}{3} \rfloor - 2 - (t-1) - .3 (\lfloor \frac{T}{3} \rfloor - 2) \bigg)\ \textrm{since}\ \sum_{s=1}^{t-1} X^1_{\pi,s} \leq t - 1
\\&=&  .7 \lfloor \frac{T}{3} \rfloor - 1 - t + 2 \times .3.
\end{eqnarray*}
Thus so long as $t \in \lbrace 1,\ldots,\lfloor\frac{T}{3}\rfloor - 2 \rbrace$ and $.7 \lfloor \frac{T}{3} \rfloor - 1 - t + 2 \times .3 \geq .5 \big(\lfloor \frac{T}{3} \rfloor - t - 1\big)$, both $CE$ and $FBayes$ will set $X^1_t = 1$ with probability at least $\frac{1}{2}$, independent of anything else (including any randomization due to the rounding of the C-E heuristic, and any randomness in the random reward sequence itself).  We note that this is an underestimate as $FBayes$ will set $X^1_t = 1$ w.p.1 in this case, but here for simplicity we only seek a bound.  Since $.7 \lfloor \frac{T}{3} \rfloor - 1 - t + 2 \times .3 \geq .5 \big(\lfloor \frac{T}{3} \rfloor - t - 1\big)$ iff $t \leq .4 \lfloor \frac{T}{3} \rfloor + .2$, and applying some elementary bounds involving the floor operator, we conclude that both $\mathbb{E}[N_{CE,.5}]$ and $\mathbb{E}[N_{FBayes,.5}]$ are at least 
$.2 \big( \lfloor \frac{T}{3} \rfloor - 2 \big)$.
\\\indent We now prove a bound on $\mathbb{E}[N_{Bayes,.5}]$.  Indeed, we will prove (by induction) that w.p.1 $X^1_{Bayes,t} = 1$ for $t = 1,\ldots,\lfloor\frac{T}{3}\rfloor - 2$, so $\mathbb{E}[N_{Bayes,.5}] = \lfloor\frac{T}{3}\rfloor - 2$.  We begin with the base case $t = 1$.  Note that w.p.1 $IP_1(S^1)$ is a knapsack problem with unique optimal solution to set $X_t = 1$ for $t = T - \lfloor \frac{T}{3} \rfloor + 3,\ldots,T$ and all other $x_t$ to zero, while $IP_1(S^0)$ is a knapsack problem with unique optimal solution to set $X_t = 1$ for $t = 1,\ldots,\lfloor\frac{T}{3}\rfloor-2$ and all other $X_t$ to zero.  It follows that $I\big( A_1(1, S^1) \big) = 0, I\big( A_1(0, S^1) \big) = 1, I\big( A_1(1, S^0) \big) = 1, I\big( A_1(0, S^0) \big) = 0, .3 \times I\big( A_1(1, S^1) \big) + .7 \times I\big( A_1(1, S^0) \big) = .7, .3 \times I\big( A_1(0, S^1) \big) + .7 \times I\big( A_1(0, S^0) \big) = .3$.  Thus $X^1_{Bayes,1} = 1$ w.p.1.  Now, let us assume $X^1_{Bayes,s} = 1$ w.p.1 for all $s \leq t$ for some $t \leq \lfloor\frac{T}{3}\rfloor-3$, and prove that $X^1_{Bayes,t+1} = 1$ w.p.1.  The proof for the induction step is essentially identical to that of the base case, and we omit the details.  Combining the above completes the proof.
\\\indent Combining the above, we conclude that for $\pi \in \lbrace CE, FBayes, Bayes \rbrace$, it holds that $\mathbb{E}[N_{\pi,.5}] \geq .2 \big( \lfloor \frac{T}{3} \rfloor - 2 \big)$.  Combining with Lemma\ \ref{stoptoomuch} then completes the proof.  $\qed$
\end{proof}
Further combining Lemmas\ \ref{stoptoomuch} and \ref{allaresubopt} with some straightforward algebra and limit calculations completes the proof of all claims made in Section\ \ref{introsec1} regarding Example 1.

\subsection{Additional commentary on whether a PTAS can exist in our model}\label{noptasstats}
Here we comment further on the potential for PTAS-type algorithms in our model.  First, we argue that due to the data-driven nature of our model, such a PTAS-type result is impossible in general as such an algorithm would facilitate the accurate estimation of exponentially rare events using only polynomially many samples.  In fact, we show that no polynomial-time algorithm can even extract $99\%$ of the optimal value.  Recall that multiple stopping is the variant of $\textsc{NRM}$ in which $L = m = 1, q = 2, a^1_{1,t} = 1$ for all $t$ w.p.1, and we denote its optimal value by $\textrm{OPT}_{\textsc{MS}}$.
\begin{lemma}\label{noptas}
There does not exist a polynomial $\textrm{poly}(\cdot)$ s.t. for some admissible policy $\textrm{A}$ for $\textsc{MS}$, and all sufficiently large $T$, on any trajectory $S \in \mathcal{S}$ one can compute decisions $\textrm{A}(S_{[t]})\ \forall\ t = 1, \dots, T$ on-the-fly, in per-decision time ${\mathcal C} \times \textrm{poly}(T)$, s.t. $\mathbb{E}\big[ \sum_{E \in {\mathcal E}} \mu(E) Z^1(E) \textrm{A}^1(E) \big] \geq .99 \times \textrm{OPT}_{\textsc{MS}}$, and this is true even when the budget $b_1$ is restricted to equalling $\lfloor \frac{T}{2} \rfloor$.
\end{lemma}
\begin{proof}{Proof of Lemma\ \ref{noptas}}
We proceed by arguing that if such an algorithm existed, one could use it to distinguish between two simulators that on all inputs agree with probability exponentially close to 1, using only a polynomial number of samples, leading directly to a contradiction.  We note that one could also frame this argument completely formally in the language of statistical minimax theorems, see e.g. \cite{tsybakov09}, but as the arguments are fairly intuitive and straightforward we instead proceed from first principles in a slightly more informal manner to maximize readability.
\\\indent Thus suppose Lemma\ \ref{noptas} is false, and let $\textrm{poly}(\cdot)$ denote the polynomial bounding the computational and simulation complexity.  Note that there exists $T_0$ s.t. $T \geq T_0$ implies $\frac{1}{2} \times 10^5 \times (T+1) \times \textrm{poly}(T) \times 2^{-T} \leq .01$.  Let us pick some even $T \geq T_0$, and define two multiple stopping problems, $\textsc{MS}^1_T$ and $\textsc{MS}^0_T$, where for both problems the o.i.p. consists only of the realized rewards and the budget is $\frac{T}{2}$.  In problem $\textsc{MS}^1_T$, in periods $1,\ldots,\frac{T}{2}$, the reward is deterministically $.5 \times 2^{-T}$.  In periods $\frac{T}{2} + 1,\ldots,T$ w.p. $2^{-T}$ all these rewards are 1, and w.p. $1 - 2^{-T}$ all these rewards are 0, independent of anything else.  In problem $\textsc{MS}^0_T$, in periods $1,\ldots,\frac{T}{2}$, the reward is deterministically $.5 \times 2^{-T}$.  In periods $\frac{T}{2} + 1,\ldots,T$ w.p. $2^{-100 T}$ all these rewards are 1, and w.p. $1 - 2^{-100 T}$ all these rewards are 0, independent of anything else.  Let $SIM^1_T$ denote the simulator corresponding to problem $\textsc{MS}^1_T$, and $SIM^0_T$ denote the simulator corresponding to problem $\textsc{MS}^0_T$.  Note that for $E \in {\mathcal S}^t$ for $t = 1,\ldots,\frac{T}{2}$, 
$SIM^0_T(E)$ outputs a trajectory consisting of $\frac{T}{2}$ rewards of value $.5 \times 2^{-T}$ and $\frac{T}{2}$ rewards of value $1$ w.p. $2^{-T}$, and outputs a trajectory consisting of $\frac{T}{2}$ rewards of value $.5 \times 2^{-T}$ and $\frac{T}{2}$ rewards of value $0$ w.p. $1 - 2^{-T}$; while $SIM^0_T(E)$ outputs a trajectory consisting of $\frac{T}{2}$ rewards of value $.5 \times 2^{-T}$ and $\frac{T}{2}$ rewards of value $1$ w.p. $2^{-100 T}$, and outputs a trajectory consisting of $\frac{T}{2}$ rewards of value $.5 \times 2^{-T}$ and $\frac{T}{2}$ rewards of value $0$ w.p. $1 - 2^{-100 T}$.  For $E \in {\mathcal S}^t$ for $t \geq \frac{T}{2} + 1$, on any input $E$ for which the reward at time $\frac{T}{2} + 1$ equals 1 both simulators w.p.1 output a trajectory consisting of $\frac{T}{2}$ rewards of value $.5 \times 2^{-T}$ and $\frac{T}{2}$ rewards of value $1$, while on any input $E$ for which the reward at time $\frac{T}{2} + 1$ equals 0 both simulators w.p.1 output a trajectory consisting of $\frac{T}{2}$ rewards of value $.5 \times 2^{-T}$ and $\frac{T}{2}$ rewards of value $0$.
\\\indent Now, suppose I tell the $\textrm{DM}$ that I am going to take one of the two simulators ($SIM^1_T$ or $SIM^0_T$), give it to the $\textrm{DM}$ without telling them which of the two simulators I selected, and ask the $\textrm{DM}$ to tell me which of the two simulators I gave them by using the given simulator at most $10^5 \times (T+1) \times \textrm{poly}(T)$ times (on whatever inputs they choose, possibly selected dynamically based on previous outputs of the simulator), where they can use whatever exogenous randomness and as much computation as they desire.  Furthermore, the $\textrm{DM}$ must determine the algorithm they will use to determine how to use the simulator (and any exogenous noise and computation) before I select the simulator, and I may use knowledge of their algorithm (although not the realization of any exogenous noise used by their algorithm) in selecting which simulator I give them.  As the $\textrm{DM}$ knows the problem is either $\textsc{MS}^1_T$ or $\textsc{MS}^0_T$, and thus knows that for $E \in {\mathcal S}^t$ for $t \geq \frac{T}{2} + 1$ the output of the simulator is deterministic (and the same for both simulators), we may w.l.o.g. assume the $\textrm{DM}$ make no such calls to the simulator (i.e. all calls they make to the simulator are for $E \in {\mathcal S}^t$ with $t = 1,\ldots,\frac{T}{2}$).  It is thus easily verified (from a straightforward union bound) that no matter which simulator the $\textrm{DM}$ has, all their (at most $10^5 \times (T+1) \times \textrm{poly}(T)$) calls to the simulator will output the same trajectory consisting of $\frac{T}{2}$ rewards of value $.5 \times 2^{-T}$ and $\frac{T}{2}$ rewards of value $0$, with probability at least $1 - 10^5 \times (T+1) \times \textrm{poly}(T) \times 2^{-T}$.  I may thus determine which simulator to give them by examining their algorithm, and selecting whichever simulator ($\textrm{SIM}^1_T$ or $\textrm{SIM}^0_T$) maximizes the probability the $\textrm{DM}$ indicates the wrong simulator conditional on all their simulator calls outputting the above same (essentially uninformative) trajectory.  Combining the above with another union bound, we conclude that the $\textrm{DM}$ cannot possibly output the correct answer (regarding the identity of the simulator) with probability greater than $
\frac{1}{2} \times \big( 1 - 10^5 \times (T+1) \times \textrm{poly}(T) \times 2^{-T} \big) + 10^5 \times (T+1) \times \textrm{poly}(T) \times 2^{-T} = \frac{1}{2} + \frac{1}{2} \times 10^5 \times (T+1) \times \textrm{poly}(T) \times 2^{-T}.$
\\\indent Now we connect the above back to Lemma\ \ref{noptas} to yield a contradiction, by showing the existence of such an algorithm $\textrm{A}$ would imply the $\textrm{DM}$ could achieve a much higher success probability.  It is easily verified that the optimal policy for $\textsc{MS}^1_T$ is to not stop at all in the first $\frac{T}{2}$ periods, and to stop in all subsequent periods, and thus for $\textsc{MS}^1_T$ one has $\textrm{OPT}_{\textsc{MS}} = \frac{T}{2} \times 2^{-T}$, while the optimal policy for $\textsc{MS}^0_T$ is to w.p.1 stop in the first $\frac{T}{2}$ periods, and for $\textsc{MS}^0_T$ one has $\textrm{OPT}_{\textsc{MS}} = \frac{T}{2} \times .5 \times 2^{-T}$.  Let $X^1_{\textrm{A},t}$ denote the (random) value output by $\textrm{A}$ at time $t$.  Then it follows from an argument very similar to that we used in our proof of Lemma\ \ref{stoptoomuch} of Appendix\ \ref{heuristicformalsec} (the details of which we omit here) that for Problem $\textsc{MS}^1_T$, 
$$\textrm{OPT}_{\textsc{MS}} - \mathbb{E}\big[ \sum_{E \in {\mathcal E}} \mu(E) Z^1(E) \textrm{A}^1(E) \big] \geq .5 \times \mathbb{E}[\sum_{t=1}^{\frac{T}{2}} X^1_{\textrm{A},t}] \times 2^{-T},$$ while for Problem $\textsc{MS}^0_T$, $$\textrm{OPT}_{\textsc{MS}} - \mathbb{E}\big[ \sum_{E \in {\mathcal E}} \mu(E) Z^1(E) \textrm{A}^1(E) \big] \geq .5 \times (\frac{T}{2} - \mathbb{E}[\sum_{t=1}^{\frac{T}{2}} X^1_{\textrm{A},t}]) \times 2^{-T} - \frac{T}{2} \times 2^{-100 T}.$$
Thus since for both $\textsc{MS}^1_T$ and $\textsc{MS}^0_T$ the optimal value is at most $\frac{T}{2} \times 2^{-T}$, it follows from the above analysis, the properties of algorithm $\textrm{A}$ (guaranteed by its contradicting Lemma\ \ref{noptas}), and a straightforward contradiction argument (the details of which we omit) that when run on Problem $\textsc{MS}^1_T$ (equivalently when given simulator $SIM^1_T$), it must hold that $T^{-1} \mathbb{E}[\sum_{t=1}^{\frac{T}{2}} X^1_{\textrm{A},t}] \leq .01;$ and when run on Problem $\textsc{MS}^0_T$ (equivalently when given simulator $SIM^0_T$), it must hold that $T^{-1} \mathbb{E}[\sum_{t=1}^{\frac{T}{2}} X^1_{\textrm{A},t}] \geq .24  - 2^{-99 T}.$ 
\\\indent It then follows from a straightforward application of Hoeffding's inequality that if the problem is $\textsc{MS}^1_T$ and the $\textrm{DM}$ independently runs $\textrm{A}$ $10^5$ times (where in each run the $\textrm{DM}$ first uses their given simulator to independently generate a single trajectory $S$ which they then input sequentially into $\textrm{A}$), which will in total use the simulator at most $10^5 \times (T+1) \times \textrm{poly}(T)$ times, the following holds.  Letting $X^{1,j}_{\textrm{A},t}$ be the corresponding quantity in run $j$,
$$\mathbb{P}\big( 10^{-5} \sum_{j=1}^{10^5} T^{-1} \sum_{t=1}^{\frac{T}{2}} X^{1,j}_{\textrm{A},t} \leq .11 \big) \geq .9.$$
Alternatively, if the problem is $\textsc{MS}^0_T$ (and the $\textrm{DM}$ acts similarly) then
$$\mathbb{P}\big( 10^{-5} \sum_{j=1}^{10^5} T^{-1} \sum_{t=1}^{\frac{T}{2}} X^{1,j}_{\textrm{A},t} \geq .13 \big) \geq .9.$$
The above thus defines a procedure by which, given an unknown simulator (either $SIM^1_T$ or $SIM^0_T$), the DM can with probability at least $.9$ correctly identify which simulator they were given (by running $\textrm{A}$ $10^5$ times independently on $10^5$ independently generated trajectories, averaging $T^{-1} \sum_{t=1}^{\frac{T}{2}} X^{1,j}_{\textrm{A},t}$ across all $10^5$ trials, and outputting $SIM^1_T$ if this empirical average is at most .11, while outputting $SIM^0_T$ otherwise).  However, since our assumptions on $T_0$ imply that $\frac{1}{2} \times 10^5 \times (T+1) \times \textrm{poly}(T) \times 2^{-T} \leq .01$, it follows from our previous analysis that the $\textrm{DM}$ cannot achieve a success probability greater than $.51$.  Thus a contradiction is reached, completing the proof.  $\qed$
\end{proof}
\ \indent To close this section, let us comment a bit further on some prior works which achieve a PTAS for related problems by exploiting special structure and assumptions that do not hold in our model.  Several such works, including \cite{duttingetal2023} and \cite{agrawaletal2014}, exploit the assumption that the columns of the associated LP are randomly permuted.  This introduces many additional symmetries and simplifications that do not hold in our model.  The work of \cite{goldberg2018polynomial} utilizes certain special structural properties of optimal stopping (which do not extend to multiple stopping even though both problems have only one resource), and (for their PTAS-type results, but not their results with additive error) must (1) make additional assumptions regarding the variance of rewards, and (2) assume that the DM is given the first two moments of the maximum reward (to implement a certain truncation) without having to estimate them from samples (where we note this second point is implicit but follows from a careful reading of their algorithm and setup).  We refer the reader to \cite{chen2021} for an exploration of how direct extensions of the results of \cite{goldberg2018polynomial} to the more general set of problems we consider here (including multiple stopping) lead to an exponential, not polynomial, time algorithm.

\subsection{Additional discussion of assumptions}\label{appassumptionsec}
Assumption\ \ref{assumption1} is common in the literature on NRM and online combinatorial optimization. Regarding Assumption \ \ref{assumption2}, the NRM literature typically assumes a fixed number of resources that does not scale with $T$, in which case we may simply take $L = m$.  In addition, several works in $\textrm{NRM}$ further assume such a column sparsity (equivalently that each product requires the usage of at most $L$ resources), with this parameter appearing e.g. in the approximation guarantees of various algorithms.  In contrast, for online matching and independent set, $m$ will correspond (roughly) to the number of nodes or edges (depending on the particular problem), and scales with the graph size.  Here Assumption \ \ref{assumption2} requires the graph to be bounded-degree, again a common assumption in the literature.  We note that Assumption \ref{assumption3} still allows for many of the $a^r_{i,t}$ values to be zero (as will be the case in e.g. maximum matching), but precludes them being in $(0,\iota)$.  Furthermore, such an assumption is natural for problems in which the relevant marginal distributions are discrete (i.e. finite-support), as is the case in our setting (since we assume the relevant probability space is finite).  Let us point out that ${\mathcal C}$ in Assumption \ \ref{assumption4} captures the complexity of simulating the underlying process, and may in general depend on $T$ and $D$ (\cite{chen2021}). Our algorithmic runtimes all scale linearly with ${\mathcal C}$.  There are otherwise no assumptions made on the o.i.p.

\subsection{Proof of Claim\ \ref{convclaim2a}}\label{convclaim2aproofsec}
\begin{proof}{Proof of Claim\ \ref{convclaim2a} :}
As the results of \cite{schmidtetal2011} hold for general sequences of errors in the gradient calculations, and as our methods have no errors in the relevant proximal step (interpreting our projection step as a specialized proximal step), it follows directly from \cite{schmidtetal2011} Proposition 2 (after translating our concave maximization to a corresponding convex minimization) that w.p.1 
\begin{equation}\label{toboundappendix1}
\frac{1}{2} (k+1)^2 \alpha \big( \textrm{OPT}_{\textsc{pen}^{\mu,\theta}} - f^{\mu,\theta}({\mathcal X}^{\mu,k}) \big)
\end{equation}
 is at most
\small
\begin{equation}\label{toboundappendix1b}\bigg( ||{\mathcal X}^{\mu,0} - X_{*,\textsc{pen}^{\mu,\theta}}|| + 2 \sum_{i=1}^k i \alpha \big|\big|
\nabla f^{\mu,\theta}\big( (1 + \beta_{i-1}) {\mathcal X}^{\mu,i-1} - \beta_{i-1} {\mathcal X}^{\mu,i-2} \big) - \hat{G}^{\mu,i-1}\big( (1 + \beta_{i-1}) {\mathcal X}^{\mu,i-1} - \beta_{i-1} {\mathcal X}^{\mu,i-2} \big)\big|\big| \bigg)^2.
\end{equation}
\normalsize
Since (1) $X \in {\mathcal Q}^3$ implies
$$||X||^2 = \sum_{E \in {\mathcal E}} \sum_{r = 0}^{q-1} \big( X^r(E) \big)^2 \leq \sum_{E \in {\mathcal E}} \big(\sum_{r = 0}^{q-1} X^r(E)\big)^2 \leq \sum_{E \in {\mathcal E}}\big(\sqrt{\mu(E)}\big)^2 = T,$$
$(2)\ (a-b)^2 \leq a^2 + b^2$ for $a,b \geq 0$, and (3) ${\mathcal X}^{\mu,0}$ and $X_{*,\textsc{pen}^{\mu,\theta}}$ both belong to ${\mathcal Q}^3$, it follows that $||{\mathcal X}^{\mu,0} - X_{*,\textsc{pen}^{\mu,\theta}}||^2 \leq 2 T$.  Combining with the fact that $(a + b)^2 \leq 2(a^2 +  b^2)$ for all $a,b$, we conclude (from (\ref{toboundappendix1b})) that (\ref{toboundappendix1}) 
is at most
$$8 T + 8 \alpha^2 \bigg( \sum_{i=1}^k i \big|\big|
\nabla f^{\mu,\theta}\big( (1 + \beta_{i-1}) {\mathcal X}^{\mu,i-1} - \beta_{i-1} {\mathcal X}^{\mu,i-2} \big) - \hat{G}^{\mu,i-1}\big( (1 + \beta_{i-1}) {\mathcal X}^{\mu,i-1} - \beta_{i-1} {\mathcal X}^{\mu,i-2} \big)\big|\big| \bigg)^2.$$
By Cauchy-Schwarz, it follows that w.p.1 (\ref{toboundappendix1}) is at most
$$8 T + 8 \alpha^2  k \sum_{i=1}^k \bigg(i \big|\big|
\nabla f^{\mu,\theta}\big( (1 + \beta_{i-1}) {\mathcal X}^{\mu,i-1} - \beta_{i-1} {\mathcal X}^{\mu,i-2} \big) - \hat{G}^{\mu,i-1}\big( (1 + \beta_{i-1}) {\mathcal X}^{\mu,i-1} - \beta_{i-1} {\mathcal X}^{\mu,i-2} \big)\big|\big| \bigg)^2,$$
which (by bounding the $i$ appearing in the sum by its largest value $k$) is itself at most
$$8 T + 8 \alpha^2  k^3 \sum_{i=1}^k \big|\big|
\nabla f^{\mu,\theta}\big( (1 + \beta_{i-1}) {\mathcal X}^{\mu,i-1} - \beta_{i-1} {\mathcal X}^{\mu,i-2} \big) - \hat{G}^{\mu,i-1}\big( (1 + \beta_{i-1}) {\mathcal X}^{\mu,i-1} - \beta_{i-1} {\mathcal X}^{\mu,i-2} \big)\big|\big|^2.$$
Taking expectations on both sides, and defining $\chi_i \stackrel{\Delta}{=} (1 + \beta_i) {\mathcal X}^{\mu,i} - \beta_{i} {\mathcal X}^{\mu,i-1}$, we conclude that $$\textrm{OPT}_{\textsc{pen}^{\mu,\theta}} - \mathbb{E}\big[f^{\mu,\theta}({\mathcal X}^{\mu,k})\big] \leq \frac{2}{\alpha (k+1)^2} \big( 8 T + 8 \alpha^2  k^3 \sum_{i=1}^k \mathbb{E}\big[ \big|\big|\nabla f^{\mu,\theta}(\chi_{i-1}) - \hat{G}^{\mu,i-1}(\chi_{i-1})\big|\big|^2 \big] \big).$$  Note that for each $i \in \lbrace 0,\ldots,k \rbrace$, $\nabla f^{\mu,\theta}(\cdot) - \hat{G}^{\mu,i}(\cdot)$ is a random function, and that the number of different (deterministic) functions that $\nabla f^{\mu,\theta}(\cdot) - \hat{G}^{\mu,i}(\cdot)$ may equal is finite.  Let us denote this set of possible functions that $\nabla f^{\mu,\theta}(\cdot) - \hat{G}^{\mu,i}(\cdot)$ may equal by $F_i$.  Let ${\mathcal Z}_i$ denote the support of $\chi_i$ (i.e. the set of all $q \times |{\mathcal E}|$-dimensional vectors in the support of $\chi_i$), where it is easily verified that ${\mathcal Z}_i$ is finite as well.  It follows from the independence of $\lbrace {\mathcal S}^{E,i}, E \in {\mathcal E} \rbrace$ and $\aleph^i$ from $\chi_i$ for each $i$ that for all $\digamma \in F_{i-1}$ and $\chi' \in {\mathcal Z}_{i-1}$, $\mathbb{P}\big( \nabla f^{\mu,\theta}(\cdot) - \hat{G}^{\mu,i-1}(\cdot) = \digamma, \chi_{i-1} = \chi' \big) = \mathbb{P}\big( \nabla f^{\mu,\theta}(\cdot) - \hat{G}^{\mu,i-1}(\cdot) = \digamma \big) \times \mathbb{P}(\chi_{i-1} = \chi')$.  Thus for each $i \in \lbrace 1,\ldots,k \rbrace$, $\mathbb{E}\big[ \big|\big|\nabla f^{\mu,\theta}(\chi_{i-1}) - \hat{G}^{\mu,i-1}(\chi_{i-1})\big|\big|^2 \big]$ equals
$$\sum_{\digamma \in F_{i-1}} \sum_{\chi' \in {\mathcal Z}_{i-1}} ||\digamma(\chi')||^2 \mathbb{P}\big( \nabla f^{\mu,\theta}(\cdot) - \hat{G}^{\mu,i-1}(\cdot) = \digamma \big) \times \mathbb{P}(\chi_{i-1} = \chi'),$$
itself equal to 
$$\sum_{\chi' \in {\mathcal Z}_{i-1}} \mathbb{P}(\chi_{i-1} = \chi') \sum_{\digamma \in F_{i-1}} \mathbb{P}\big( \nabla f^{\mu,\theta}(\cdot) - \hat{G}^{\mu,i-1}(\cdot) = \digamma \big) ||\digamma(\chi')||^2,$$
itself at most
$$\max_{\chi' \in {\mathcal Z}_{i-1}} \sum_{\digamma \in F_{i-1}} \mathbb{P}\big( \nabla f^{\mu,\theta}(\cdot) - \hat{G}^{\mu,i-1}(\cdot) = \digamma \big) ||\digamma(\chi')||^2.$$
Next, we observe that (\ref{inQ4a}) implies ${\mathcal X}^{\mu,j} \in {\mathcal Q}^4(0,1)$ for all $j$. As it is easily verified that $\beta \in [0,1]$ and $X,Y \in {\mathcal Q}^4(0,1)$ implies $(1 + \beta) X - \beta Y \in {\mathcal Q}^4(-1,2)$, we conclude that w.p.1 $\chi_{i-1} \in {\mathcal Q}^4(-1,2)$ for all $i$.  Thus for each $i \in \lbrace 1,\ldots,k \rbrace$, $\mathbb{E}\big[ \big|\big|\nabla f^{\mu,\theta}(\chi_{i-1}) - \hat{G}^{\mu,i-1}(\chi_{i-1})\big|\big|^2 \big]$ is at most
$$\sup_{\chi' \in {\mathcal Q}^4(-1,2)} \sum_{\digamma \in F_{i-1}} \mathbb{P}\big( \nabla f^{\mu,\theta}(\cdot) - \hat{G}^{\mu,i-1}(\cdot) = \digamma \big) ||\digamma(\chi')||^2,$$
itself equal to $\sup_{\chi' \in {\mathcal Q}^4(-1,2)} \mathbb{E}\big[ \big|\big|\nabla f^{\mu,\theta}(\chi') - \hat{G}^{\mu,i-1}(\chi')\big|\big|^2\big]$.  Combining the above with the fact that $\hat{G}^{\mu,j}(\cdot)$ has the same distribution (as a random function) for all $j$ and some straightforward algebra completes the proof.  $\qed$
\end{proof}

\subsection{Proof of Corollary\ \ref{convclaim2}}\label{convclaim2proofsec}
To prove Corollary\ \ref{convclaim2}, we bound the suprema appearing in the convergence guarantees of Claims\ \ref{convclaim1} and \ref{convclaim2a}, and plug our bounds back into those claims to complete the proof.  First, let us prove an auxiliary concentration result, bounding the moments of the difference of two sums which arise naturally in our analysis (one a noisy approximation of the other) using Hoeffding's inequality.  

\begin{claim}\label{concentrateclaim1}
\small{
$$\forall\ a < b, \sup_{X \in {\mathcal Q}^4(a,b), S \in {\mathcal S}, i \in \lbrace 1,\ldots,m \rbrace}  \mathbb{E}\bigg[ \bigg( \sum_{t = 1}^T \sum_{r=1}^{q-1} a^r_{i}(S_{[t]}) \frac{X^r(S_{[t]})}{\sqrt{\mu(S_{[t]})}} - \frac{T}{\eta_2} \sum_{t \in \aleph^1} \sum_{r=1}^{q-1}a^r_{i}(S_{[t]}) \frac{X^r(S_{[t]})}{\sqrt{\mu(S_{[t]})}} \bigg)^2 \bigg] \leq  (b-a)^2 \eta_2^{-1} T^2.$$}
\end{claim}
\begin{proof}{Proof :} By the tail integral formula for higher moments (itself following from integration by parts), we have that
$$\mathbb{E}\bigg[ \bigg( \sum_{t = 1}^T \sum_{r=1}^{q-1} a^r_{i}(S_{[t]}) \frac{X^r(S_{[t]})}{\sqrt{\mu(S_{[t]})}} - \frac{T}{\eta_2} \sum_{t \in \aleph^1} \sum_{r=1}^{q-1}a^r_{i}(S_{[t]}) \frac{X^r(S_{[t]})}{\sqrt{\mu(S_{[t]})}} \bigg)^2 \bigg]\ \textrm{equals}$$
$$2 \int_0^{\infty} x \mathbb{P}\bigg( \bigg| \sum_{t = 1}^T \sum_{r=1}^{q-1} a^r_{i}(S_{[t]}) \frac{X^r(S_{[t]})}{\sqrt{\mu(S_{[t]})}} - \frac{T}{\eta_2} \sum_{t \in \aleph^1} \sum_{r=1}^{q-1}a^r_{i}(S_{[t]}) \frac{X^r(S_{[t]})}{\sqrt{\mu(S_{[t]})}} \bigg| > x \bigg) dx.$$
Let us apply Hoeffding's inequality (which is also valid for sampling without replacement), which (since $\sum_{r=1}^{q-1}a^r_{i}(S_{[t]}) \frac{X^r(S_{[t]})}{\sqrt{\mu(S_{[t]})}} \in [a,b]$) implies that 
$$
 \mathbb{P}\bigg( \bigg| \sum_{t = 1}^T \sum_{r=1}^{q-1} a^r_{i}(S_{[t]}) \frac{X^r(S_{[t]})}{\sqrt{\mu(S_{[t]})}} - \frac{T}{\eta_2} \sum_{t \in \aleph^1} \sum_{r=1}^{q-1}a^r_{i}(S_{[t]}) \frac{X^r(S_{[t]})}{\sqrt{\mu(S_{[t]})}} \bigg| > x \bigg) \leq 2 \exp\big( - \frac{2}{(b-a)^2} \eta_2 T^{-2} x^2 \big).$$
Combining with the fact that (by standard calculus arguments) $\int_0^{\infty} x \exp\big( - \frac{2}{(b-a)^2} \eta_2 T^{-2} x^2 \big) dx = \frac{1}{2} \big( \frac{2}{(b-a)^2} \eta_2 T^{-2} \big)^{-1}$ and applying some straightforward algebra completes the proof. $\qed$
\end{proof}

Next, let us apply the above to bound the suprema of interest.

\begin{claim}\label{boundnorm1}
$$
\sup_{X \in {\mathcal Q}^4(-1,2)} 
\mathbb{E}\Bigg[ \bigg|\bigg| \nabla f^{\mu,\theta}(X) - 
\hat{G}^{\mu,1}(X)\bigg|\bigg|^2 \Bigg] \leq 72 \iota^{-2} L^2 \theta^{-2} T^3 \eta_2^{-1} q + 8 \iota^{-2} \eta_1^{-1} L^2 T q.
$$
If $\eta_2 = T$ (i.e. the relevant sum is computed exactly), then
$$\sup_{X \in {\mathcal Q}^4(-1,2)} 
\mathbb{E}\Bigg[ \bigg|\bigg| \nabla f^{\mu,\theta}(X) - 
\hat{G}^{\mu,1}(X)\bigg|\bigg|^2 \Bigg] \leq 8 \iota^{-2} \eta_1^{-1} L^2 T q.
$$
\end{claim}
\begin{proof}{Proof :} Let us fix $X \in {\mathcal Q}^4(-1,2)$.  First suppose $\eta_2 < T$.  Let us fix $E \in {\mathcal E}, l \in \lbrace 0,\ldots,q-1 \rbrace$.  Then $\nabla f^{\mu,\theta}(X)_{E,l} - \hat{G}^{\mu,1}(X)_{E,l}$ equals
$$\sqrt{\mu(E)} \Bigg( Z^l(E) - 2 \iota^{-1} \sum_{S' \in {\mathcal S} : E \subseteq S'} \frac{\mu(S')}{\mu(E)} \sum_{i=1}^m a^l_i(E) \phi'_{\theta}\bigg(\sum_{t=1}^T \sum_{r=1}^{q-1} a^r_i(S'_{[t]}) \frac{X^r(S'_{[t]})}{\sqrt{\mu(S'_{[t]})}} - b_i\bigg) \Bigg)$$
$$- \sqrt{\mu(E)} \Bigg( Z^l(E) - 2 \iota^{-1} \sum_{i=1}^m a^l_{i}(E) \times \eta^{-1}_1 \sum_{S' \in {\mathcal S}^{E,1}}\phi'_{\theta}\bigg( \frac{T}{\eta_2} \sum_{t \in \aleph^1} \sum_{r=1}^{q-1} a^r_{i}(S'_{[t]}) \frac{X^r(S'_{[t]})}{\sqrt{\mu(S'_{[t]})}} - b_i \bigg) \Bigg).$$
It follows that $\bigg| \nabla f^{\mu,\theta}(\overline{X})_{E,l} - \hat{G}^{\mu,1}(\overline{X})_{E,l} \bigg|$ is at most
\begin{equation}\label{makeeasier1}2 \iota^{-1} \sqrt{\mu(E)} \bigg| \sum_{S' \in {\mathcal S} : E \subseteq S'} \frac{\mu(S')}{\mu(E)} \sum_{i=1}^m a^l_i(E) \phi'_{\theta}\bigg(\sum_{t=1}^T \sum_{r=1}^{q-1} a^r_i(S'_{[t]}) \frac{X^r(S'_{[t]})}{\sqrt{\mu(S'_{[t]})}} - b_i\bigg)
\end{equation}
\begin{equation}\label{makeeasier2}
\ \ \ \ \ \ \ \ \ \ \ \ \ \ \ \ \ \ \ \ \ \ \ \ \ - \eta^{-1}_1 \sum_{S' \in {\mathcal S}^{E,1}} \sum_{i=1}^m a^l_{i}(E) \phi'_{\theta}\bigg( \frac{T}{\eta_2} \sum_{t \in \aleph^1} \sum_{r=1}^{q-1} a^r_{i}(S'_{[t]}) \frac{X^r(S'_{[t]})}{\sqrt{\mu(S'_{[t]})}} - b_i \bigg) \bigg|.
\end{equation}
Notice that (\ref{makeeasier2}) has two sources of randomization that (\ref{makeeasier1}) does not : the randomization over trajectories in ${\mathcal S}^{E,1}$ and the randomization over time indices in $\aleph^1$.  Thus to make the analysis simpler, we add and subtract an intermediate term to reduce bounding the above difference to bounding two differences, but where the two terms in each difference differ in only one source of randomization.  Thus (adding and subtracting $\eta^{-1}_1 \sum_{S' \in {\mathcal S}^{E,1}} \sum_{i=1}^m a^l_{i}(E) \phi'_{\theta}\bigg(\sum_{t=1}^T \sum_{r=1}^{q-1} a^r_i(S'_{[t]}) \frac{X^r(S'_{[t]})}{\sqrt{\mu(S'_{[t]})}} - b_i\bigg) $, applying the triangle inequality, and using the fact that $(a+b)^2 \leq 2 (a^2 + b^2)$ along with linearity of expectation) we conclude that for any fixed $X,E,l$, $\mathbb{E}\bigg[\bigg( \nabla f^{\mu,\theta}(X)_{E,l} - \hat{G}^{\mu,1}(X)_{E,l} \bigg)^2\bigg]$ is at most
\small
$$8 \iota^{-2} \mu(E) \mathbb{E}\bigg[\bigg( \sum_{S' \in {\mathcal S} : E \subseteq S'} \frac{\mu(S')}{\mu(E)} \sum_{i=1}^m a^l_i(E) \phi'_{\theta}\bigg(\sum_{t=1}^T \sum_{r=1}^{q-1} a^r_i(S'_{[t]}) \frac{X^r(S'_{[t]})}{\sqrt{\mu(S'_{[t]})}} - b_i\bigg) 
$$
\begin{equation}\label{boundnorm1eq3}
\ \ \ \ \ \ \ \ \ \ \ \ \ \ \ \ \ \ \ \ \ \ \ \ \ \ \ \ \ \ \ \ \ \ \ - \eta^{-1}_1 \sum_{S' \in {\mathcal S}^{E,1}} \sum_{i=1}^m a^l_{i}(E) \phi'_{\theta}\bigg(\sum_{t=1}^T \sum_{r=1}^{q-1} a^r_i(S'_{[t]}) \frac{X^r(S'_{[t]})}{\sqrt{\mu(S'_{[t]})}} - b_i\bigg)  \bigg)^2\bigg] 
\end{equation}
$$+ 8 \iota^{-2} \mu(E) \mathbb{E}\bigg[\bigg( \eta^{-1}_1 \sum_{S' \in {\mathcal S}^{E,1}} \sum_{i=1}^m a^l_{i}(E) \phi'_{\theta}\bigg(\sum_{t=1}^T \sum_{r=1}^{q-1} a^r_i(S'_{[t]}) \frac{X^r(S'_{[t]})}{\sqrt{\mu(S'_{[t]})}} - b_i\bigg)$$
\begin{equation}\label{boundnorm1eq4}
\ \ \ \ \ \ \ \ \ \ \ \ \ \ \ \ \ \ \ \ \ \ \ \ \ \ \ \ \ \ \ \ \ \ \ \ \ \ \ \ - \eta^{-1}_1 \sum_{S' \in {\mathcal S}^{E,1}} \sum_{i=1}^m a^l_{i}(E) \phi'_{\theta}\bigg( \frac{T}{\eta_2} \sum_{t \in \aleph^1} \sum_{r=1}^{q-1} a^r_{i}(S'_{[t]}) \frac{X^r(S'_{[t]})}{\sqrt{\mu(S'_{[t]})}} - b_i \bigg)   \bigg)^2 \bigg] 
\end{equation}
\normalsize
We now bound (\ref{boundnorm1eq3}) and (\ref{boundnorm1eq4}), beginning with (\ref{boundnorm1eq3}).  To bound (\ref{boundnorm1eq3}), we first bound the tail of the difference (i.e. the r.v. within the expectation), and use that to bound the expectation.  Thus let us first bound (for $x > 0$ and $X \in {\mathcal Q}^4(-1,2)$)
\begin{eqnarray*}
\ &\ &\ (a) : \mathbb{P}\Bigg( \Bigg|
\sum_{S' \in {\mathcal S} : E \subseteq S'} \frac{\mu(S')}{\mu(E)} \sum_{i=1}^m a^l_i(E) \phi'_{\theta}\bigg(\sum_{t=1}^T \sum_{r=1}^{q-1} a^r_i(S'_{[t]}) \frac{X^r(S'_{[t]})}{\sqrt{\mu(S'_{[t]})}} - b_i\bigg) 
\\&\ &\ \ \ \ \ \ \ \ \ \ \ \ \ \ \ - \eta^{-1}_1 \sum_{S' \in {\mathcal S}^{E,1}} \sum_{i=1}^m a^l_{i}(E) \phi'_{\theta}\bigg(\sum_{t=1}^T \sum_{r=1}^{q-1} a^r_i(S'_{[t]}) \frac{X^r(S'_{[t]})}{\sqrt{\mu(S'_{[t]})}} - b_i\bigg) \Bigg| > x \Bigg).
\end{eqnarray*}
We bound this quantity using Hoeffding's inequality.  Indeed, observe that $\eta^{-1}_1 \sum_{S' \in {\mathcal S}^{E,1}} \sum_{i=1}^m a^l_{i}(E) \phi'_{\theta}\bigg(\sum_{t=1}^T \sum_{r=1}^{q-1}a^r_i(S'_{[t]}) \frac{X^r(S'_{[t]})}{\sqrt{\mu(S'_{[t]})}} - b_i\bigg)$ is the average of $\eta_1$ i.i.d. r.v.s, each of which lies in $[0, L]$ (by Claim\ \ref{huberprops} and the definition of $L$), and each of which has expected value $\sum_{S' \in {\mathcal S} : E \subseteq S'} \frac{\mu(S')}{\mu(E)} \sum_{i=1}^m a^l_i(E) \phi'_{\theta}\bigg(\sum_{t=1}^T \sum_{r=1}^{q-1} a^r_i(S'_{[t]}) \frac{X^r(S'_{[t]})}{\sqrt{\mu(S'_{[t]})}} - b_i\bigg)$.  We may thus apply Hoeffding's inequality, and conclude that (a) is at most $2 \exp\big( - 2 \eta_1 L^{-2} x^2 \big)$.  It then follows from the tail integral formula for higher moments that (\ref{boundnorm1eq3}) is at most $32 \iota^{-2} \mu(E) \int_0^{\infty} x \exp\big( - 2 \eta_1 L^{-2} x^2 \big) dx$, which by some straightforward calculus equals $8 \iota^{-2} \mu(E) \eta_1^{-1} L^2$.  We conclude that (\ref{boundnorm1eq3}) is at most $8 \iota^{-2} \mu(E) \eta_1^{-1} L^2$.
\\\indent We next bound (\ref{boundnorm1eq4}).  To bound (\ref{boundnorm1eq4}), let us bound 
$$
(b) : \mathbb{E}\bigg[\Bigg( \eta^{-1}_1 \sum_{S' \in {\mathcal S}^{E,1}} \sum_{i=1}^m a^l_{i}(E) \phi'_{\theta}\bigg(\sum_{t=1}^T \sum_{r=1}^{q-1} a^r_i(S'_{[t]}) \frac{X^r(S'_{[t]})}{\sqrt{\mu(S'_{[t]})}} - b_i\bigg)$$
$$\ \ \ \ \ \ \ \ \ \ \ \ \ \ \ \ \ \ \ \ \ \ \ \ \ \ \ \ \ \ - \eta^{-1}_1 \sum_{S' \in {\mathcal S}^{E,1}} \sum_{i=1}^m a^l_{i}(E) \phi'_{\theta}\bigg( \frac{T}{\eta_2} \sum_{t \in \aleph^1} \sum_{r=1}^{q-1} a^r_{i}(S'_{[t]}) \frac{X^r(S'_{[t]})}{\sqrt{\mu(S'_{[t]})}} - b_i \bigg)   \Bigg)^2 \bigg].
$$
We proceed by : (1) conditioning on ${\mathcal S}^{E,1}$ to remove one layer of randomness, (2) using the triangle inequality and fact that $a^l_i(E) \in [0,1]$ to further reduce the problem to bounding the (conditional) expectation of the square of a double summation of several difference terms (one for each $S' \in {\mathcal S}^{E,1}$ and $i \in a^{l,+}(E)$), (3) using Cauchy-Schwarz to bound that term by the (conditional) expectation of a double summation of several individual squared difference terms (at the cost of a multiplicative $\eta_1 L$ corresponding to the number of terms in the double sum), (4) bounding this expression by taking a supremum over the (conditional) expectation of the square of any one of those individual difference terms (at the cost of another multiplicative $\eta_1 L$ again corresponding to the number of terms in the double sum), using the Lipschitz continuity of $\phi'_{\theta}$ (which appears in this difference term) to further bound the expectation of the square of any one of those individual difference terms (at the cost of a multiplicative $\theta^{-2}$ coming from the Lipschitz constant), and (5) recognizing the resulting conditional expectation of a squared difference term can be bounded using Claim\ \ref{concentrateclaim1}.
\\\indent Let $\sigma({\mathcal S}^{E,1})$ denote the $\sigma$-field generated by ${\mathcal S}^{E,1}$.  Then by the triangle inequality and fact that $a^l_i(E) \in [0,1]$, (b) is at most
\small
$$\eta_1^{-2} \mathbb{E}\Bigg[ \mathbb{E}\bigg[\Bigg( \sum_{S' \in {\mathcal S}^{E,1}} \sum_{i \in a^{l,+}(E)} \bigg| \phi'_{\theta}\bigg(\sum_{t=1}^T \sum_{r=1}^{q-1} a^r_i(S'_{[t]}) \frac{X^r(S'_{[t]})}{\sqrt{\mu(S'_{[t]})}} - b_i\bigg) - \phi'_{\theta}\bigg( \frac{T}{\eta_2} \sum_{t \in \aleph^1} \sum_{r=1}^{q-1} a^r_{i}(S'_{[t]}) \frac{X^r(S'_{[t]})}{\sqrt{\mu(S'_{[t]})}} - b_i \bigg) \bigg|   \Bigg)^2 \Bigg| \sigma({\mathcal S}^{E,1}) \bigg] \Bigg],
$$
\normalsize
which by Cauchy-Schwarz is at most
\small
$$
\eta_1^{-1} L \times \mathbb{E}\Bigg[ \mathbb{E}\bigg[ \sum_{S' \in {\mathcal S}^{E,1}} \sum_{i \in a^{l,+}(E)} \Bigg(\phi'_{\theta}\bigg(\sum_{t=1}^T \sum_{r=1}^{q-1} a^r_i(S'_{[t]}) \frac{X^r(S'_{[t]})}{\sqrt{\mu(S'_{[t]})}} - b_i\bigg) - \phi'_{\theta}\bigg( \frac{T}{\eta_2} \sum_{t \in \aleph^1} \sum_{r=1}^{q-1} a^r_{i}(S'_{[t]}) \frac{X^r(S'_{[t]})}{\sqrt{\mu(S'_{[t]})}} - b_i \bigg) \Bigg)^2  \Bigg| \sigma({\mathcal S}^{E,1}) \bigg] \Bigg].
$$
\normalsize
It then follows from definitions, some straightforward reasoning about conditional expectations, and the fact that the relevant double summation has at most $\eta_1 L$ terms (where $\eta_1^{-1} L \times \eta_1 L = L^2$) that (b) is at most
\small
$$
L^2 \times \sup_{X \in {\mathcal Q}^4(-1,2), S' \in {\mathcal S}, i \in \lbrace 1,\ldots,m \rbrace} \mathbb{E}\Bigg[ \Bigg(\phi'_{\theta}\bigg(\sum_{t=1}^T \sum_{r=1}^{q-1} a^r_i(S'_{[t]}) \frac{X^r(S'_{[t]})}{\sqrt{\mu(S'_{[t]})}} - b_i\bigg) - \phi'_{\theta}\bigg( \frac{T}{\eta_2} \sum_{t \in \aleph^1} \sum_{r=1}^{q-1} a^r_{i}(S'_{[t]}) \frac{X^r(S'_{[t]})}{\sqrt{\mu(S'_{[t]})}} - b_i \bigg) \Bigg)^2 \Bigg].
$$
\normalsize
Applying Claim\ \ref{huberprops} (in particular the fact that $\phi'_{\theta}$ is $\theta^{-1}$-Lipschitz), we conclude that
(b) is at most
$$
L^2 \theta^{-2} \times \sup_{X \in {\mathcal Q}^4(-1,2), S' \in {\mathcal S}, i \in \lbrace 1,\ldots,m \rbrace} \mathbb{E}\Bigg[ \Bigg(\sum_{t=1}^T \sum_{r=1}^{q-1} a^r_i(S'_{[t]}) \frac{X^r(S'_{[t]})}{\sqrt{\mu(S'_{[t]})}} - \frac{T}{\eta_2} \sum_{t \in \aleph^1} \sum_{r=1}^{q-1} a^r_{i}(S'_{[t]}) \frac{X^r(S'_{[t]})}{\sqrt{\mu(S'_{[t]})}}  \Bigg)^2 \Bigg],
$$
which by Claim\ \ref{concentrateclaim1} is at most $9 L^2 \theta^{-2} T^2 \eta_2^{-1}$, and we conclude that (\ref{boundnorm1eq4}) is at most $72 \iota^{-2} \mu(E) L^2 \theta^{-2} T^2 \eta_2^{-1}$.
\\\indent Combining our bounds for (\ref{boundnorm1eq3}) and (\ref{boundnorm1eq4}) with the definition of $||X||$ and fact that $\sum_{E \in {\mathcal E}} \mu(E) = T$ completes the proof.  The only difference when $\eta_2 = T$ is that the error term $72 \iota^{-2} \mu(E) L^2 \theta^{-2} T^2 \eta_2^{-1}$ vanishes, and that case thus follows from a nearly identical argument. $\qed$
\end{proof}
\ \indent Next, we prove a bound on $L^{\mu,\theta}$.

\begin{claim}\label{smoothbound1}
$L^{\mu,\theta} \leq 2 \iota^{-1} \theta^{-1} \sqrt{ U_r L W }.$
\end{claim}
Our proof proceeds by simplifying the relevant expressions using the triangle inequality, Cauchy-Schwarz, and Jensen's inequality, and then recognizing that a key term can be bounded by $W$. 
\begin{proof}{Proof of Claim\ \ref{smoothbound1}} 
By definition, the desired statement is equivalent to the statement that $\forall\ X,Y \in {\mathcal R}^{q \times |{\mathcal E}|},$
\begin{equation}\label{toshowsmooth1}
\sum_{E \in {\mathcal E} , l \in \lbrace 0,\ldots,q - 1 \rbrace} \big( \nabla f^{\mu,\theta}(X)_{E,l} - \nabla f^{\mu,\theta}(Y)_{E,l} \big)^2 \leq 4 \iota^{-2} \theta^{-2} U_r L W \sum_{E \in {\mathcal E}, l \in \lbrace 0,\ldots,q - 1 \rbrace} \big( X^l(E) - Y^l(E) \big)^2.
\end{equation}
Let us fix $E \in {\mathcal E}$ and $l \in \lbrace 0,\ldots,q - 1 \rbrace$, and examine $\frac{\iota^2}{4 \mu(E)}\big( \nabla f^{\mu,\theta}(X)_{E,l} - \nabla f^{\mu,\theta}(Y)_{E,l} \big)^2$, which by Claim\ \ref{whatisderiv1}, definitions, and some straightforward algebra equals
{\small
$$\Bigg( \sum_{S' \in {\mathcal S} : E \subseteq S'} \frac{\mu(S')}{\mu(E)} \sum_{i \in a^{l,+}(E)} a^l_i(E)\bigg( \phi'_{\theta} \big( \sum_{(r,t) \in {\mathcal RT}_i(S')} a^r_i(S'_{[t]}) \frac{X^r(S'_{[t]})}{\sqrt{\mu(S'_{[t]})}} - b_i \big)
- \phi'_{\theta} \big(  \sum_{(r,t) \in {\mathcal RT}_i(S')} a^r_i(S'_{[t]}) \frac{Y^r(S'_{[t]})}{\sqrt{\mu(S'_{[t]})}} - b_i \big) \bigg) \Bigg)^2,$$}
which by Claim\ \ref{huberprops} and the triangle inequality is at most
$$(a) : \theta^{-2} \Bigg( \sum_{S' \in {\mathcal S} : E \subseteq S'} \frac{\mu(S')}{\mu(E)} \sum_{i \in a^{l,+}(E)} \nonumber
\bigg| \sum_{(r,t) \in {\mathcal RT}_i(S')} \big( a^r_i(S'_{[t]}) \frac{X^r(S'_{[t]})}{\sqrt{\mu(S'_{[t]})}} - a^r_i(S'_{[t]}) \frac{Y^r(S'_{[t]})}{\sqrt{\mu(S'_{[t]})}} \big) \bigg| \Bigg)^2.
$$
By Cauchy-Schwarz, 
\small
$$\bigg| \sum_{(r,t) \in {\mathcal RT}_i(S')} \big( a^r_i(S'_{[t]}) \frac{X^r(S'_{[t]})}{\sqrt{\mu(S'_{[t]})}} - a^r_i(S'_{[t]}) \frac{Y^r(S'_{[t]})}{\sqrt{\mu(S'_{[t]})}} \big) \bigg|
\leq \sqrt{|{\mathcal RT}_i(S')|} \sqrt{ \sum_{(r,t) \in {\mathcal RT}_i(S')} \big( a^r_i(S'_{[t]}) \frac{X^r(S'_{[t]})}{\sqrt{\mu(S'_{[t]})}} - a^r_i(S'_{[t]}) \frac{Y^r(S'_{[t]})}{\sqrt{\mu(S'_{[t]})}} \big)^2},$$
\normalsize
and thus (also since $a^r_i(\cdot) \in [0,1]$) (a) is at most 
$$\theta^{-2} U_r \Bigg( \sum_{S' \in {\mathcal S} : E \subseteq S'} \frac{\mu(S')}{\mu(E)} \sum_{i \in a^{l,+}(E)} \sqrt{ \sum_{(r,t) \in {\mathcal RT}_i(S')} \big( \frac{X^r(S'_{[t]})}{\sqrt{\mu(S'_{[t]})}} - \frac{Y^r(S'_{[t]})}{\sqrt{\mu(S'_{[t]})}} \big)^2} \Bigg)^2.$$
By the fact that $\sum_{S' \in {\mathcal S} : E \subseteq S'} \frac{\mu(S')}{\mu(E)} = 1$ and Jensen's inequality, we conclude that (a) is at most
$$
\theta^{-2} U_r \sum_{S' \in {\mathcal S} : E \subseteq S'} \frac{\mu(S')}{\mu(E)} \Bigg(\sum_{i \in a^{l,+}(E)} \sqrt{ \sum_{(r,t) \in {\mathcal RT}_i(S')} \big( \frac{X^r(S'_{[t]})}{\sqrt{\mu(S'_{[t]})}} - \frac{Y^r(S'_{[t]})}{\sqrt{\mu(S'_{[t]})}} \big)^2} \Bigg)^2.
$$
Applying Cauchy-Schwarz to $\Bigg(\sum_{i \in a^{l,+}(E)} \sqrt{ \sum_{(r,t) \in {\mathcal RT}_i(S')} \big( \frac{X^r(S'_{[t]})}{\sqrt{\mu(S'_{[t]})}} - \frac{Y^r(S'_{[t]})}{\sqrt{\mu(S'_{[t]})}} \big)^2} \Bigg)^2$, we conclude that for all $E \in {\mathcal E}$ and $l \in \lbrace 0,\ldots,q-1 \rbrace$, $\frac{\iota^2}{4 \mu(E)}\big( \nabla f^{\mu,\theta}(X)_{E,l} - \nabla f^{\mu,\theta}(Y)_{E,l} \big)^2$ is at most 
$$
\theta^{-2} U_r L \sum_{S' \in {\mathcal S} : E \subseteq S'} \frac{\mu(S')}{\mu(E)} \sum_{i \in a^{l,+}(E)}  \sum_{(r,t) \in {\mathcal RT}_i(S')} \big( \frac{X^r(S'_{[t]})}{\sqrt{\mu(S'_{[t]})}} - \frac{Y^r(S'_{[t]})}{\sqrt{\mu(S'_{[t]})}} \big)^2.
$$
Combining the above, and (for $E' \in {\mathcal E}$) letting ${\mathcal D}(r,E')$ denote $\big( \frac{X^r(E')}{\sqrt{\mu(E')}} - \frac{Y^r(E')}{\sqrt{\mu(E')}} \big)^2,$ we conclude that $\sum_{E \in {\mathcal E}, l \in \lbrace 0,\ldots,q-1 \rbrace} \big( \nabla f^{\mu,\theta}(X)_{E,l} - \nabla f^{\mu,\theta}(Y)_{E,l} \big)^2$ is at most
\begin{equation}\label{step1aaa}
4 \iota^{-2} \theta^{-2} U_r L \sum_{E \in {\mathcal E}} \sum_{l=0}^{q-1} \sum_{S' \in {\mathcal S} : E \subseteq S'} \mu(S') \sum_{i \in a^{l,+}(E)}  \sum_{(r,t) \in {\mathcal RT}_i(S')} {\mathcal D}(r,S'_{[t]}).
\end{equation}
Notice that we may equivalently decompose (\ref{step1aaa}) as 
$$
\sum_{r = 0}^{q-1} \sum_{E' \in {\mathcal E}} C_{r,E'} {\mathcal D}(r,E')
$$
for some appropriate set of $C_{r,E'}$ (whose values may depend on the terms of (\ref{step1aaa}) in a complicated way).  We now bound bound these $C_{r,E'}$ by carefully analyzing exactly when ${\mathcal D}(r,E')$ appears in (\ref{step1aaa}), and with what prefactors.  Thus let us fix some $E' \in {\mathcal E}$, and suppose $E' \in {\mathcal S}^t$ for some $t \in \lbrace 1,\ldots,T \rbrace$.  Then it is easily verified that $C_{r,E'}$ equals
$$4 \iota^{-2} \theta^{-2} U_r L \sum_{E \in {\mathcal E}} \sum_{l=0}^{q-1} \sum_{S' \in {\mathcal S} : E \subseteq S'} \mu(S') \sum_{i \in a^{l,+}(E)}  I\big( (r,t) \in {\mathcal RT}_i(S')\ \textrm{and}\ S'_{[t]} = E' \big).
$$
By (1) observing that $(r,t) \in {\mathcal RT}_i(S')\ \textrm{and}\ S'_{[t]} = E'$ iff $i \in a^{r,+}(E')\ \textrm{and}\ S'_{[t]} = E'$, (2) observing that $\sum_{i \in a^{l,+}(E)} I\big(i \in a^{r,+}(E')\big) = |a^{l,+}(E) \bigcap a^{r,+}(E')|$, (3) interchanging the order of summation w.r.t. $l$ and $S'$, and (4) rewriting the requirement $E \subseteq S'$ directly into the sum using an indicator,
we conclude that $C_{r,E'}$ equals
$$
4 \iota^{-2} \theta^{-2} U_r L \sum_{E \in {\mathcal E}} \sum_{S' \in {\mathcal S}} I(E \subseteq S'\ \textrm{and}\ S'_{[t]} = E') \mu(S') \sum_{l=0}^{q-1}|a^{l,+}(E) \bigcap a^{r,+}(E')|.
$$
By interchanging the order of summation of $E$ and $S'$, and rewriting the requirements $S'_{[t]} = E'$ and $E \subseteq S'$ into the indices of summation, we conclude that $C_{r,E'}$ equals
$$
4 \iota^{-2} \theta^{-2} U_r L \sum_{S' \in {\mathcal S} : S'_{[t]} = E'} \mu(S') \sum_{E \in {\mathcal E} : E \subseteq S'}\sum_{l=0}^{q-1}|a^{l,+}(E) \bigcap a^{r,+}(E')|,
$$
itself equal to
$$
4 \iota^{-2} \theta^{-2} U_r L \sum_{S' \in {\mathcal S} : S'_{[t]} = E'} \mu(S') \sum_{u=1}^T\sum_{l=0}^{q-1}|a^{l,+}(S'_{[u]}) \bigcap a^{r,+}(S'_{[t]})|.
$$
Since by definition $\sum_{u=1}^T\sum_{l=0}^{q-1}|a^{l,+}(S'_{[u]}) \bigcap a^{r,+}(S'_{[t]})| \leq W$ and $\sum_{S' \in {\mathcal S} : S'_{[t]} = E'} \mu(S') = \mu(E')$, we conclude that 
$$C_{r,E'} \leq 4 \iota^{-2} \theta^{-2} U_r L W \mu(E').$$
Combining the above, we conclude that $\sum_{E \in {\mathcal E}, l \in \lbrace 0,\ldots,q-1 \rbrace} \big( \nabla f^{\mu,\theta}(X)_{E,l} - \nabla f^{\mu,\theta}(Y)_{E,l} \big)^2$ is at most
$$
4 \iota^{-2} \theta^{-2} U_r L W \sum_{r = 0}^{q-1} \sum_{E' \in {\mathcal E}} \mu(E') \big( \frac{X^r(E')}{\sqrt{\mu(E')}} - \frac{Y^r(E')}{\sqrt{\mu(E')}} \big)^2.
$$
Cancelling the $\mu(E')$ terms completes the proof.  $\qed$
\end{proof}
\ \indent We now complete the proof of Corollary\ \ref{convclaim2}, which follows directly by plugging in the above bounds for the relevant suprema into the convergence bounds of Claims\ \ref{convclaim1} and \ref{convclaim2a}.
\begin{proof}{Proof of Corollary\ \ref{convclaim2} :}
First, let us analyze the case $\beta_k = 0$ for all $k \geq 1$.   It follows from Claims\ \ref{whatisderiv1} and \ref{huberprops}, definitions, and the fact that $\sum_{E \in {\mathcal E}} \mu(E) = T$ that $$\sup_{X \in {\mathcal Q}^4(0,1)} ||\nabla f^{\mu,\theta}(X)||^2  \leq 4 \iota^{-2} L^2 T q.$$  It follows from Jensen's inequality that $$\sup_{X \in {\mathcal Q}^4(0,1)}\big|\big| \mathbb{E}\big[\nabla f^{\mu,\theta}(X) - \hat{G}^{\mu,1}(X) \big] \big|\big|
\leq \sqrt{\sup_{X \in {\mathcal Q}^4(0,1)} \mathbb{E}\big[ \big|\big|\nabla f^{\mu,\theta}(X) - 
\hat{G}^{\mu,1}(X)\big|\big|^2\big] }.$$
Combining the above with Claim\ \ref{boundnorm1} and some straightforward algebra, and plugging those bounds into the convergence guarantee of Claim\ \ref{convclaim1}, completes the proof in this case.  The case $\beta_k = \frac{k-1}{k+2}$ for all $k \geq 1$ follows by plugging the bounds of Claims \ref{boundnorm1} and \ref{smoothbound1} into the convergence guarantee of Claim\ \ref{convclaim2a}.  $\qed$
\end{proof}

\subsection{Proof of Claim\ \ref{algclaim1}}\label{APPalgclaim1sec}
We begin by making some preliminary observations regarding $\textrm{R}$.  Here we say that a given entry of $\Upsilon$ has been ``assigned a value" (``initialized with a value") if that entry is overwritten and assigned a value in ${\mathcal R}^q$ (initialized with a value in ${\mathcal R}^q$ due to the corresponding $k$ equalling $-1$ or $0$).  We refer to the last step of routine $\textrm{R},$ in which an entry of $\Upsilon$ is explicitly assigned a value, as the time at which the call to $\textrm{R}(E,k)$ ``computes $\Upsilon(E,k)$".  
\begin{observation}\label{algobs1}
\ \begin{enumerate}
\item For any given $E$ and $k \geq 1$, $\textrm{R}(E,k)$ only makes recursive calls to $\textrm{R}(E',k')$ for $k' < k$, and thus the recursive definition of $\textrm{R}(E,k)$ is well-defined and each call to $\textrm{R}(E,k)$ terminates in finite time. \label{algobs1a}
\item For any given $E$ and $k \geq 1$, during the first call to $\textrm{R}(E,k)$, before the function computes $\Upsilon(E,k)$ it will hold that ${\mathcal S}^{E,k-1}$ has been generated, and $\Upsilon(S'_{[t]},j)$ has been assigned a value (or initialized with a value) for all $S' \in {\mathcal S}^{E,k-1}, t \in \aleph^{k-1} \bigcap \bigcup_{i \in a^+(E)} {\mathcal T}_i(S')$ and $j \in \lbrace -1,\ldots,k-1 \rbrace$.  Also, $\Upsilon(E,j)$ will have been assigned a value for all $j \in \lbrace -1,\ldots,k-1 \rbrace$. \label{algobs1b}
\item During the first call to $\textrm{R}(E,k)$, when the function computes $\Upsilon(E,k)$, it only uses values of $\Upsilon$ which have already been assigned a value.\label{algobs1c}
\item For any given $E$ and $k \geq 1$, $\Upsilon(E,k)$ is assigned a value at the end of the execution of the first call to $\textrm{R}(E,k)$, and is not (re)assigned a value at any other time. \label{algobs1d}
\item For any given $E$ and $k \geq 1$, ${\mathcal S}^{E,k-1}$ is generated and stored at the start of the first call to $\textrm{R}(E,k)$, and never generated again.  \label{algobs1e}
\end{enumerate}
\end{observation}
\begin{proof}{Proof of Observation\ \ref{algobs1} :}
We proceed by induction on $k$, with base case $k = 1$.  (\ref{algobs1a}) follows from the fact that all recursive calls made in $\textrm{R}(E,k)$ are to $\textrm{R}(E',k-1)$ for different values of $E'$.  (\ref{algobs1b}) follows from the initialization of $\Upsilon(E',-1)$ and $\Upsilon(E',0)$, and fact that ${\mathcal S}^{E,k-1}$ is generated at the start of the first call to $\textrm{R}(E,k)$.  (\ref{algobs1c}) follows from (\ref{algobs1b}).  (\ref{algobs1d}) follows from the ``If $\Upsilon(E,k) = \emptyset$" statement at the start of $\textrm{R}(E,k)$.  (\ref{algobs1e}) follows from the same logic as (\ref{algobs1d}).  For the induction case, suppose the induction is true for all $j \in \lbrace 1,\ldots,k \rbrace$ for some $k \geq 1$.  We now prove the induction also holds for $k+1$.  (\ref{algobs1a}) follows for the same reason as in the base case.  To prove (\ref{algobs1b}), observe that for each $S' \in {\mathcal S}^{E,k}, t \in \aleph^{k} \bigcap \bigcup_{i \in a^+(E)} {\mathcal T}_i(S')$, before the function computes $\Upsilon(E,k+1)$, the for loop ensures that either $\Upsilon(S'_{[t]},k)$ has already been assigned a value (which must mean that $\textrm{R}(S'_{[t]},k)$ has completed an execution), or $\textrm{R}(S'_{[t]},k)$ is called.  Either way, before the function computes $\Upsilon(E,k+1)$, $\textrm{R}(S'_{[t]},k)$ has completed an execution.  The desired result then follows by applying the induction hypothesis to (\ref{algobs1b}) and (\ref{algobs1d}), along with the fact that ${\mathcal S}^{E,k}$ is generated at the start of the first call to $\textrm{R}(E,k+1)$.  An identical logic demonstrates that $\Upsilon(E,j)$ will have been assigned a value for all $j \in \lbrace -1,\ldots,k \rbrace.$  (\ref{algobs1c}) follows from the already proven (\ref{algobs1b}) and a straightforward inspection of the manner in which the value of $\Upsilon(E,k+1)$ is set.  (\ref{algobs1d}) and (\ref{algobs1e}) follow from the same logic as the base case.  Combining the above completes the proof.  $\qed$
\end{proof}
\begin{proof}{Proof of Claim\ \ref{algclaim1} :}
Let us proceed by induction, showing that any entry $\Upsilon(E,k)$ assigned a value must be assigned value ${\mathcal X}^k(E)$ as computed by \textrm{Algorithm} \ref{alg:meta-gd}.  Let us begin with the base case $k = 1$ (the cases $k = -1,0$ are trivial).  By Observation\ \ref{algobs1}, for any $E$, it suffices to consider the first time $\textrm{R}(E,1)$ is called.  During this call, due to the initialization of $\Upsilon(E',-1)$ and $\Upsilon(E',0)$ for all $E' \in {\mathcal E}$, no recursive calls are made, and $\Upsilon(E,1)$ is set equal to the $q$-dimensional vector
$$
\Pi_{{\mathcal Q}^0}\bigg( (1 + \beta_{0}) \Upsilon(E,0) - \beta_{0} \Upsilon(E,-1) + \alpha \hat{G}^{0}\big( (1 + \beta_{0}) \Upsilon^{0} - \beta_{0} \Upsilon^{-1}\big)_E \bigg).$$
It is easily verified that by construction ${\mathcal X}^1(E)$ has the same value.
\\\indent Now, suppose that for some $k \geq 1$, and all $j \leq k$ and $E \in {\mathcal E}$, any entry $\Upsilon(E,j)$ assigned a value must be assigned value ${\mathcal X}^j(E).$  Consider the first time $\textrm{R}(E,k+1)$ is called.  By Observation\ \ref{algobs1} and the induction hypothesis, at the end of the execution of $\textrm{R}(E,k+1)$, the value of $\Upsilon(E,k+1)$ will be set to
$$\Pi_{{\mathcal Q}^0}\bigg( (1 + \beta_{k}) \Upsilon(E,k) - \beta_{k} \Upsilon(E,k-1) + \alpha \hat{G}^{k}\big( (1 + \beta_{k}) \Upsilon^{k} - \beta_{k} \Upsilon^{k-1}\big)_E \bigg).$$
It is easily verified that by construction ${\mathcal X}^{k+1}(E)$ has the same value, completing the proof.  Further note that finite termination of $\textrm{R}(E,k)$, along with the fact that an entry of $\Upsilon$ assigned a value is never overwritten (i.e. the value is permanent), both follow from  Observation\ \ref{algobs1}.  Combining the above completes the proof.  $\qed$
\end{proof}

\subsection{Proof of Lemma\ \ref{complexity1b}}\label{APPcomplexity1bsec}
To avoid the need to discuss manipulations at the level of data structures, for $E \in {\mathcal E}$ and $S' \in \mathcal{S}$ s.t. $E \subseteq S'$, 
we assume the set of times $\bigcup_{i \in a^+(E)} \mathcal{T}_i(S')$ can be extracted in unit time given $S' \in {\mathcal S}.$

\begin{proof}{Proof of Lemma\ \ref{complexity1b} :}
To ensure our complexity analysis of $\textrm{R}$ is as tight as possible, note that we may equivalently define $\hat{G}^k(X)_{E,l}$ as follows (since all other terms vanish) :  
\begin{equation}\label{ghatkeeq}
\hat{G}^k(X)_{E,l} = 
Z^l(E) - 2 \iota^{-1} \sum_{i \in a^{l,+}(E)} a^l_{i}(E) \times \eta^{-1}_1 \sum_{S' \in {\mathcal S}^{E,k}}\phi'_{\theta}\bigg( \frac{T}{\eta_2} \sum_{t \in \aleph^k \bigcap {\mathcal T}_i(S')} \sum_{r=1}^{q-1} a^r_{i}({S'}_{[t]}) X^r({S'}_{[t]}) - b_i \bigg).
\end{equation}
To prevent any confusion about the runtime of manipulating $\aleph^{k-1} \bigcap \bigcup_{i \in a^+(E)} {\mathcal T}_i(S')$, we treat two cases separately : the case $\eta_2 = T$ (i.e. no subsampling of the sum), and the general case.  First, suppose $\eta_2 = T$.  In that case, $\aleph^{k-1} \bigcap \bigcup_{i \in a^+(E)} {\mathcal T}_i(S') = \bigcup_{i \in a^+(E)} {\mathcal T}_i(S')$.  In addition, $|\bigcup_{i \in a^+(E)} {\mathcal T}_i(S')| \leq U \times L \times q$ for each $S' \in {\mathcal S}^{E,k-1}$.  For $E \in {\mathcal E}$, let $\overline{Z}(E)$ denote the $q$-dimensional vector with component $l$ equal to $Z^l(E)$.  Then a call to $\textrm{R}(E,k)$ constitutes (in the worst case) :
\begin{itemize}
\item $O(\eta_1)$ calls to $\textrm{SIM}$ (at a total cost of $O(\eta_1 \times {\mathcal C})$);
\item $O(\eta_1)$ units of computational time to extract the set $\bigcup_{i \in a^+(E)} {\mathcal T}_i(S')$ for all $\eta_1$ of the $S' \in {\mathcal S}^{E,k-1}$;
\item the evaluation of at most $O(\eta_1 \times U \times L \times q)$ if statements (each costing one unit of time) and $O(\eta_1 \times U \times L \times q)$ recursive calls to $\textrm{R}(\cdot,k-1)$ (where these recursive calls will play a key role in the complexity analysis);
\item a single calculation to compute $\Upsilon(E,k)$, whose complexity we evaluate as follows.  There is a single projection onto ${\mathcal Q}^0$ (taking $O\big(q \log(q)\big)$ units of time); $O(q)$ units of computational time to query the values $\Upsilon(E,k-1), \Upsilon(E,k-2), \lbrace Z^r(E), r = 0,\ldots,q-1 \rbrace$ (which have already been computed); $O(q)$ multiplications and additions to compute $$(1 + \beta_{k-1}) \Upsilon(E,k-1) - \beta_{k-1} \Upsilon(E,k-2) + \alpha \overline{Z}(E).$$  Next, for each $S' \in {\mathcal S}^{E,k-1}$ and $i \in a^+(E)$, we pre-compute 
\begin{equation}\label{precomputeme1}
\sum_{t \in {\mathcal T}_i(S')} \sum_{r=1}^{q-1} a^r_{i}(S'_{[t]}) \big( (1 + \beta_{k-1}) \Upsilon(S'_{[t]},k-1)_r - \beta_{k-1} \Upsilon(S'_{[t]},k-2)_r \big) - b_i.
\end{equation}  For each such $S'$ and $i$, this will require querying $O(q \times |{\mathcal T}_i(S')|)$ values which will already be accessible (the $a_i(\cdot), \Upsilon(S'_{[t]},k-1)_r,\Upsilon(S'_{[t]},k-2)_r$); and performing $O(q \times |{\mathcal T}_i(S')|)$ additions and multiplications.  Thus the total time to precompute (\ref{precomputeme1}) for all $S' \in {\mathcal S}^{E,k-1}$ and $i \in a^+(E)$ is $O\big( \eta_1 \times |a^+(E)| \times q \times |{\mathcal T}_i(S')|\big)$, which (since $|a^+(E)| \leq L \times q$ and $|{\mathcal T}_i(S')| \leq U$) is itself $O\big( \eta_1 L q^2 U \big)$.  Next, we must evaluate $\phi'_{\theta}$ on each such sum, which (as the sums are precomputed) will take $O(\eta_1 U)$ time (since each individual evaluation of $\phi'_{\theta}$ takes $O(1)$ time).  For each $i \in a^+(E)$ we must then sum 
$\phi'_{\theta}\bigg(\sum_{t \in {\mathcal T}_i(S')} \sum_{r=1}^{q-1} a^r_{i}(S'_{[t]}) \big( (1 + \beta_{k-1}) \Upsilon(S'_{[t]},k-1)_r - \beta_{k-1} \Upsilon(S'_{[t]},k-2)_r \big) - b_i\bigg)$ over $S'$, which (as the relevant terms are precomputed) will (in total) take $O(\eta_1 L q)$ time.  For each $l \in \lbrace 0,\ldots,q-1 \rbrace$ we must then combine the relevant quantities above to compute the right-hand-side of (\ref{ghatkeeq}), which (as for each $l$ one must sum $|a^{l,+}(E)|$ pre-computed terms, with term $i$ multiplied by $a^l_i(E) \times \eta_1^{-1}$) requires (in total) time $O(q L)$.  As the above can then be combined with $O(q)$ additional additions and multiplications to compute $\Upsilon(E,k)$, this may be seen to lead a computational time (to compute $\Upsilon(E,k)$) of $O\big(\eta_1 U L q^2 \big).$
\end{itemize}
\ \indent Combining the above, we find that in the case $\eta_2 = T$, for an absolute constant $c$, $$\textrm{COMPLEXITY}(k) \leq c \eta_1 U L q^2 {\mathcal C} + c \eta_1 U L q \times \textrm{COMPLEXITY}(k-1).$$  It then follows from a straightforward induction that $$\textrm{COMPLEXITY}(k) \leq c \eta_1 U L q^2 {\mathcal C} \sum_{j=0}^{k-1} (c \eta_1 U L q)^j + (c \eta_1 U L q)^k,$$ 
and combining with some straightforward algebra completes the proof in this case.
\\\indent Next, consider the general case.  After making the $O(\eta_1)$ calls to $\textrm{SIM}$ (again at a total cost of $O(\eta_1 {\mathcal C})$), we compute the set $\aleph^{k-1} \bigcap \bigcup_{i \in a^+(E)} {\mathcal T}_i(S')$ as follows.  For each $t \in \aleph^{k-1}$, $i \in a^+(E)$, $S' \in {\mathcal S}^{E,k-1}$, and $r \in \lbrace 0,\ldots,q-1 \rbrace$, we query whether $a^r_i(S'_{[t]}) \neq 0$.  In total this takes $O(\eta_1 \eta_2 L q^2)$ units of computational time.  Next, we use the bound $|\aleph^{k-1} \bigcap \bigcup_{i \in a^+(E)} {\mathcal T}_i(S')| \leq |\aleph^{k-1}| = \eta_2$ to conclude that we must evaluate at most $O(\eta_1 \eta_2)$ if statements and make $O(\eta_1 \eta_2)$ recursive calls to $\textrm{R}(\cdot,k-1)$.  We then again make a single calculation to compute $\Upsilon(E,k)$, which can be implemented in computational time at most $O\big(\eta_1 \eta_2 L q^2\big)$ (by an argument very similar to that in the case $\eta_2 = T$, and the details of which we omit).  Combining the above with an argument very similar to that in the case $\eta_2 = T$ then completes the proof.  $\qed$
\end{proof}

\subsection{Proof of Claim\ \ref{NRMAIN3citenoweq1}}\label{NRMAIN3citenoweq1proofsec}
First, we show that $f^{\theta}(X)$ is close to $f(X)\ \forall\ X$.
\begin{claim}\label{thetaclose1}
For all $X \in {\mathcal Q}^2$, $\big| f^{\theta}(X) - f(X) \big| \leq \iota^{-1} V \theta$.
\end{claim}
\begin{proof}{Proof of Claim\ \ref{thetaclose1} :}
It follows from Claim\ \ref{huberprops} (see Section\ \ref{appconvergence1sec}), and the fact that $\phi_{\theta}(x) = x^+ = 0$ for $x \leq 0$, that that for any $S \in {\mathcal S}$ and $i \in \lbrace 1,\ldots,m \rbrace$,
$$
\big|\phi_{\theta}\big( \sum_{t=1}^T \sum_{r=1}^{q-1} a^r_i(S_{[t]}) X^r(S_{[t]}) - b_i \big) - \big( \sum_{t=1}^T \sum_{r=1}^{q-1} a^r_i(S_{[t]}) X^r(S_{[t]}) - b_i \big)^+\big| \leq \frac{1}{2} \theta I\big( \sum_{t=1}^T \sum_{r=1}^{q-1} a^r_i(S_{[t]}) X^r(S_{[t]}) > b_i \big).
$$
It follows that $\frac{1}{2} \iota \big| f(X) - f^{\theta}(X) \big|$ equals
$$\bigg| \sum_{S \in {\mathcal S}} \mu(S) \bigg( \sum_{i=1}^m \big( \sum_{t=1}^T \sum_{r=1}^{q-1} a^r_i(S_{[t]}) X^r(S_{[t]}) - b_i \big)^+ - \sum_{i=1}^m \phi_{\theta}\big( \sum_{t=1}^T \sum_{r=1}^{q-1} a^r_i(S_{[t]}) X^r(S_{[t]}) - b_i \big) \bigg) \bigg|,$$
itself at most 
$$\sum_{S \in {\mathcal S}} \mu(S)  \sum_{i=1}^m \bigg| \big(\sum_{t=1}^T \sum_{r=1}^{q-1} a^r_i(S_{[t]}) X^r(S_{[t]}) - b_i \big)^+ - \phi_{\theta}\big( \sum_{t=1}^T \sum_{r=1}^{q-1} a^r_i(S_{[t]}) X^r(S_{[t]}) - b_i \big)  \bigg|,$$
itself at most 
$$\frac{1}{2} \theta  \sum_{S \in {\mathcal S}} \mu(S)  \sum_{i=1}^m I\big( \sum_{t=1}^T \sum_{r=1}^{q-1} a^r_i(S_{[t]}) X^r(S_{[t]}) > b_i \big) \leq \frac{1}{2} \theta V  \sum_{S \in {\mathcal S}} \mu(S)\ \ =\ \ \frac{1}{2} \theta V.$$
Combining the above completes the proof.  $\qed$
\end{proof}
We now complete the proof of Claim\ \ref{NRMAIN3citenoweq1}.
\begin{proof}{Proof of Claim\ \ref{NRMAIN3citenoweq1} :}
First, we claim that $\big| \textrm{OPT}_{\textsc{pen}} - \textrm{OPT}_{\textsc{pen}^{\theta}} \big| \leq \iota^{-1} V \theta.$  Indeed, it follows from Claim\ \ref{thetaclose1} that $\textrm{OPT}_{\textsc{pen}}$ equals
$$f\big( X_{*,\textsc{pen}}\big)\ \geq\ f\big( X_{*,\textsc{pen}^{\theta}}\big)\ \geq\ f^{\theta}\big( X_{*,\textsc{pen}^{\theta}}\big) - \iota^{-1} V \theta\ =\ \textrm{OPT}_{\textsc{pen}^{\theta}} - \iota^{-1} V \theta.$$
As a nearly identical symmetric argument proves the other direction, this completes the proof.  Again applying Claim\ \ref{thetaclose1}, it thus holds that for $X \in {\mathcal Q}^2$, $\textrm{OPT}_{\textsc{pen}}  - f(X) \leq \big( \textrm{OPT}_{\textsc{pen}^{\theta}} + \iota^{-1} V \theta \big) - \big( f^{\theta}(X) - \iota^{-1}  V \theta \big).$  Simplifying completes the proof.  $\qed$
\end{proof}  

\subsection{Analysis of complexity of $\textrm{F}$}\label{compfsec}
We may implement $\textrm{F}(\textrm{A}_{\textsc{pen}})$ efficiently by maintaining $m$ counters, with counter $i$ containing (when the value of $\textrm{F}(\textrm{A}_{\textsc{pen}})^r(S_{[t]})$ must be set) the quantity $b_i - \sum_{u=1}^{t-1} \sum_{l=1}^{q-1} a^l_i(S_{[u]}) \textrm{F}(X)^l(S_{[u]}) - \sum_{l=1}^{r - 1} a^l_i(S_{[t]}) \textrm{F}(X)^l(S_{[t]})$.  Note that in any given time period, all necessary counter updates (across all $r \in \lbrace 1,\ldots,q-1 \rbrace$ and relevant $i \in \lbrace 1,\ldots,m \rbrace$) may be implemented in time $O(L \times q)$ (as only $L$ of the counters will need to be updated for each $r \in \lbrace 1,\ldots,q-1 \rbrace$), and each update requires $O(1)$ additions and multiplications).  It follows from the definition of $\textrm{F}$ that (in addition to the time to call $\textrm{A}_{\textsc{pen}}$) implementing 
$\textrm{F}(\textrm{A}_{\textsc{pen}})$ will thus require an additional $O(L \times q)$ time units of computation. 

\subsection{Proof of Lemma\ \ref{nrmroundlemma1}}\label{APPnrmroundlemma1sec}
\begin{proof}{Proof of Lemma\ \ref{nrmroundlemma1} :}
We will prove the desired inequality for each individual $S \in {\mathcal S}$, i.e. we will prove that for all $S \in {\mathcal S}$, 
$\mathbb{E}\big[ \sum_{i=1}^m \big( \sum_{t=1}^T \sum_{r=1}^{q-1} a^r_i(S_{[t]}) \textrm{H}(X)^r(S_{[t]}) - b_i \big)^+ \big]$ is at most  
$$\sum_{i=1}^m \big( \sum_{t=1}^T \sum_{r=1}^{q-1} a^r_i(S_{[t]}) X^r(S_{[t]}) - b_i \big)^+ + \sqrt{\frac{\pi}{2}} \min\big( \sqrt{L q m} \sqrt{T} , V \sqrt{T} + m T \exp( - \frac{\nu^2}{2} T ) \big).$$  
We break the proof into two parts (one for each term of the min), and begin by proving that 
$$\mathbb{E}\big[ \sum_{i=1}^m \big( \sum_{t=1}^T \sum_{r=1}^{q-1} a^r_i(S_{[t]}) \textrm{H}(X)^r(S_{[t]}) - b_i \big)^+ \big] - 
\sum_{i=1}^m \big( \sum_{t=1}^T \sum_{r=1}^{q-1} a^r_i(S_{[t]}) X^r(S_{[t]}) - b_i \big)^+ \leq \sqrt{L q m} \sqrt{T}.$$  
At a high level, the proof proceeds by a standard analysis of the relevant randomized rounding ($\textrm{H}$) to derive a $\Theta(\sqrt{T})$ error using Hoeffding's inequality.  This yields an error bound for each resource of the form $\sqrt{T_i}$ for some appropriate $T_i$, and we then combine with an extra step at the end in which we solve a simple optimization problem to derive a tighter bound for $\sum_{i=1}^m \sqrt{T_i}$ (the overall error since the penalty term in our objective sums over all resources) using the assumed L-sparsity.
\\\indent We first simplify the term we wish to bound a bit using the triangle inequality and continuity of relevant functions.  By the triangle inequality, the definition of ${\mathcal T}_i(S)$, and Lipschitz continuity of $g(z) \stackrel{\Delta}{=} (z - b_i)^+$,
$$
(a) : \mathbb{E}\big[ \sum_{i=1}^m \big( \sum_{t=1}^T \sum_{r=1}^{q-1} a^r_i(S_{[t]}) \textrm{H}(X)^r(S_{[t]}) - b_i \big)^+ \big] - 
\sum_{i=1}^m \big( \sum_{t=1}^T \sum_{r=1}^{q-1} a^r_i(S_{[t]}) X^r(S_{[t]}) - b_i \big)^+$$
is at most
$$(b) : \sum_{i=1}^m \mathbb{E}\big[ \big| \sum_{t \in {\mathcal T}_i(S) } \sum_{r=1}^{q-1} a^r_i(S_{[t]}) \textrm{H}(X)^r(S_{[t]}) - \sum_{t \in {\mathcal T}_i(S) } \sum_{r=1}^{q-1} a^r_i(S_{[t]}) X^r(S_{[t]})\big| \big].
$$
Let us fix $i \in \lbrace 1,\ldots,m \rbrace$, and use Hoeffding's inequality to bound
$$(c) : 
\mathbb{E}\big[ \big| \sum_{t \in {\mathcal T}_i(S) } \sum_{r=1}^{q-1} a^r_i(S_{[t]}) \textrm{H}(X)^r(S_{[t]}) - \sum_{t \in {\mathcal T}_i(S) } \sum_{r=1}^{q-1} a^r_i(S_{[t]}) X^r(S_{[t]})\big| \big].
$$
Observe that by construction $\mathbb{E}\big[ \sum_{r=1}^{q-1} a^r_i(S_{[t]}) \textrm{H}(X)^r(S_{[t]}) \big] = \sum_{r=1}^{q-1} a^r_i(S_{[t]}) X^r(S_{[t]})$ and $\sum_{r=1}^{q-1} a^r_i(S_{[t]}) \textrm{H}(X)^r(S_{[t]}) \in [0,1]$.  Thus as the rounding is independent (across different $S_{[t]}$), we may directly apply Hoeffding's inequality to conclude that for all $x > 0$,
$$\mathbb{P}\big( \big| \sum_{t \in {\mathcal T}_i(S) } \sum_{r=1}^{q-1} a^r_i(S_{[t]}) \textrm{H}(X)^r(S_{[t]}) - \sum_{t \in {\mathcal T}_i(S) } \sum_{r=1}^{q-1} a^r_i(S_{[t]}) X^r(S_{[t]})\big| > x \big) \leq 2 \exp(- \frac{ 2 x^2 }{|{\mathcal T}_i(S)|}).
$$
Combining with the tail-integral form of the expectation of a non-negative random variable, along with known results for Gaussian integrals, we conclude that (c) is at most $\sqrt{\frac{\pi}{2}} \sqrt{ |{\mathcal T}_i(S)|}$.  Combining with (b), we conclude that (a) is at most $\sqrt{\frac{\pi}{2}} \sum_{i=1}^m \sqrt{ |{\mathcal T}_i(S)| }$.  As (by computing the sum in two different ways and using the definitions of $L$ and $q$) $\sum_{i=1}^m {\mathcal T}_i(S) \leq L q T$, we further conclude that (a) is at most the value of the following optimization problem : 
$$ \max \sqrt{\frac{\pi}{2}} \sum_{i=1}^m \sqrt{ T_i }\ \ \ \textrm{s.t.}\ \ \ \sum_{i=1}^m T_i \leq L q T\ \ ;\ \ T_i \geq 0\ \forall i.$$
Due to the concavity of the square root function, it is straightforward to show that at optimality $T_i = \frac{L q T}{m}$ for all $i$, and the optimal value is $\sqrt{\frac{\pi}{2}} \times m \times \sqrt{\frac{L q T}{m}} = \sqrt{\frac{\pi}{2}} \times \sqrt{m L q T}$.  Combining the above completes the proof of the desired inequality.
\\\\We now prove the bound corresponding to the other term in the min, i.e. we prove that for all $S \in {\mathcal S}$, 
$\mathbb{E}\big[ \sum_{i=1}^m \big( \sum_{t=1}^T \sum_{r=1}^{q-1} a^r_i(S_{[t]}) \textrm{H}(X)^r(S_{[t]}) - b_i \big)^+ \big]$ is at most  
$$\sum_{i=1}^m \big( \sum_{t=1}^T \sum_{r=1}^{q-1} a^r_i(S_{[t]}) X^r(S_{[t]}) - b_i \big)^+ + \sqrt{\frac{\pi}{2}} \big( V \sqrt{T} + m T \exp( - \frac{\nu^2}{2} T ) \big).$$
We proceed by using a similar analysis to that used for the first term of the min to get a $\Theta(\sqrt{T})$ error bound for each of the (at most) $V$ resources which are at least half-saturated on that $S$.  We then use a different analysis for the remaining $m - V$ resources.  In particular, we use Hoeffding's inequality to argue that each of those $m - V$ (which we simply bound by $m$) other inequalities is exponentially unlikely to be saturated (where the exponent contains $\nu$).  If they are not saturated, the result is zero error.  If they are saturated (which for each inequality is exponentially unlikely), we bound the worst-case error by its largest possible value $T$.  Combining the above will complete the proof.
\\\indent Let $V(S,X) \stackrel{\Delta}{=} \lbrace i : \sum_{t=1}^T \sum_{r=1}^{q-1} a^r_i(S_{[t]}) X^r(S_{[t]}) \geq \frac{b_i}{2} \rbrace$.  By the definition of $V$, $|V(S,X)| \leq V$.  Note that (a) equals
\begin{eqnarray*}
\ &\ &\ (b')\ \ \mathbb{E}\big[ \sum_{i \in V(S,X)} \big( \sum_{t=1}^T \sum_{r=1}^{q-1} a^r_i(S_{[t]}) \textrm{H}(X)^r(S_{[t]}) - b_i \big)^+ \big] - 
\sum_{i \in V(S,X)} \big( \sum_{t=1}^T \sum_{r=1}^{q-1} a^r_i(S_{[t]}) X^r(S_{[t]}) - b_i \big)^+
\\&\ &\ (b'')\ \ + \mathbb{E}\big[ \sum_{i  \notin V(S,X)} \big( \sum_{t=1}^T \sum_{r=1}^{q-1} a^r_i(S_{[t]}) \textrm{H}(X)^r(S_{[t]}) - b_i \big)^+ \big].
\end{eqnarray*}
By essentially the same analysis used in the proof of the first inequality, we find that (b') is at most $V \sqrt{\frac{\pi}{2}} \sqrt{T}$.  We now analyze (b'').  Let us fix $i \notin V(S,X)$.  Observe that
$$\big( \sum_{t=1}^T \sum_{r=1}^{q-1} a^r_i(S_{[t]}) \textrm{H}(X)^r(S_{[t]}) - b_i \big)^+ \leq T \times I\big( \sum_{t=1}^T \sum_{r=1}^{q-1} a^r_i(S_{[t]}) \textrm{H}(X)^r(S_{[t]}) > b_i \big).$$
By the definition of $V(S,X)$ it holds that $$E[\sum_{t=1}^T \sum_{r=1}^{q-1} a^r_i(S_{[t]}) \textrm{H}(X)^r(S_{[t]})] = \sum_{t=1}^T \sum_{r=1}^{q-1} a^r_i(S_{[t]}) X^r(S_{[t]}) \leq \frac{b_i}{2}.$$  We may thus apply Hoeffding's inequality in essentially the same manner used in the proof of the first inequality to conclude that 
\begin{eqnarray*}
\ &\ &\ \mathbb{P}\big( \sum_{t=1}^T \sum_{r=1}^{q-1} a^r_i(S_{[t]}) \textrm{H}(X)^r(S_{[t]}) > b_i \big)
\\&\leq& \mathbb{P}\bigg( \sum_{t=1}^T \sum_{r=1}^{q-1} a^r_i(S_{[t]}) \textrm{H}(X)^r(S_{[t]}) - E\big[\sum_{t=1}^T \sum_{r=1}^{q-1} a^r_i(S_{[t]}) \textrm{H}(X)^r(S_{[t]})\big] > \frac{1}{2} b_i \bigg)
\\&\leq& \exp\big( - \frac{ 2 (\frac{1}{2} b_i)^2}{T} \big)\ \ =\ \ \exp\big( - \frac{b_i^2}{2 T} \big)\ \ \leq\ \ \exp( - \frac{\nu^2}{2} T ).
\end{eqnarray*}
Combining the above with the linearity of expectation completes the proof.
\end{proof}

\subsection{Proof of Lemma\ \ref{nrmroundlemma2}}\label{APPnrmroundlemma2sec}
\begin{proof}{Proof of Lemma\ \ref{nrmroundlemma2} :}
We begin by proving the first part of the lemma.  Let $\Gamma$ denote the set of $(r,t)$ s.t. $\textrm{F}(X)^r(S_{[t]})$ is non-integer.  Then as $\textrm{F}(X)^r(S_{[t]}) \leq X^r(S_{[t]})$, for $(r,t) \in \Gamma$ it must hold that $\textrm{F}(X)^r(S_{[t]}) \in \big(0 , X^r(S_{[t]}) \big)$, and thus (by construction of $\textrm{F}$) 
$$\min_{i \in a^{r,+}(S_{[t]})} \frac{b_i - \sum_{u=1}^{t-1} \sum_{l=1}^{q-1} a^l_i(S_{[u]}) \textrm{F}(X)^l(S_{[u]}) - \sum_{l=1}^{r - 1} a^l_i(S_{[t]}) \textrm{F}(X)^l(S)}{a^r_i(S_{[t]})} \in \big(0, X^r(S_{[t]}) \big).$$
For $(r,t) \in \Gamma$, Let ${\mathcal I}(r,t)$ denote the set of $i$ at which $\min_{i \in a^{r,+}(S_{[t]})} \frac{b_i - \sum_{u=1}^{t-1} \sum_{l=1}^{q-1} a^l_i(S_{[u]}) \textrm{F}(X)^l(S_{[u]}) - \sum_{l=1}^{r - 1} a^l_i(S) \textrm{F}(X)^l(S_{[t]})}{a^r_i(S_{[t]})}$ attains its (strictly positive) value (i.e. the set of minimizers), and note that one must have $|{\mathcal I}(r,t)| \geq 1$ for all $(r,t) \in \Gamma$.  It follows from a straightforward contradiction that for all $(r,t) \in \Gamma$ and $i \in {\mathcal I}(r,t)$, 
{\small
$$\sum_{u=1}^{t-1} \sum_{l=1}^{q-1} a^l_i(S_{[u]}) \textrm{F}(X)^l(S_{[u]}) + \sum_{l=1}^{r - 1} a^l_i(S) \textrm{F}(X)^l(S_{[t]}) < b_i\ \textrm{and}\ \sum_{u=1}^{t-1} \sum_{l=1}^{q-1} a^l_i(S_{[u]}) \textrm{F}(X)^l(S_{[u]}) + \sum_{l=1}^{r} a^l_i(S_{[t]}) \textrm{F}(X)^l(S) = b_i.$$}
Thus any given $i \in \lbrace 1,\ldots,m \rbrace$ appears in ${\mathcal I}(r,t)$ for at most one $(r,t)$.  Furthermore, $\sum_{u=1}^{t-1} \sum_{l=1}^{q-1} a^l_i(S_{[u]}) \textrm{F}(X)^l(S_{[u]}) + \sum_{l=1}^{r} a^l_i(S_{[t]}) \textrm{F}(X)^l(S_{[t]}) = b_i$ implies that inequality is saturated.  Thus there can be at most $V$ such indices, and hence at most $V$ such $(r,t)$ pairs. 
\\\indent We now prove the second part of the lemma.  For fixed $S \in {\mathcal S}$ and $X \in {\mathcal Q}^2$, let ${\mathcal V}$ denote $\lbrace i : \sum_{t=1}^T \sum_{r=1}^{q-1} a^r_i(S_{[t]}) X^r(S_{[t]}) \geq \frac{1}{2} b_i \rbrace$.  Note that $$\sum_{i \in {\mathcal V}} \sum_{t=1}^T \sum_{r=1}^{q-1} a^r_i(S_{[t]}) X^r(S_{[t]}) \geq \frac{1}{2} \sum_{i \in {\mathcal V}} b_i,$$ which itself implies that 
$$\sum_{i=1}^m \sum_{t=1}^T \sum_{r=1}^{q-1} a^r_i(S_{[t]}) X^r(S_{[t]}) \geq \frac{1}{2} T |{\mathcal V}| \nu.$$  However, 
\begin{eqnarray*}
\sum_{i=1}^m \sum_{t=1}^T \sum_{r=1}^{q-1} a^r_i(S_{[t]}) X^r(S_{[t]}) &=& \sum_{t=1}^T \sum_{r=1}^{q-1} X^r(S_{[t]}) \sum_{i=1}^m  a^r_i(S_{[t]}) 
\\&\leq& L \sum_{t=1}^T \sum_{r=1}^{q-1} X^r(S_{[t]})\ \textrm{since}\ \sum_{i=1}^m  a^r_i(S_{[t]}) \leq L
\\&\leq& L T\ \textrm{since}\ \sum_{r=1}^{q-1} X^r(S_{[t]}) \leq 1.
\end{eqnarray*}
We conclude that $L T \geq \frac{1}{2} T |{\mathcal V}| \nu$.  Dividing both sides by $T \nu$ implies $|{\mathcal V}| \leq \frac{2 L}{\nu}$.  As the argument holds for general $S \in {\mathcal S}$, we conclude that we may take $V \leq \frac{2 L}{\nu}$, completing the proof.  $\qed$
\end{proof}

\subsection{Proof of Theorem\ \ref{NRMAINCOMBOPT} for \textsc{mwmlp}}\label{MWMLPproofsec}
\begin{proof}{Proof of Theorem\ \ref{NRMAINCOMBOPT} for \textsc{mwmlp} :}
As explained in Section\ \ref{applicationssec}, $\textsc{mwmlp}$ can be put in the $\textsc{lp}$ framework, with $q = 2, T = \lfloor \frac{1}{2} \Delta n \rfloor, m = n;$ and one can take $L = 2, U = \Delta, U_r = \Delta, W = 2 \Delta, V = n, \iota = 1$.  The proof then follows directly from Theorem\ \ref{NRMAIN2}, applied with $\epsilon' = \frac{2}{\Delta} \epsilon.  \qed$
\end{proof}

\end{document}